\documentclass[a4paper,11pt]{amsart}
\usepackage[latin1]{inputenc}
\usepackage{amsmath,amsthm,amsfonts,amscd,amssymb,eucal,latexsym,mathrsfs}
\usepackage[all]{xy}
\usepackage{latexsym, amssymb, amscd, mathrsfs, graphics, fullpage,graphicx}

\newtheorem{theorem}{Theorem}
\newtheorem{corollary}[theorem]{Corollary}
\newtheorem{lemma}[theorem]{Lemma}
\newtheorem{proposition}[theorem]{Proposition}
\newtheorem{claim}[theorem]{Claim}
\theoremstyle{definition}
\newtheorem{remark}[theorem]{Remark}
\newtheorem{remarks}[theorem]{Remarks}

\newtheorem{example}[theorem]{Example}

\newtheorem{question}[theorem]{Question}

\newtheorem{definition}[theorem]{Definition}

 \DeclareMathOperator{\Aut}{Aut}

\DeclareMathOperator{\supp}{ supp}

\newcommand{\reg}{r}

\def\min{\mathop{\mathrm{min}}\nolimits}
\def\max{\mathop{\mathrm{max}}\nolimits}

\newcommand{\CI}{{\mathrm{\bf C}}}\newcommand{\RI}{{\mathrm{\bf R}}}\newcommand{\QI}{{\mathrm{\bf Q}}}
\newcommand{\ZI}{{\mathrm{\bf Z}}}\newcommand{\NI}{{\mathrm{\bf N}}}

\newcommand{\Id}{\mathrm{Id}}\newcommand{\cA}{{\mathcal A}}\newcommand{\cH}{{\mathcal H}}
\newcommand{\cD}{{\mathcal D}}

\newcommand{\cF}{{\mathcal F}}

\newcommand{\id}{{\rm Id}}

\newcommand{\G}{\Gamma}

\newcommand{\TI}{{\bf T}}

\newcommand{\GL}{\mathrm {GL}}
\newcommand{\SL}{\mathrm {SL}}\newcommand{\PSL}{\mathrm {PSL}}

\newcommand{\FI}{{\bf{F}}}

\newcommand{\tr}{\mathop{\mathrm{Tr}}}

\newcommand{\del}{\partial}

\newcommand{\conv}{\mathrm{Conv}}
\newcommand{\Con}{\mathrm{Con}}

\newcommand{\sq}{${7\over 4}$}
\newcommand{\td}{$1^+$}
\newcommand{\sqm}{{7\over 4}}
\newcommand{\tdm}{{1^+}}

\date{\today}

\title{Intermediate rank and property RD}

\author{Sylvain Barr\'e}
\address{\hskip-\parindent
Sylvain Barr\'e, Universit\'e de Bretagne Sud, BP 573, Centre Yves Coppens, Campus de Tohannic, 56017 Vannes, France}
\email{sylvain.barre@univ.ubs.fr}
\author{Mika\"el Pichot}
\address{\hskip-\parindent
Mika\"el Pichot, Max-Planck-Institut f\"ur Mathematik, Vivatsgasse 7, 53111 Bonn, Germany}
\email{pichot@ihes.fr}

\begin{document}

\begin{abstract}
We introduce concepts of intermediate  rank for countable groups that ``interpolate" between consecutive values of the classical  (integer-valued) rank.    
Various classes of groups are proved to have intermediate rank behaviors. We are especially interested in  interpolation between rank 1 and rank 2. 

For instance, we construct groups ``of rank \sq".  
Our setting is essentially that of non positively curved spaces, where concepts of intermediate rank include  polynomial branching, local rank, and mesoscopic rank. 

The resulting framework has interesting connections to operator algebras.  We prove property RD  in many cases where intermediate rank occurs. This gives a new family of groups satisfying the Baum-Connes conjecture.  We  prove that the reduced $C^*$-algebras of groups of rank \sq\ have stable  rank 1.
\end{abstract}

\maketitle

The paper is organized along the following thematic lines.
\begin{itemize}
\item[A)] Rank interpolation from the viewpoint of property RD; 
\item[B)] Triangle polyhedra and the classical rank;
\item[C)] Polynomial and exponential branching,  branchings and property RD; 
\item[D)] Local rank, rank \sq, existence and classification results; 
\item[E)] Triangle polyhedra and property RD;  
\item[F)]  Applications to the Baum-Connes conjecture; 
\item[G)] $C^*$-algebraic  rank, stable rank, real rank; 
\item[H)] Mesoscopic rank. Mixed local rank.  
\end{itemize}



\section{Introduction and statement of the results}


\subsection{Definition of property RD} Let $\G$ be a countable group endowed with a length  $\ell$. One says that $\G$ has Property RD with respect to $\ell$ if there is a polynomial $P$ such that for any $r\in \RI_+$ and  $f,g\in \CI\G$  with $\supp(f)\subset B_r$ one has 
\[
\|f*g\|_2\leq P(r)\|f\|_2\|g\|_2
\]  
where $B_r=\{x\in \G,~\ell(x)\leq r\}$ is the ball of radius $r$ in $\G$. For example groups of polynomial growth have property RD as the number of decompositions $xy=z$ for fixed $z\in \G$ and $x,y\in \G$ with $\ell(x)\leq r$, is polynomial in $r$. 
 Property RD was introduced by Jolissaint in \cite{Jol-def} after the work of Haagerup 
 \cite{Haa} on reduced $C^*$-algebras of free groups. 
 
In the case of amenable groups  property RD implies polynomial growth \cite{Jol-def,Con,Val-bc}. 
 Thus  $\SL_3(\ZI)$, for example, does not have property RD  because it contains amenable subgroups which are not of polynomial growth  (see \cite{Jol-def}, this is the only  obstruction to property RD known so far).  
 
 Here are two fundamental examples of groups with property RD: 

\begin{enumerate}
\item  free groups on finitely many generators have property RD (with respect to the usual word length), as was proved by Haagerup in \cite{Haa},  
\item groups acting freely isometrically on Bruhat-Tits buildings of type $\tilde A_2$ (also called triangle buildings)  have property RD with respect to the length induced from the 1-skeleton, as was proved by Ramagge, Robertson and Steger in \cite{RRS}. 
 \end{enumerate}
Haagerup's result was generalized by Jolissaint  \cite{Jol-def} and de la Harpe \cite{Harpe-rd} to every hyperbolic groups in the sense of Gromov. 
The result of Ramagge, Robertson and Steger provided the first occurrences of property RD in ``higher rank"  situations. This was extended in \cite{Laf-rd}, where Lafforgue proves that all cocompact lattices in $\SL_3(\RI)$ (and $\SL_3(\CI)$)   have property RD. 

We refer to Section 2 for details and further recent developments.
 We are interested here in interpolation between  (1) and (2).

\subsection{Triangle polyhedra}
 
The \emph{rank} of a non positively curved metric space $X$ is an asymptotic invariant of $X$  usually defined as the dimension of maximal flats in $X$. A \emph{flat} in $X$ is the image of an isometric embedding of an Euclidean space $\RI^k$, $k\geq 2$.  In \cite[page 127]{Gro-asym}  seven definitions of the rank  are discussed.  As mentioned there, they   express the idea that $X$  behaves ``hyperbolically above the rank".

  
 \begin{definition}\label{def1}
We call  \emph{triangle polyhedron} a non positively curved (i.e.\ CAT(0))  
 simplicial complex $X$ of dimension 2 without boundary and whose faces are {\bf equilateral} triangles. 
 A countable group is called a \emph{triangle group}  if it admits a proper and isometric action on a triangle polyhedron. Proper means that  stabilizers are uniformly finite. 
  \end{definition}

The  role of triangle polyhedra below is to allow rank interpolation within a tractable geometrical framework. The word triangle refers to the Coxeter diagram of flats and follows the classical terminology for Tits buildings \cite{Ronan}. All flats  in $X$ are isometric to Euclidean planes $\RI^2$ tessellated by equilateral triangles (i.e.\ they are of type $\tilde A_2$).  Note that in the literature the terminology `triangle groups' can also refer to a class (different from the above) of reflection groups of the plane (Euclidean or hyperbolic, see e.g.\ \cite[V.38]{Harpe}).

Examples of triangle groups include $\tilde A_2$-groups, i.e.\ groups which act freely and simply transitively on triangle buildings  (see   \cite{Cart} and \cite{Ronan} for a general reference on Tits buildings) as well as many hyperbolic groups (in particular all triangle groups that satisfy the `girth $>6$' local condition \cite{Gro-asym}).  Note that by the \emph{no flat criterion} (see \cite[page 176]{Gro-asym}) the classical integer-valued rank (defined above)  detects precisely hyperbolicity among triangle polyhedra. 
For clarity we adopt the following 
 convention concerning the integer-valued rank. 
 
 A triangle polyhedron is said to have 
 \begin{itemize}
 \item \emph{rank 1} if it is hyperbolic, 
 \item \emph{rank 2} if it is symmetric, i.e., if it is a triangle building (\cite{BB,rg2}).
 \end{itemize}
For other triangle polyhedra, the integer-valued rank is undefined. Other concepts of rank shall be substituted to it.

Polyhedra in Definition \ref{def1} are not assumed to be locally finite a priori. Note that  every countable group admits a triangle presentation (by adding generators to a given presentation that split relations  into length 3 relations) but this presentation does not define a triangle polyhedron in general. 
Cohomological arguments implies that $\SL_3(\ZI)$ is not  a triangle group (see also Theorem \ref{th1}).   

\subsection{Intermediate rank as restricted branching growth} As we will see intermediate rank behaviours can be exhibited at the  microscopic,  mesoscopic, and macroscopic---or rather asymptotic---scale. We  first discuss the latter along with the notion of \emph{polynomial branching} in triangle polyhedra.

\begin{definition}\label{def2}
A triangle polyhedron $X$ is said to have \emph{polynomial branching} if there exists a polynomial $P$ such that for any simplicial geodesic segment $\gamma$ in $X$, the number of flat equilateral triangles in $X$ with base $\gamma$ is bounded by $P(r)$, where $r$ is the length of $\gamma$. One says that a triangle group has polynomial branching  if it admits a proper and isometric action on a triangle polyhedron of polynomial branching.
\end{definition}

In other words we restrict the \emph{branching of flats} in $X$ to be polynomial. This essentially captures spaces whose rank is ``not too far from to 1".   
 For instance triangle polyhedra which are hyperbolic or which have  isolated flats are of polynomial branching.  Triangle buildings, for which the branching of flats is exponential,  are not.

 In a similar way we define   \emph{subexponential branching}  by replacing the above polynomial $P$ by some given subexponential function.  Triangle polyhedra which are not of subexponential branching are said to be of \emph{exponential branching}. These notions  are asymptotic in nature and, as  in \cite{Gro-asym} (see page 127), can be detected at infinity.

Definition 2 is generalized to arbitrary countable groups endowed with a length (e.g.\ finitely generated groups with the word length) in Section \ref{reminder}. While being in the non amenable setting this generalization  is strongly reminiscent of the theory of \emph{growth of groups} (for which  we refer to \cite{Harpe} and references therein). It  relies on tools that arise from the study of property RD,  allowing  to merely retain the polynomial growth of `flats' rather than sharp flatness.  Subexponential and exponential branching can be defined analogously. The proof that  geodesic polynomial  branching coincide with the above Definition \ref{def2} for triangle groups is given in Section \ref{rd-triangles} (Proposition \ref{polrk-same}). In the amenable case the branching growth captures the classical notion of growth of groups (see Section \ref{sec22}). 
Regarding property RD the following holds. 

\begin{proposition}\label{prop3}
Let $\G$ be a countable group endowed with a length $\ell$. If $\G$ has polynomial branching with respect to $\ell$ then it has  property RD with respect to $\ell$.
\end{proposition}

Examples of groups with polynomial branching are relatively hyperbolic groups with respect to a finite family of groups of polynomial growth. In this case property RD was already known by a theorem of Chatterji and Ruane \cite{Chat-Ruane}.
Subexponential branching does not imply property RD in general, but it does imply some useful subexponential variations of it (see the end of Section \ref{sec22}).

\subsection{Polyhedra of rank \sq} We now turn to constructions of  triangle groups of intermediate rank. The groups exhibited in this Section  will be called \emph{groups of rank \sq}. As we will  explain, the rank conditions that prevail in these constructions are \emph{local}. Nonetheless,   large rank tends to propagate to the asymptotic level and groups of rank \sq, at least for some of them  (see Theorem \ref{th4}),  have exponential branching in the sense of the previous subsection (yet their rank is  strictly lower than that of triangle buildings).

Recall that the link of a CAT(0) complex $X$ at a point $A$ is the set of directions  at $A$ in $X$, endowed with the angular metric (see e.g.\ \cite[p. 103]{BH}). 
The tangent cone $\Con_A X$ of $X$ at $A$ is  the CAT(0) cone over the link of $A$ endowed with the angular metric. This is a CAT(0) space as well (\cite[p. 190]{BH}).
In this paper links are always assumed to be connected.

For a (geometrically finite) CAT(0) complex $X$ the \emph{local rank} at a point $A$    represents the 
\emph{proportion of flats in the tangent cone $\Con_A X$ of $X$ at $A$}. The notion of local rank \sq\ defined below is mediating between local rank 2 (the highest value) and a local version of the isolated flat property called local rank  \td, which indicates low local rank and that we introduce first. In \cite{cras}, the first author introduced a complex which satisfies this property (this is further discussed in Section \ref{class-rg74}).
  
Observe that in a triangle polyhedron links at vertices can be a priori  any graph with girth 6 (this is equivalent to the CAT(0) condition, edges have length $\pi/3$) so that  the proportion of tangent flats at a vertex corresponds to the proportion of 6-cycles in its link.


We say that  the local rank of a  2-dimensional CAT(0) complex $X$ to be $\leq$ \td\  if the following condition is satisfied.

  \begin{itemize}
\item \emph{Local rank $\leq \tdm$.} For every vertex $A$ of $X$, every segment (not necessarily simplicial) of length $\pi$ in the link $L$ at $A$ is included in at most one cycle  of length $2\pi$ in $L$.
\end{itemize}
The complex $X$ is said to have \emph{local rank \td}  if it has rank $\leq \tdm$ and if the following  condition is satisfied: 
\begin{itemize}
\item \emph{Local rank $\geq \tdm$.} For every vertex $A$ of $X$,   every edge in the link $L$ at $A$, as well as every pair of thick vertices at distance no greater than $\pi$,   is included in at least one cycle of length $2\pi$ of $L$ (where we say that a vertex is  thick if its valency  is at least 3).  
\end{itemize}
The condition local rank $\leq \tdm$ is a localization of the isolated flat property.  It is not hard to see that triangle polyhedra of rank $\leq$ \td\ have isolated flats. In particular they have polynomial branching.

 In the case of triangle buildings links corresponds to \emph{projective planes}.  They have local rank 2 in the sense that an incidence graph of a projective plane  is a \emph{spherical building}: compare \cite{these,rg2} and see also \cite{BB}  for a semi-local definition of rank 2. A 2-dimensional CAT(0) complex $X$ is said to have local rank 2 if the following condition, which characterizes spherical buildings, is satisfied (see \cite{rg2}).

  \begin{itemize}
\item \emph{Local rank 2.} For every vertex $A$ of $X$, every segment (not necessarily simplicial) of length  $\pi$ in  the link $L$ of $A$,  is included in at least one cycle of length $2\pi$. 
\end{itemize}

\noindent For triangle polyhedra, well-known examples of link of local rank 2 (i.e., of spherical buildings) include the incidence graph of the Fano plane (see  Figure 1.3 in \cite{Ronan}).  

The only link which is both of rank $\leq \tdm$ and of rank 2 is the circle of length $2\pi$ (cf.\ \cite{rg2}). Thus the only 2-dimensional CAT(0) complex of local rank $\leq \tdm$ and local rank 2 is $\RI^2$. Note that this is a (thin, reducible) building, which indeed has low rank but still,  large rank relatively to itself.   
  
  \bigskip
  
 We say  that a triangle polyhedron has \emph{local rank \sq}, or merely \emph{rank \sq} when no confusion can arise,   if its links at each  vertex are isomorphic to the following graph, henceforth denoted $L_{\sqm}$.

\begin{figure}[htbp]
\centerline{\includegraphics[width=4.5cm]{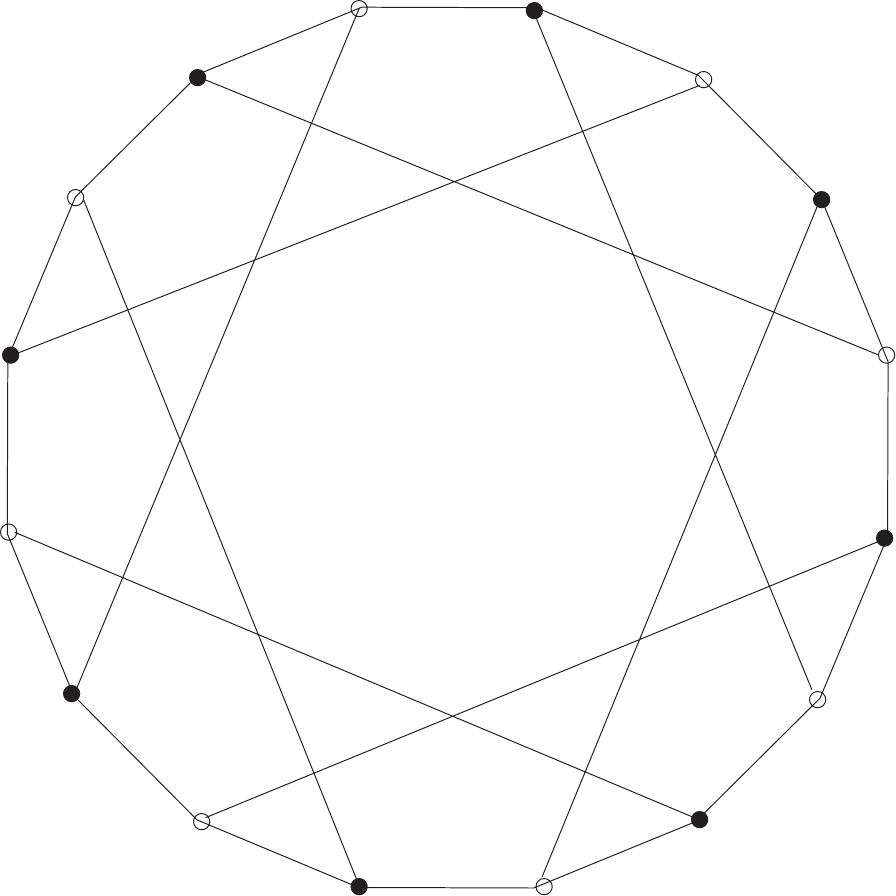}}
\caption{Rank \sq\ for triangle polyhedra}\label{fig1}
\end{figure}


This graph should be compared to the above-mentioned incidence graph of the Fano plane. It belongs to the family of so-called Generalized Pertersen graph (this is  $GP(8,3)$). We refer to \cite{Cox} (see in particular Fig. 3.3.c on page 22) for further informations. Justifications for its use in rank interpolation  can be found in Proposition \ref{trivalent} of Section \ref{class-rg74} below and the paragraph following it. 

Our main results on the geometric structure of triangle polyhedra of rank \sq\ are Theorem \ref{class} and the results in Subsections \ref{localsq} and \ref{asymsq}, that we summarize as follows. By  \emph{complex of rank \sq} we mean a compact CW-complex with triangle faces whose universal cover is a polyhedron of rank \sq\ (see Definition \ref{comp74}).

\begin{theorem}\label{th4}
 There are precisely 13 orientable  complexes of rank \sq\ with one vertex (their universal covers  are triangle polyhedra of rank \sq).
Moreover, 
\begin{enumerate} 
\item three of them have an abelianization of non zero rank (which can be 1 or 2) and in particular they don't have the property T of Kazhdan,
\item all of them contain copies of the free abelian group $\ZI^2$ and satisfy the following additional property: for any copy of 
 $\ZI^2$ in $\G$, there is a $\gamma\in \G$ such that the pairwise intersection of the subgroups  $\gamma^n \ZI^2\gamma^{-n}$, $n\in \ZI$,  is reduced to the identity. 
\end{enumerate}
\end{theorem}


Tits showed  in \cite{Tits74}  that a triangle polyhedron all of whose links correspond to a projective plane are buildings (see also \cite{rg2,BB}). It follows that they have exponential branching.  
 In the case of general triangle polyhedra, the combination of germs of flats can be very intricate and, depending on their relative position, does not necessarily ``integrate'' to actual flats in $X$. We make this precise in  Section \ref{localsq} where we study the local flat structure imposed  by the rank \sq\ 
(this should be compared to the results in  Section \ref{meso}), and give concrete examples of complexes of rank \sq\ which have exponential flat branching.  We do not know if a triangle polyhedron of rank \sq\ contains flats in general (see in particular Question \ref{Qbest}).%

In \cite{Gar} Garland introduced another local invariant for non positively curved polyhedra of  dimension 2 (and larger), called the \emph{$p$-adic curvature}. Given a vertex $A$ of a polyhedron $X$ the $p$-adic curvature of $X$ at $A$ is  defined to be the first non-zero eigenvalues $\lambda_1$ of the Laplacian on the link of $A$. Then he proved his famous vanishing cohomology results under the assumption $\lambda_1>1/2$, which eventually lead to the $\lambda_1>1/2$ criterion for property T \cite{Zuk,Pan,BS}. 
The first eigenvalue of $L_\sqm$  is 
\[
\lambda_1(L_\sqm)=0.42...<1/2.
\] 
We do not know whether there are groups of rank \sq\ which have property T (in fact Item (1) in Theorem \ref{th4}   first seemed quite unexpected to us). Property T for triangle buildings was first established in \cite{Cart} and the  proof given there, based on (local) spherical analysis, fails to apply in our context because the automorphism group of $L_{7/4}$ is not sufficiently transitive (see Section \ref{localsq}).



\subsection{The Baum-Connes conjecture for triangle groups} Our general criterion for proving property RD in the above framework is the following result.

\begin{theorem}\label{th1} 
Let $\G$ be a triangle group and let $\ell$ be the length on $\G$ induced by the 1-skeleton of a triangle polyhedron $X$ on which  $\G$ acts isometrically and properly. Then $\G$ has property RD with respect to $\ell$.
\end{theorem}

In particular  the (thirteen)  groups of rank \sq\ described in Theorem \ref{th4} have property RD.   The classical scheme for establishing property RD  consists in reducing the convolution product to partial convolutions over simpler triangles and  we proceed exactly in the same way in the present paper (see Section \ref{reminder} and the references therein). As in \cite{RRS,Laf-rd}  triangles will be reduced to flat equilateral triangles.    
Our contribution is in Section \ref{rd-triangles} and concerns the geometrical part of the proof. 

Note that symmetric spaces tools (e.g.\ the retraction onto apartments that was useful for buildings in \cite{RRS} or computations in $\SL_3$ as in \cite{Laf-rd})  are not available in our context. 
According to a conjecture of Valette \cite[page 66]{Val-bc} property RD should hold  for every groups properly isometrically and cocompactly on an affine building or a Riemannian symmetric space. If true, as Theorem \ref{th1}  suggests, it might hold  even more generally in situations where  rank interpolations is available. Understanding to what extend Theorem \ref{th1} generalizes to groups acting on other type of (say, geometrically finite but not necessarily symmetric) CAT(0) simplicial complexes  is an interesting open problem (compare this to Subsection \ref{meso-intr}).

Theorem \ref{th1} has the following consequence, which is a straightforward application of Lafforgue's Theorem \cite{Laf-bc}. 

\begin{corollary} \label{BC-cor}
Let $\G$ be a countable group admitting a proper, isometric, and cocompact action on a triangle polyhedron $X$. Then $\G$ satisfies the Baum-Connes conjecture, i.e.\ the Baum-Connes assembly map 
\[
\mu_\reg : K_*^\mathrm{top} (\G) \to K_*(C^*_\reg(\G))
\]
is an isomorphism. 
\end{corollary}

See \cite{BCH,Sk,Val-bc} for information on the Baum-Connes conjecture. 
We simply comment here that Lafforgue  considered in \cite{Laf-T} a strengthening of property T which holds  for cocompact lattices in $\SL_3(\QI_p)$ (in particular $\tilde A_2$-groups) but fails for every hyperbolic groups. This version of property T can be seen as an obstruction for proving the Baum-Connes conjecture with coefficients using Banach $KK$-theory (see \cite{Laf-T} where a proof of the Baum-Connes conjecture with coefficients for any hyperbolic groups is announced). It would  be interesting in that respect to determine `up to what rank' (necessarily $<2$) Banach $KK$-theory techniques can be applied in the framework of triangle polyhedra to get  Baum-Connes with coefficients. Note  that the construction of  Kasparov's element $\gamma$ and of the homotopy between $\gamma$ and 1 in (asymptotic versions of) $KK^\mathrm{ban}_\G(\CI,\CI)$ are technically easier  to perform in the context of triangle polyhedra than in the general case \cite{Laf-bc} of (strongly) bolic spaces, and coefficients appearing in the homotopy should be controllable to some extend (see \cite{Laf-bc}  and  Section 4 in \cite{Sk}).  The first-named author proved in \cite{cras} that some groups of local rank $\tdm$ have the Haagerup property and thus satisfy the Baum-Connes conjecture with coefficients by Higson-Kasparov's Theorem \cite{HigsonKasp} (and $\gamma=1$ in Kasparov's  $KK_\G(\CI,\CI)$ by \cite{Tu}). The proof of strengthened property T for cocompact lattices of $\SL_3(\QI_p)$ in \cite{Laf-T} relies on symmetric spaces tools and it is not clear at all that the same holds when the rank is (even slightly) lower, e.g.\ for some groups of rank \sq.

\newcommand{\sr}{\mathrm{sr}}
\subsection{$C^*$-algebraic  rank} Let $A$ be a unital $C^*$-algebra. The stable rank $\sr(A)$ of $A$ is an invariant of $A$ taking values in $\{1,2,\ldots\}\cup \{\infty\}$ which was introduced by Rieffel \cite{Ri83}.  
In the commutative case $\sr(A)$ behaves  as a dimension. Thus for a compact space $X$ and $A=C(X)$ the $C^*$-algebra of complex-valued function on $X$ one has 
\[
\sr(A)=\lfloor \dim X/2\rfloor+1.
\] 
In particular 
\[
\sr(C_\reg^*(\ZI^2))=2
\]
where $C_\reg^*(\ZI^2)\simeq C(\TI^2)$ is the $C^*$-algebra of the abelian free group $\ZI^2$.  We are interested here in the stable rank of reduced $C^*$-algebras of non amenable countable groups where, as opposed to the commutative case, an interpretation of $\sr(A)$ as a ``dimension" of $A$ is far less evident.     
In another direction, we mention that the case of nuclear (simple) algebras received much attention recently in connection to Elliott's classification program (see e.g.\ \cite{Toms} and references).  Villadsen \cite{Vil} constructed for any integer $n$ a  simple, separable and unital AH-algebra  of stable rank $n$.

We investigate here the relationships between the ``asymptotic dimension of $\G$" (especially from the intermediate rank point of view) and the stable rank of $C^*_\reg(\G)$, in the case of triangle groups. 

   A unital $C^*$-algebra $A$ has stable rank 1 if and only if the group $\GL(A)$ of invertible elements of $A$ is norm dense in $A$. There are well-known structural consequences of the stable rank 1 condition (see \cite{Bla}), especially concerning non stable $K$-theory properties of $A$. For instance the map $U(A)/U(A)_0\to K_1(A)$ from the quotient of the unitary group of $A$ by the connected component of the identity to the first $K$-theory group of $A$ is an isomorphism.

       In \cite{DHR} Dykema, Haagerup and R\o rdam proved that if $\G_1$ and $\G_2$ are two countable groups with $|\G_1|\geq 2$ and $|\G_2|\geq 3$ then
 \[
 \sr(  C_\reg^*(\G_1*\G_2))=1.
 \]
In particular  for the free groups $F_n$ on $n\geq 2$ generators one has $\sr(  C_\reg^*(F_n))=1$. In \cite{DH} Dykema and de la Harpe generalized these results and proved that if $\G$ is a torsion free non elementary hyperbolic group, or a  cocompact lattice in a real, noncompact, simple, connected Lie group of real rank one with trivial center, one has 
 \[
 \sr(  C_\reg^*(\G))=1.
 \]
We also mention that the rank of a group and the stable rank of its reduced $C^*$-algebra are known to be related to each other in the realm of Lie groups.  In \cite{Sud97} Sudo proved that for a  connected noncompact real semisimple Lie group $G$ the stable rank of 
$\sr(C^*_\reg(G))$ is 1 if the real rank of $G$ is 1, while it is 2 if the real rank of $G$ is $\geq 2$. It is unknown if a similar dichotomy holds true for cocompact lattices in real Lie groups   (see Problem 1.8 in \cite{DH}). The $p$-adic case is open as well, and in particular we don't know what the stable rank of the reduced $C^*$-algebra of  $\tilde A_2$-groups  is.

\begin{theorem}\label{th7} Let $X$ be a  complex of rank \sq\ and let  $\G=\pi_1(X)$ be the fundamental group of $X$. Then the reduced $C^*$-algebra $C^*_\reg(\G)$ of $\G$ has stable rank 1.
\end{theorem}

The proof of this result occupies  Section \ref{pfth7}. We  use a sufficient condition for stable rank 1 of Dykema and de la Harpe \cite{DH}, which is recalled at the beginning of Section \ref{pfth7}. 

All previously known reduced group $C^*$-algebras with stable rank 1 were related to free products or hyperbolicity.  In our case we know from Theorem \ref{th4} that there exist groups of rank \sq\ with  exponential branching and containing 
  infinitely many subgroups isomorphic to  $\ZI^2$ (each of them further satisfying the conditions in item (2) of this Theorem).  
 These groups are neither Gromov hyperbolic nor they are decomposable as non trivial free products. 
 In many respect they are actually \emph{closer to $\tilde A_2$-groups than to hyperbolic groups} (we remark that our proof of Theorem \ref{th7}, however, definitely fails in the rank 2 case).  

\newcommand{\rr}{\mathrm{rr}}
\newcommand{\sa}{\mathrm{sa}}

Another invariant of $A$,  the real rank,  was defined by Brown and Perdesen in \cite{BP}.  It is denoted $\rr(A)$ and takes values in $\{0,1,\ldots\}\cup \{\infty\}$.
  A unital $C^*$-algebra $A$ has real rank 0 if and only if  $\GL(A_{\sa})$ is dense in $A_{\sa}$, where the subscript $_\sa$ denotes the self-adjoint subspace of $A$.
  In the commutative case one has $\rr(C(X))=\dim X$, where $X$ is a compact space, and in general the following relation holds for a $C^*$-algebra $A$ (see \cite{Bla}):
  \[
  \rr(A)< 2 \sr(A).
  \]
  Thus the real rank of the reduced $C^*$-algebra of fundamental groups of compact complexes of rank \sq\ is at most 1. Let us now show that it is 1.

Recall  the following conjecture of Kaplansky and Kadison: for any torsion-free countable group $\G$, the reduced $C^*$-algebra of $\G$ has no idempotent besides 0 and 1.
As is well-known, this conjecture is a consequence of the surjectivity of the Baum-Connes assembly map  $\mu_\reg$ (see \cite{BCH,Val-bc}).   Thus, the absence of non trivial projection in the $C^*$-algebras considered in Theorem \ref{th7} is part of  Corollary \ref{BC-cor}.  
On the other hand in a  $C^*$-algebra of real rank 0 every self-adjoint element can be approximated by self-adjoint elements with finite spectrum, so in particular real rank 0 implies the existence of many non trivial projections \cite{Bla}.  

Summarizing, the following is a straightforward consequence of Theorem \ref{th7}, Corollary \ref{BC-cor}, and known facts. 

\begin{corollary}\label{corth7} Let $X$ be a complex of rank \sq\ and let  $\G=\pi_1(X)$ be the fundamental group of $X$. Then the reduced $C^*$-algebra $C^*_\reg(\G)$ of $\G$ has real rank 1.
\end{corollary}

It would be interesting to have a direct proof of this result.

\newcommand{\ptA}{A}

\subsection{Mesoscopic rank, mixed local rank}\label{meso-intr} Let us now come to our last  (and most refined) notion of intermediate rank.  Interpolation  occurs here  from  local  to global, i.e., we aim at measuring the proportion of ``flat pieces" of a space which are strictly in between the microscopic and the macroscopic scale.  

Let $X$ be a CAT(0) space of dimension 2 (without boundary) and $\ptA$ be a point of $X$. Consider the function $\varphi_\ptA : \RI_+ \to \NI$ which associated to an $r\in \RI_+$ the number of flat disks in $X$ of center $\ptA$  \emph{which are not included in a flat of $X$}. It can be seen as a way to measure the quantity of flats in $X$ which are situated strictly in between  local flats (i.e.\ flats of the tangent cones at $\ptA$) and  global ones (i.e.\ isometric copies of $\RI^2$ in $X$ containing $\ptA$).  We call $\varphi_\ptA$ the \emph{mesoscopic rank profile} of $X$ at $\ptA$ (or simply \emph{mesoscopic profile} for short). 

In the case of triangle polyhedra of extremal rank (rank 1 or rank 2) the mesoscopic profile trivializes  as follows. 

\begin{proposition} Let $X$ be a triangle polyhedron. If $X$ is hyperbolic (i.e.\ has rank 1) then its mesoscopic profile at every point is compactly supported. On the other hand $X$ is a buildings (i.e.\ has rank 2) if and only if its mesoscopic profile vanishes identically at every point.
\end{proposition}

Phenomenons start to appear for polyhedra of rank $\tdm$. The following graph, as one can show,  is the mesoscopic profile at some vertex of the triangle polyhedron of rank \td\ constructed in \cite{cras}. 

\begin{figure}[htbp]
\centerline{\includegraphics[width=8.5cm]{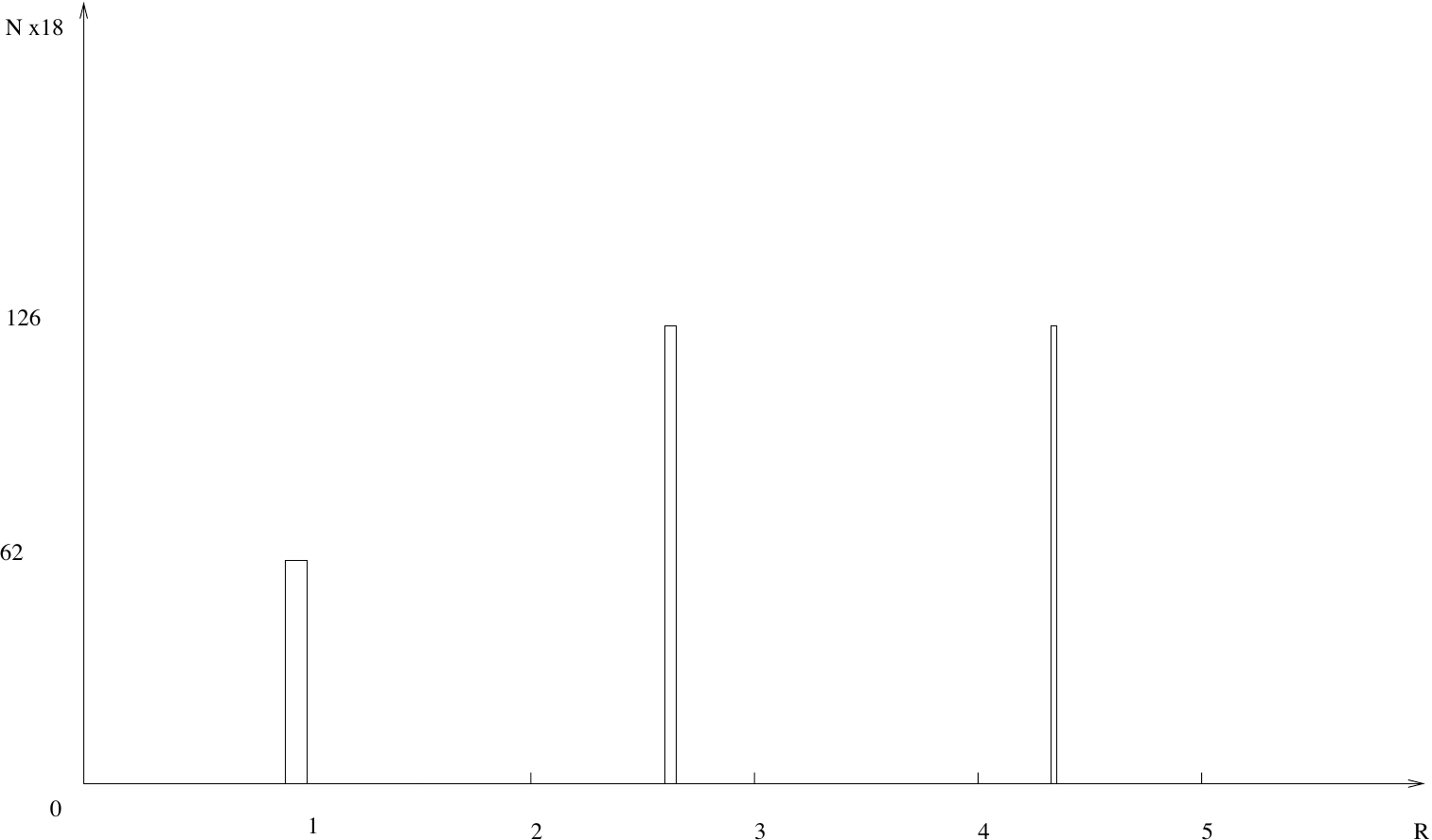}}
\caption{Mesoscopic  profile of the rank $\tdm$ polyhedron of \cite{cras}}\label{profile}
\end{figure}

\noindent Here the mesoscopic profile is is bounded (by $2268=126\times18$) and its support is an infinite union of disjoint intervals whose length tends to 0 at infinity.  

Intermediate rank at the mesoscopic scale is defined as follows.

\begin{definition}\label{def7}
A CAT(0) space $X$ of dimension 2 is said to be \emph{of mesoscopic rank} at a point $A$ if  the support of $\varphi_\ptA$ contains a neighborhood of infinity. 
\end{definition}

The meaning of Definition \ref{def7} is clear:  in a space $X$ of mesoscopic rank one can   continuously rescale   the radius   of  disks centered at $A$ 
which are flat but not included in flats of $X$ from some constant $C$ up to $\infty$. 
It is trivial to work out examples of 2-dimensional CAT(0) spaces $X$ without boundary which are of mesoscopic rank at \emph{some} point $A$.  
What we aim to construct here are $X$ for  which the set of $A$ satisfying this property is a (uniform) \emph{lattice} in $X$. 
We say that a countable group is \emph{of mesoscopic rank} if it admits a free and cocompact isometric action on a CAT(0) space $X$ of dimension 2 which is of mesoscopic rank at some  point. 

\medskip
Here are several sources of examples:

\medskip

\emph{(a) Surgery of exotic buildings.}  Recall that an affine building is said to be \emph{exotic} if it is cocompact but not classical, i.e.\ not associated to an algebraic group over a local field (see e.g.\ \cite{Ronan86,RonanTits,Cart,toulouse}). Fundamental work of Tits \cite{Tits74} led to the complete classification of affine buildings of dimension $\geq 3$: they all are classical.   
    The situation is entirely different in dimension 2 (see \cite{henri,quasiper} and references). 
        In Section 3 of \cite{toulouse},  the first author constructed an exotic triangle buildings $\tilde P$   which is the universal cover of a compact complex $P$ with two vertices (links at these vertices are trivalent, i.e.\ they are associated  to the Fano plane). 
    In fact   the fundamental group $\pi_1(P)$ of $P$ has finite index in the automorphisms group of $\tilde P$ (\cite[Th\'eor\`eme 7]{toulouse}).

By surgery on the complex $P$, one can  construct a compact complex $V_{\bowtie}$ with 8 vertices, whose universal cover is a CAT(0) space of dimension 2 without boundary (see Section \ref{meso} for details).   
This complex has \emph{mixed local rank}:  2 of its vertex have rank 2 and the 6 others have rank \td. 

We call the fundamental group $\G_{\bowtie}=\pi_1(V_{\bowtie})$ of $V_{\bowtie}$ a \emph{group of friezes} and its universal cover $\tilde V_{\bowtie}$ a \emph{complex of friezes}   (`groupe de frizes' and `complexe de frizes'  in french).

\medskip

\emph{(b) Polyhedra of rank \sq.}  In the classification of orientable complexes of rank \sq\ with one vertex of Section \ref{class-rg74} there is one, namely  $V_0^1$  in the notations introduced there, that has quite a distinctive intermediate rank property: in a sense that  will be  make precise in Sections \ref{class-rg74}  and \ref{meso}, its universal cover has \emph{the ``maximum" asymptotic branching
  within the range allowed by $L_\sqm$}.  As the local analysis in Section \ref{localsq}  will show,  that this upper-bound is ``attained"  is  quite  remarkable.  Its proof is the first step towards Item (b) of the following theorem.

\begin{theorem}\label{th11}The following groups are of mesoscopic rank.
\begin{itemize}
\item[(a)] The group of friezes $\G_{\bowtie}$ acting on the complex of friezes.
\item[(b)] The fundamental group of the  complex  $V_0^1$ (which is of rank \sq) acting on its universal cover $\tilde V_0^1$.
\end{itemize}
\end{theorem}

Even more, these complexes  have \emph{exponential mesoscopic rank} in the sense that their  mesoscopic profile  converges exponentially to infinity at infinity. 
Recall that $\tilde V_0^1$ in (b)  has all links isomorphic to $L_{\sqm}$ and  is \emph{transitive on vertices},  showing that  extremely homogeneous local data, that precludes  in particular mixed the local rank and spaces built out of different shapes, may still create  `singular' flats disks in $X$ at the mesoscopic scale (homogeneity of $L_{\sqm}$  is studied in Subsection \ref{localsq}).   This cannot happen for the \emph{most} homogeneous local data (i.e.\ spherical buildings), as we already saw. 

\bigskip

The  group $\G_{\bowtie}$ is a triangle group, as one can see after a suitable subdivision of the  complex of friezes.  Thus it satisfies property RD and the Baum-Connes conjecture.  
However,   the proof of Theorem \ref{th1} is \emph{very sensitive to the ambient geometry}, and establishing  property RD for $\Gamma_{\bowtie}$ directly (without subdivising) would  further increase the technical difficulties of Section \ref{rd-triangles}. In fact it is while looking for a way to bypass this technicalities in the case of the \emph{Wise group} that we first encountered mesoscopic rank phenomenon.  We shall now briefly  discuss this as a conclusion.

Write $\G_W$ for Wise's non Hopfian group, as constructed in \cite{Wise}. Recall that $\G_W$ is the fundamental group of a compact complex whose universal cover is a 2-dimensional CAT(0) space $\tilde W$. Then one can prove that:

\begin{itemize}
\item[(c)]\emph{Wise's group $\G_W$ is of mesoscopic rank.}
\end{itemize}

\noindent The proof  is omitted here (it is  similar to that of $\G_{\bowtie}$). In \cite{questions} Section 5.2  the question is raised of whether $\G_W$ has property RD or not (see also 6.6 in \cite{questions},  where $\G_W$ is proposed as a possible counter-example to property RD). As noted there, $\G_W$ does not acts on any cube complex (which implies property RD by Theorem 0.4 of \cite{Chat-Ruane}) nor it is relatively hyperbolic (it is actually of exponential branching). 

 Inspection of  the proof Theorem \ref{th1}  in the case of $\tilde W$ reveals that the situation is slightly worse than in the case of the non-subdivised $V_{\bowtie}$ but we believe, nevertheless,  that  $\G_W$ has property RD.  
(This has now been established in \cite{rdwise}.)
 
  Theorem \ref{th11} is proved in Section \ref{meso}. The proof of (b) is illustrated on Figure \ref{bibifurq}. 

\bigskip

\noindent \emph{Acknowledgements.} The second author was supported by an EPDI Post-doctoral Fellowship and is grateful to IHES, the Isaac Newton Institute for Mathematical Sciences, and the Max-Planck-Institut f\"ur Mathematik, for their hospitality over the 2005-2006 and 2006-2007 academic years, while the present work was carried out. 

\setcounter{tocdepth}{1}
\tableofcontents

\section{Property RD and polynomial branching}\label{reminder}

\renewcommand{\t}{{\mathfrak{s}}}

Let $\G$ be a countable group. A  triple $(x,y,z)\in \G^3$  such that $xy=z$ is called a triangle in $\G$.  Given a set $\t$ of triangles  in $\G$,  finitely supported functions $f,g\in \CI\G$, and $z\in \G$,   define $f*_\t g(z)$ by the expression 
\[
f*_\t g (z)=\sum_{(x,y,z)\in \t} f(x)g(y)
\]
if a triangle of the form $(x,y,z)$ belongs to $\t$, and 0 otherwise. The convolution product over the family of all triangles in $\G$ is  written $f*g$.  Let $\ell$ be a   length on $\G$, i.e.\  a non negative function $\ell$ on $\G$ such that $\ell(e)=0$,  $\ell(x)=\ell(x^{-1})$ and $\ell(xy)\leq \ell(x)+\ell(y)$ for $x,y\in \G$. Then $\t$ is said to have property RD with respect to $\ell$ if one can find a polynomial $P$ such that for any $r\in \RI_+$ and $f,g\in \CI\G$ with $\supp(f)\subset B_r$ one has
\[
\|f*_\t g\|_2\leq P(r)\|f\|_2\|g\|_2.
\]
If $*_\t=*$, i.e.\ if $\t$ consists of all triangles in $\G$,  then the group $\G$ is said to have property RD with respect to $\ell$ (see \cite{Jol-def}). For finitely generated groups this is  independent of the choice of $\ell$ among  word metrics associated to a finite generating sets,  so we simply speak of  property RD for  $\G$  in that case.     Note that it is sufficient to check  the above inequality on non negative functions $f,g\in \RI_+\G$. 
A standard approach to prove  property RD for $\G$ consists in  reducing $\ell^2$ estimates over $*$ to  estimates over simpler partial convolutions $*_\t$ (see \cite{RRS,Laf-rd,Chat-phd,Talbi}).  
In this section we prove Lemma \ref{rd-retr} and Lemma \ref{rd-dual}, which are the two main known tools for reducing convolution products, and we introduce  a notion of polynomial branching for countable groups endowed with a length. 

\begin{remark} Our  basic framework for this section will be that  of a countable group endowed with a length. This has the advantage of simplifying the exposition without hiding the important issues and this is  well adapted to groups  acting freely and simply transitively on the vertex set of a triangle polyhedron (which is the case, for instance, of the groups of rank 7/4 constructed in Section \ref{class-rg74}). As shown in \cite{RRS}  the appropriate tools for generalizing these results to  non necessarily simply transitive action are \emph{transitive groupoids}. This is discussed in more details at the end of Subsection \ref{groupoids}.
\end{remark}

\subsection{Statement and proof of Lemma \ref{rd-retr}.} \label{sub21}

Let $\G$ be a countable group and $\ell$ be a length on $\G$. A \emph{3-path} from the identity $e$ to a $z\in \G$ is a triple  $\gamma=(a_3,a_2,a_1)$ in $\G^3$ such that $z=a_3a_2a_1$.

\begin{definition}\label{HK}   
 A  $\G$-indexed family of 3-paths in $\G$, i.e.\ family $C=(C_z^r)_{z\in \G,~~r\in \NI^*}$  where  $C_z^r$ is a set of 3-paths from $e$ to $z$ in $\G$ for every $z\in \G$ and $r\in \RI^*$, is said to have \emph{polynomial growth} if there exists a polynomial $p_1$ such that for any $r\in \RI_+$ and  any $z\in \G$ one has 
$\#C_z^r\leq p_1(r)$.
\end{definition}

Let $\t$ and $\t^-$ be two sets of triangles in $\G$ and $C=(C_z^r)_{z\in \G,~~r\in \RI_+}$ be a $\G$-indexed set of 3-paths. For $(u,v,w)\in \t^-$ and $r\in \RI_+$ define $D^r_{(u,v,w)}$ to be the set of triple $(a,b,c)$ in $\G^3$ such that $(b^{-1},u,a)\in C^r_{b^{-1}ua}$, $(c^{-1},v,b)\in C^r_{c^{-1}vb}$ and $(c^{-1},w,a)\in C^r_{c^{-1}wa}$. 

Given $x\in \G$ we often write $|x|$ for $\ell(x)$. 

\begin{definition} \label{retr-tr}
One says that   \emph{$\t^-$  is a retract of $\t$ along $C$} if there exists a polynomial $p_2$ such that for every $(x,y,z)\in \t$ there exists  $(u,v,w)\in \t^-$ with $|u|\leq p_2(|x|)$ and   $(a,b,c)\in D^{|x|}_{(u,v,w)}$ such that $b^{-1}ua=x$ and $c^{-1}wa=z$.
\end{definition}

The idea of retracting to simpler sets of triangles originates in \cite{RRS}. Definition \ref{HK} corresponds to Property $H_\delta$ and a part of Property $K_\delta a$ in \cite{Laf-rd}, and 
Definition \ref{retr-tr} is another part of Property $K_\delta a$ in \cite{Laf-rd} (see also \cite{Chat-phd}  and Section 1.3 in \cite{Talbi}).  
 Our assumptions here are actually slightly weaker (in particular we do not to assume triangles in the retraction $\t^-$  to be ``balanced" at this stage,  cf.\ Subsection \ref{RRS-dual}).

The following lemma  was first proved by Haagerup in \cite[Lemma 1.4]{Haa} in the case of finitely generated free groups, where the set $\t$ of all triangles consists of \emph{tripod triangles} (i.e.\ triangles which retract to $\t^-=\{(e,e,e)\}$). It has then been extended in \cite{Jol-def,Harpe-rd,RRS,Laf-rd}.
The statement below corresponds to Proposition 2.3 and a part of Theorem 2.5 in \cite{Laf-rd} (compare \cite{RRS}), and the proof given below reproduces the arguments on pages 258 and 259 of this paper.

\begin{lemma}\label{rd-retr} Let $\G$ be a group and $\ell$ be a length on $\G$. Fix a $\G$-indexed set of path $C$ in $\G$,   a family $\t$ of triangles  in $\G$ and  a retract $\t^-$ of $\t$ along $C$ as in Definition \ref{retr-tr}.  
Assume that $C$ has polynomial growth.  Then  for any $r\in \RI_+$ and $f,g\in \RI_+\G$ with $\supp(f)\subset B_r$ there exist two functions $i,j\in \RI_+\G$ with $\supp(i)\subset B_{p_2(r)}$ such that 
\[
\|f*_\t g\|_2\leq \sqrt{p_1(r)}\|i*_{\t^-} j\|_2
\]
and $\|i\|_2\leq \sqrt{p_1(r)}\|f\|_2$, $\|j\|_2\leq \sqrt{p_1(r)}\|g\|_2$ where $p_1,p_2$ are as above. Thus property RD holds for $\t$ provided it does for $\t^-$. 
\end{lemma}
\begin{proof} Let $f,g,h\in \RI_+\G$ with $\supp(f)\subset B_r$. As any triangle in $\t$ can be retracted to $\t^-$ one has 
\begin{align*}
\langle f*_\t g\mid h\rangle &=\sum_{(x,y,z)\in \t} f(x)g(y)h(z)\\
&\leq \sum_{(u,v,w)\in \t^-, ~~|u|\leq p_2(r)} \sum_{(a,b,c)\in D^{r}_{(u,v,w)}} f(b^{-1}ua)g(c^{-1}vb)h(c^{-1}wa)\\ 
&=\sum_{(u,v,w)\in \t^-}\tr(R_uS_vT_w)\end{align*}
where $R_u$ is the operator on $\ell^2(\G)$ defined for $|u|\leq p_2(r)$ by 
\[
\langle R_u\delta_a\mid \delta_b\rangle=f(b^{-1}ua)
\]
if $(b^{-1},u,a)\in C_{b^{-1}ua}^{r}$ and 0 elsewhere. The operators $S_v,T_w$ are defined similarly for any $v,w$ using $g,h$. 
Set 
\[
i(u)=\|R_u\|_2,~~ j(v)=\|S_v\|_2~~~\mathrm{and}~~~k(w)=\|T_w\|_2.
\]
As $\tr(R_uS_vT_w)\leq \|R_u\|_2\|S_v\|_2\|T_w\|_2$ one has
\[
\langle f*_\t g\mid h\rangle \leq \sum_{(u,v,w)\in \t^-}i(u)j(v)k(w)=\langle i*_{\t^-} j\mid k\rangle.
\]
On the other hand
\[
\|i\|_2^2=\sum_{a,b,u|~  (b^{-1},u,a)\in C^{r}_{b^{-1}ua}}f(b^{-1}ua)^2\leq \sum_{x}\#C_x^{r} f(x)^2\leq p_1(r)\|f\|_2^2
\]
and similarly $\|j\|_2^2\leq p_1(r)\|g\|_2^2$ and $\|k\|_2^2\leq p_1(r)\|h\|_2^2$.  The lemma follows from the Cauchy-Schwarz inequality for $h=f*_\t g$.
\end{proof}

\subsection{Polynomial branchings}\label{sec22}
Let $\G$ be a countable group and $\ell$ be a length on $\G$. Let $\kappa\geq 1,\delta\geq 0$. A 3-path $(a_3,a_2,a_1)$ in $\G$ is said to be  $(\kappa,\delta)$-geodesic if $|a_1|+|a_2| +|a_3|\leq \kappa|a_3a_2a_1|+\delta$. 

\begin{definition}\label{polrk}
We say that  $\G$ has \emph{polynomial branching} with respect to $\ell$ if there exists $\kappa\geq 1$, $\delta\geq 0$, a family  $C=(C_z^r)_{z\in \G,~~r\in  \NI^*}$ of sets $C_z^r$ of $(\kappa,\delta)$-geodesic 3-paths  from $e$ to $z$ with polynomial growth (see Definition \ref{HK}),  a subset $\t$ of triangle in $\G$ which is a retract along $C$ of the family of all triangles in $\G$, and a  polynomial $p_3$ such that the for every $z$ in $\G$ and every $r\in \RI_+$ the number of triangles in $\t$ of the form $(x,y,z)$ with $|x|\leq r$ is no greater than $p_3(r)$. 
 \end{definition}

The present paper is mostly concerned with the case where $\kappa$ can be chosen to be equal to 1 in the above definition, a property we call \emph{geodesic polynomial branching}.


\begin{proposition}\label{pol-rd}
Let $\G$ be countable group with length $\ell$.  If $\G$ has polynomial branching with respect to $\ell$ then it has property RD with respect to $\ell$. 
\end{proposition}
\begin{proof} Let $\t$ be a family of triangles in $\G$ as given by  Definition \ref{polrk}. By Proposition \ref{rd-retr} there exists a polynomial $p_1$   such that for any $r\in \RI_+$ and $f,g\in \RI_+\G$ with $\supp(f)\subset B_r$ there exist two functions $i,j\in \RI_+\G$ with $\supp(i)\subset B_{\kappa r+\delta}$ such that 
\[
\|f* g\|_2\leq \sqrt{p_1(r)}\|i*_{\t} j\|_2
\]
and $\|i\|_2\leq \sqrt{p_1(r)}\|f\|_2$, $\|j\|_2\leq \sqrt{p_1(r)}\|g\|_2$. We then have
\begin{align*}
\|i*_\t j\|_2^2&=\sum_{z\in \G} \left (\sum_{(x,y,z)\in \t,~~~|x|\leq \kappa  r+\delta} i(x)j(y)\right )^2\\
&\leq p_3(\kappa  r+\delta) \sum_z\sum_{(x,y,z)\in \t,~~~|x|\leq p_2(r)} i(x)^2j(y)^2\\
&=p_3(\kappa  r+\delta) \sum_{x,~~|x|\leq p_2(r)} i(x)^2 \sum_{(x,y,z)\in \t} j(x^{-1}z)^2\\
&\leq p_3(\kappa  r+\delta)\|i\|_2^2\|j\|_2^2.
\end{align*}
where $p_3$ is a polynomial as in Definition \ref{polrk}. Hence,
\[
\|f* g\|_2\leq \sqrt{p_1(r)p_3(\kappa  r+\delta)}\|i\|_2\|j\|_2\leq p_1(r)^{3\over 2}p_3(\kappa  r+\delta)^{1\over 2}\|f\|_2\|g\|_2
\]
which proves the Proposition.
\end{proof}

\begin{corollary}
An amenable group has polynomial branching if and only if it has polynomial growth \cite{Harpe}. 
\end{corollary}

\begin{proof} Recall that an amenable group has property RD if and only if it has polynomial growth \cite{Jol-def}. So if $\G$ is amenable of polynomial branching then it is has polynomial growth by Proposition \ref{pol-rd}.
(Note that polynomial branching is stable under direct but not semi-direct products.)
The converse is easily seen by choosing $C_z^r$ to be reduced to $\{(e,z,e)\}$ if $|z|\leq r$ and equal to $\{(zx^{-1},x,e): |x|\leq r\}$ if $r<|z|$,  and $\t$ to be the set of all triangles in the ball of radius $r$ of $\G$. 
\end{proof}

 Examples of groups with  property RD which don't have polynomial branching  include groups acting freely and simply transitively on triangles buildings (i.e., $\tilde A_2$ groups). Examples of non amenable groups with polynomial branching are as follows.

\begin{proposition}\label{ex} Let $\G$ be a finitely generated group. 
 If $\G$ is hyperbolic then it has polynomial branching (of degree 0). 
More generally, if $\G$ is hyperbolic relatively to finitely generated subgroups $\{\Lambda_{1},\ldots \Lambda_{n}\}$ of polynomial growth (for the induced length function)  have polynomial branching.
\end{proposition}

\begin{proof}
In the case of hyperbolic groups,  the required properties are satisfied if we choose $C_z^r$ to be the set of $\delta$-geodesic 3-paths from $e$ to $z$ of the form $(b,e,a)$ with $|a|\leq r$, for some $\delta$ large enough and $\t^-=\{(e,e,e)\}$. (See \cite{Harpe-rd}.)

So let $\G$ be a group which is hyperbolic relatively to $\{\Lambda_{1},\ldots \Lambda_{n}\}$. In fact we assume more generally  that $\G$ is (*)-relatively hyperbolic with respect to $\{\Lambda_{1},\ldots \Lambda_{n}\}$ in the sense of  Drutu and Sapir  \cite[Definition 2.8]{DrutuSapir}, for the length $\ell$ coming from some finite generating set of $\G$, and argue exactly  as in the beginning of the proof of Theorem 3.1 in \cite{DrutuSapir}. Thus   for $z\in \G$ we fix a simplicial geodesic $g_z$ from $e$ to $z$ and let    $C_z^r$ for  $r\in \RI_+$ be is the set of triples $(b,h,a)$ which are called central decompositions in \cite{DrutuSapir} (see Definition 3.3), where  $r$ is fixed to be equal to $r_1$ in their paper. Then one has  that $C=(C_z^r)_{z,r}$ has polynomial growth (see Lemma 3.3 in \cite{DrutuSapir}) and that all triangles of $\G$ retract along $C$ to the family $\t$ of all triangles in the subgroups $\Lambda_{i}$ (see the begining of the proof of Theorem 3.1).  Thus, if the subgroups $\Lambda_i$ have polynomial growth with respect to $\ell$ then the number of triangles in $\t$ with fixed basis $z$ is polynomial. This  proves the proposition.  
\end{proof}

Note that if in the above proof  the  $\Lambda_i$  are not of polynomial growth but  have property RD, then  one can still apply Lemma \ref{rd-retr} to deduce (following \cite{DrutuSapir}) that $\G$ has property RD as well (see Theorem 1.1 of \cite{DrutuSapir}). Accordingly, not all relatively hyperbolic groups with property RD have polynomial branching.

\begin{definition}
We say that a group $\G$  has  \emph{subexponential branching}  with respect to a length $\ell$ if all conditions in  Definition \ref{polrk} are satisfied except perhaps for the polynomial growth assumption on $p_1$ and $p_3$, which we now allow to be  subexponential,   i.e.\ $p_1$ and $p_2$ are non negative functions on $\RI_+$ such that $\lim_rp_1(r)^{1/r}=\lim_rp_3(r)^{1/r}=1$. 
\end{definition}

Since for $p_1,p_3$ of subexponential growth the function $p_1(r)^{3\over 2} p_3(\kappa r+\delta)^{1\over 2}$  has subexponential growth as well, the proof of Proposition \ref{pol-rd} shows that groups with subexponential branching satisfy a subexponential variation of property  RD where in the definition  $P$ is replaced by a function, say $\eta$, of subexponential growth. Amenable groups with subexponential property RD have subexponential growth. Indeed denoting  $\chi_s$  the characteristic function of the ball of radius $s$, one has
\[
|B_s| \leq \|\chi_s\|_\reg\leq \eta(s)\|\chi_s\|_2=\eta(s)\sqrt{|B_s|}
\]
where the first inequality follows from $|\sum_{x\in \G} f(x)|\leq \|f\|_\reg$
for every $f\in \CI\G$ by weak containment of the trivial representation in  the regular representation of $\G$ (we write $\|f\|_\reg$ for the norm of $f\in \CI\G$ acting by convolution on $\ell^2(\G)$). Thus  indeed $\G$ has subexponential growth (compare \cite{Jol-def,Con,Val-bc}) and in particular subexponential branching coincide with subexponential growth in the amenable setting. 
On the other hand arguing as in Proposition \ref{ex} we obtain that groups that are relatively hyperbolic with respect to groups of subexponential growth have subexponential branching. Taking free product $A*B$ of groups $A,B$ of intermediate growth (e.g.\ the groups of Grigorchuk, see \cite{Harpe} and the references therein) shows that the class of groups with (optimal) subexponential branching $\kappa$ vary when the growth of $\kappa$ varies (relying upon examples by A. Erschler).  Note also that for groups which are relatively hyperbolic with respect to $\{\Lambda_1,\ldots, \Lambda_n\}$, following \cite{DrutuSapir},  subexponential property RD is equivalent to subexponential property RD for the $\Lambda_i$. We do not know, however, the answer to the following `intermediate branching growth problem' in the case of triangle polyhedra.

\begin{question}
Are there finitely generated groups admitting a proper and cocompact  action on a triangle polyhedron $X$ which have geodesic intermediate (i.e.\ subexponential but not polynomial)  branching with respect to the length induced from the 1-skeleton of $X$?
\end{question}

In fact $\tilde A_2$ groups, and some other triangle groups constructed below, have geodesic exponential branching in the following sense.

\begin{definition} A group $\G$ is said to have exponential branching with respect to a length $\ell$ if it is not of subexponential branching.
\end{definition}

One defines similarly subexponential geodesic branching and exponential geodesic branchings. 

For future use (Section \ref{pfth7}) we end this subsection with a  discussion of the $\ell^2$ spectral radius property (see \cite{HRV,DH}) and related applications of property RD  to random walks on groups \cite{Grigo}. 
In Section 5 of \cite{Grigo} Grigorchuk and Nagnibeda considered the \emph{operator growth function of $\G$},  defined as $F_\reg(z)=\sum_{n} a_n z^n$ with coefficients 
\[
a_n =\sum_{|x|=n} u_x
\]
where $u_x$, $x\in \G$, is the canonical unitary corresponding to $\G$ in $C^*_\reg(\G)$ under the regular representation.
The radius of convergence $\rho_\reg$ of $F_\reg$ satisfies 
\[
{1\over {\rho_\reg}}=\limsup_{n\to \infty} \|a_n\|_\reg^{1/n}\leq \limsup_{n\to \infty} |S_n|^{1/n}={1\over {\rho}}
\] 
where $\rho$ is the usual (inverse) exponential growth rate of $\G$ with respect to $\ell$. Conjecture 2 in \cite{Grigo} states that $\G$ is amenable if and only if $\rho=\rho_\reg$.  It is proved in \cite{Grigo} that $\rho=\rho_\reg=1$ for amenable groups and that $\rho_\reg=\sqrt{\rho} <1$ for non amenable hyperbolic groups (recall that $\G$ is amenable if and only if $\rho=1$ by Kesten criterion).  Valette noted that the proof given of \cite{Grigo} was only using property RD for $\G$, and in fact  that \emph{radial property RD} was sufficient. This allows for instance to include every $\tilde A_n$-groups, for $n\geq 2$, in the above Conjecture 2 which have radial property RD thanks to the work of Valette  \cite{Val-rad} and  \'Swi\c atkowski \cite{Swiat-rad}.   

In the same way \emph{for a non amenable group $\G$  satisfying radial subexponential property RD with respect to $\ell$ one has $\rho_\reg=\sqrt{\rho} <1$}, so $\G$ satisfies conjecture 2 in \cite{Grigo}. The proof is exactly as in \cite{Grigo} (see also \cite{HRV}): by radial subexponential property RD we have 
\[
\|a_n\|_\reg\leq \kappa(n)\|a_n\|_2=\kappa(n)\sqrt{|S_n|}
\]
so 
\[
{1\over {\rho_\reg}}\leq \limsup_{n\to \infty}\kappa(n)^{1\over n}\sqrt{|S_n|}^{1\over n}={1\over \sqrt{\rho}}.
\]
As $\|a_n\|_2\leq \|a_n\|_\reg$ always holds $\rho_\reg\leq \sqrt{\rho}$ as well. The same argument also shows the $\ell^2$ spectral radius property for every element in the group algebra of $\G$ provided $\G$ has subexponential property RD, i.e.\  the spectral radius of every element $a\in \CI\G$ acting by convolution on $\ell^2(\G)$  is equal to
\[
\lim_{n\to\infty} \|a^{*n}\|_2^{1/n},
\]
since the radius of the support of $n$-th convolution product $a^{*n}$ is at most $n$ times the radius of the support of $a$. Summarizing we have  the following result  (compare Proposition 8 in \cite{HRV}, Section 3 in \cite{DH}, and Proposition 4 in \cite{Grigo}).

\begin{proposition}\label{spectr}
If  $\G$ is  a countable group with radial subexponential property RD,  then 
$\rho_\reg=\sqrt{\rho}$ and thus $\G$ satisfies conjecture 2 in \cite{Grigo}. If moreover  $\G$ has subexponential property RD (in particular if it has property RD), then it satisfies the  $\ell^2$ spectral radius property.
\end{proposition}

As noted at the end of Section 3 of \cite{DH} finitely generated groups of subexponential growth provide examples of groups with the $\ell^2$-spectral radius property which don't have property RD. All known examples of groups with the $\ell^2$-spectral radius property seems, however, to have subexponential property RD.

\subsection{Statement and proof of Lemma \ref{rd-dual}}\label{RRS-dual} We now recall the (crucial) analytical argument of Ramagge, Robertson and Steger    \cite[Lemma 3.2]{RRS}  for establishing property RD in the case of triangle buildings. For  a family $\t$ of triangles in $\G$ we call \emph{dual of $\t$} the family $\t^*$ of triangles of the form $(x^{-1},u,ux^{-1})$ and $(y,v^{-1},vy^{-1})$ with common basis $ux^{-1}=vy^{-1}$ whenever $(x,y,z)$ and $(u,v,z)$ are two triangle in $\t$ with common basis  $z$ (cf.\  property $K_\delta b$ in \cite{Laf-rd} or Definition 1.30 in \cite{Talbi}).  The proof below is contained in Lemma 3.2 of \cite{RRS}, see also the top of p. 260 in \cite{Laf-rd} or \cite{Talbi}.
One says that a family $\t$ of triangles is \emph{balanced} if there is a polynomial $p_4$ such that for  every $(x,y,z)\in \t$ one has $\max\{|y|,|z|\}\leq p_4(|x|)$.

\begin{lemma}\label{rd-dual}
Let $\G$ be a countable group endowed with a length $\ell$, and let $\t$ be a balanced family of triangles in $\G$. There exists a polynomial $p_4$ such that for $r\in \RI_+$ and  $f,g\in \RI_+\G$ with $\supp(f)\subset B_r$ one has 
\[
\|f*_\t g\|_2^2\leq \|\check f*_{\t^*} f\|_2\|(g\chi_{B_{p_4(r)}})*_{\t^*}\check g\|_2
\] where $\check h(z)=h(z^{-1})$ for $h\in \CI\G$. So  property RD holds for $\t$ provided it does for $\t^*$.
\end{lemma}

\begin{proof} For $f,g\in \RI_+\G$ with $\supp(f)\subset B_r$ one has
\begin{align*}
\|f*_\t g\|_2^2&=\sum_{z\in \G} \sum_{(x,y,z)\in \t}\sum_{(u,v,z)\in \t} f(x)g(y)f(u)g(v)\\
&\leq \sum_{z'\in \G} \sum_{x^{-1}uz'\in \t^*}\sum_{yv^{-1}z'\in \t^*,\ |y|\leq p_4(r)} f(x)g(y)f(u)g(v)\\
&=\sum_{z'\in \G} (\check f*_{\t^*}f)(z') ((g\chi_{B_{p_4(r)}})*_{\t^*} \check g)(z')
\end{align*}
as $|y|\leq p_4(r)$ for $|x|\leq r$ as $\t$ is balanced. The Lemma follows from the Cauchy-Schwarz inequality.
\end{proof}

\label{groupoids}
Let us conclude this section by recalling the generalization of the above to \emph{transitive} groupoids \cite{RRS}. This notably allows to prove property RD for countable groups whose length is coming from a general free isometric actions on metric spaces (rather than vertex-transitive actions). So let $\G$ be a countable group acting freely on a metric space $(X,d)$ and consider,  following \cite{RRS},   the countable groupoid $G=X\times_\G X$ of base $G^{(0)}=X/\G$. Let $\ell$ be the length on $G$ defined by $\ell([x,y])=d(x,y)$ for $[x,y]\in G$ and $B_r=\{[x,y]\in G,~\ell([x,y])\leq r\}$. Then $(G,\ell)$ is said to have property RD if the usual convolution estimate (with respect to  the groupoid law in $G$) is satisfied for $f,g\in \RI_+G$ with $\supp(f)\subset B_r$. All definitions presented in this section (in particular retractions along $\G$-indexed family of path of Subsection \ref{sub21} and the above dualization procedure)  extends to the case of $(G,\ell)$, and straightforward generalizations of Lemma \ref{rd-retr} and Lemma \ref{rd-dual} provide criteria for proving property RD for $(G,\ell)$.  In turn property RD for $(G,\ell)$ is easily seen to imply property RD for $\G$ with respect to the length induced from $d$ on one of its orbit in $X$ (see e.g.\ \cite[Prop. 2.1]{Laf-rd}).  

Note however that an extension of the techniques presented in this section to other---non transitive, but say, $r$-discrete and locally compact---groupoids is an  open problem in general, compare \cite{Laf-groupoids}  and the last sections of \cite{survey}.

\section{Proof of Theorem \ref{th1}}\label{rd-triangles}

The proof relies on several preliminary lemmas. 
Throughout the section we let $X$ be a fixed triangle polyhedron (Definition \ref{def1}). 
A curve between two vertices $A$ and $B$ of $X$ is said to be
\begin{itemize}
\item a \emph{geodesic segment} if its length equals the CAT(0) distance between $A$ and $B$. By the CAT(0) property there is a unique geodesic segment  between any two points  of $X$. 
\item a  \emph{simplicial geodesic segment} if it is simplicial, i.e.\ included in the 1-skeleton of $X$, and if its length coincide with the simplicial length between $A$ and $B$ in $X$, where the length of every edge in $X$ is normalized to 1.
\end{itemize}
A geodesic segment is called \emph{singular} if it is simplicial (up to parallelism this coincides with the usual definition in case $X$ is symmetric \cite[p. 322]{BH}).

\begin{definition} Let  $\gamma$ be a geodesic segment between two vertices $A$ and $B$ of $X$. One calls \emph{simplicial convex hull}   of $\gamma$ the union, denoted $\conv(\gamma)$,  of all triangles of $X$ whose three vertices belong  to simplicial geodesic segments from $A$ to $B$ in $X$. 
\end{definition}

By $\RI^2$ we mean the Euclidean plane endowed with the tessellation by equilateral triangles. Isometries are assumed to preserve the simplicial structures.   
 A \emph{flat} in $X$ is the image of an isometric embedding in $X$ of the Euclidean plane $\RI^2$. 
 A \emph{flat topological disk} in $X$ is the image an isometric embedding in $X$ of a topological disk of $\RI^2$. In particular
 a \emph{flat equilateral triangle} is the image of an isometric embedding of an equilateral triangle of $\RI^2$. 

Let $D$ be an open topological disk in $X$ with piecewise linear topological boundary of $\Delta$. Let $s$ be a point in $\Delta$ and $L$ be the link of $s$ in $X$. The disk $D$ determine a path $c$ in  $L$  from the two points  of $L$ corresponding  the incoming and outgoing segments of $\Delta$ at $s$. The angle between these segments, i.e.\ the angular length of  $c$ in $L$, is called the  \emph{internal angle} of $D$ at $s$ and is denoted by $\theta_s$.

\begin{lemma}\label{conv-triangles} 
Let  $\gamma$ be a geodesic segment between two vertices of $X$. Then   there exist finite sets $J$ and $J^\circ$ such that 
 \[
 \conv(\gamma)=\bigcup_{i\in J} G_i \cup \bigcup_{i\in J^\circ} S_i
\]
where 
\begin{enumerate}
\item $G_i$, $i\in J$, is  a  closed flat topological disks of $X$ which is,  under an isometry with a closed disk of $\RI^2$, a union of minimal galleries (see \cite{Ronan}) between two given vertices of $\RI^2$,
\item $S_i$, $i\in J^\circ$, are singular geodesic segments included in $\gamma$,
\item $S_i\cap S_j$, $i,j\in J^\circ$, is empty,  while $G_i\cap G_j$, $i,j\in J$, and $G_i\cap S_j$, $i\in J$, $j\in J^\circ$, are either empty or reduced to a vertex of $\gamma$.
\end{enumerate}  
\end{lemma}

\begin{proof}
Let $A_0, A_1,\ldots A_n$ be the set of vertices of $\gamma=[A_0,A_n]$. Denote by $I$ the set of integers $i\in [0\ldots n-1]$ for which the segment $[A_i,A_{i+1}]$ is non singular and let $I^\circ$ be the complement of $I$ in $[0\ldots n-1]$.  For each $i\in I$ let $G_i^0$  be the  gallery  from $A_i$ to $A_{i+1}$ in $X$,  which can be defined in this context as the union   the triangles of $X$ whose interior intersects  $]A_i,A_{i+1}[$.  We call the set 
\[
n_\gamma =\bigcup_{i\in I} G_i^0 \cup \bigcup_{i\in I^\circ} [A_i,A_{i+1}] 
\]
the \emph{nerve} of the simplicial convex hull of $\gamma$.  Note that   $G_i^0$,  $i\in I$, is a flat disk 
satisfying property (1) of the lemma. 
Let $J^\circ$ the set of $i\in I^\circ$ such that $i-1\notin I^\circ$. For $j\in J^\circ$ we denote by $S_j$  the union of segments $[A_i,A_{i+1}]$ for $i\in I^\circ$ such that $[j,j+1,\ldots i]\subset I^\circ$. 

Call a vertex $A_i$, $i\in I$, regular if $i-1\in I$ and if the distance in the link $L_i$ of $A_i$ in $X$ between the two edges corresponding to $n_\gamma$, say $e_i$ and $f_i$,  equals $2\pi/3$. Let $I_r$ be the set of $i\in I$ such that $A_i$ is regular and let  $J$ be the complement of $I_r$ in $I$. 

For every regular vertex $A_i$, $i\in I_r$ choose two edges $h_i^0$ and $h_i^1$ of the link $L_i$ such that the family $\{e_i,h_i^0,h_i^1,f_i\}$ forms a connected path in $L_i$
(there might be two such paths) and denote by $t_i^0$ and $t_i^1$ the triangles in $X$ containing $A_i$ and having  $h_i^0$ and $h_i^1$ respectively as basis (where we identified the link $L_i$ with the simplicial sphere of radius 1 in $X$). 

Let $j\in J$. Consider the largest integer $k< n$ such that for all integer $i< n$ with $j<i\leq k$ one has $i\in I$ and the vertex $A_i$ is regular. Consider the set $G_j^1$ defined as 
\[
G_j^1=\bigcup_{j<i\leq k} G_i^0~~\cup ~~\bigcup_{j<i\leq k} \{t_i^0\cup t_i^1\}.
\]
It is easy to see that $G_i^1$ is a flat disk which satisfies (1). Denote $B_j=A_{k+1}$ and fix, for every $j\in J$, an isometry $\varphi_j$ between $G_j^1$ and a closed disk $F_i$ of $\RI^2$. Let also $\tilde A_j$ and $\tilde B_j$ be the points in $\RI^2$ corresponding to $A_j$ and $B_j$ under this isometry and note that $F_i$ is included in the simplicial convex closure  $E_j$ of $\tilde A_j$ and $\tilde B_j$ in $\RI^2$ (which is a parallelogram).

Let $\cF_j$ be the (finite) set of closed disk of $\RI^2$ containing $F_j$ and which are union of minimal gallery from $\tilde A_j$ to $\tilde B_j$ (so every disk in $\cF$ is a subset of $E_i$). Consider the set $\cF_j^0$  of elements of $\cF_j$ which admit an isometric embedding in $X$ which coincide with $\varphi_j^{-1}$ on $F_j$. For each $j$ choose a disk in $\cF_j^0$ which has the maximal number of triangles and denote by $G_j$ its corresponding embedding in $X$. We will   show that the closed subset 
\[
C=\bigcup_{j\in J} G_j \cup \bigcup_{j\in J^\circ} S_j
\]
of $X$ coincide with the convex closure of $\gamma$.

Let $\cH_0$ be the set of simplicial geodesics of $C$ from $A_0$ to $A_n$ and $\cH_1$ be the set (which may be infinite a priori) of simplicial geodesics of $X$ from $A_0$ to $A_n$ which are not included in $C$. We will show that $\cH_1$ is empty. Note that by construction every simplicial geodesic of $C$ from $A_0$ to $A_n$  is the union of a simplicial geodesic in the flat disk $G_j$, $j\in J$ and the CAT(0) geodesic  $S_j$, $j\in J^\circ$.   For every $g_0\in \cH_0$ and $g_1\in \cH_1$ there exists by the CAT(0) property a finite family $\cD=\{D_0,\ldots, D_m\}$ of non empty topological disks of $X$ with disjoint interiors, which are union of triangles, which are filling $g_0\cup g_1$ in the sense that the subset $(g_0\cup g_1)\bigcup (D_0\cup\ldots \cup D_m)$ of $X$ is contractile, and such that for every $i=0\ldots m$ the intersection of the topological boundary $\del D_i$ of $D_i$ with $g_0\cap g_1$ consists of two points.
Moreover up to modifying  $g_0$ among elements of $\cH_0$ one  can choose  (for every $g_1\in \cH_1$) a $g_0$ such that the interior of $D_0\cup\ldots \cup D_m$ is disjoint from the interior of $C$. Let $\cA$ be the set of triples $(g_0,g_1,\cD)$ satisfying these conditions (so $\cA\to \cH_1$ which maps $(g_0,g_1,\cD)$ to $g_1$ is surjective).

Assume that $\cH_1$ is non empty and pick a  $(g_0,g_1,\cD)$ in $\cA$ such that the number of triangles of $\cD$ is minimal among all elements of $\cA$. Let $D$ be a  disk  in $\cD$. By construction the topological boundary $\del D$ of $D$ is included in $g_0\cup g_1$. 

Let $s$ be a 
vertex of $\del D\cap g_0$ which does not belong to $g_0\cap g_1$. Let us show that the internal angle $\theta_s$ of $D$ at $s$ is at least $\pi$. By (1) $\theta_s\neq \pi/3$ so  $\theta_s\geq 2\pi/3$.   Assume that  $\theta_s=2\pi/3$ and denote by $(x,s,t)$ and $(t,s,y)$ the corresponding triangles in $D$, where  $[x,s]$ and $[s,y]$ are two consecutive edges of $g_0$ because $s\notin g_0\cap g_1$. By definition of $S_j$, $j\in J^\circ$, the point $s$ does not belong to  $S_j$ (neither its interior nor its extremities) as this would contradict the fact that $\gamma$ is geodesic and the definition of $C$. In particular $s\in G_j$ for some $j\in J$ which in turn implies that $s\in \del G_j$.  Indeed otherwise the path $[x,t]\cup [t,y]$ would create with two edges of $C$ a cycle of length at most $\pi+2\pi/3<2\pi$ in the link of $s$ in $X$, contradicting the CAT(0) property.  Furthermore one has $s\neq A_j$ and $s\neq B_j$ as otherwise (one at least of) these points would be regular. It follows that $[x,s]\cup [s,y]$ is included in the boundary of $G_j$. Now by (1) the internal angle of $G_j$ at $s$ is at most $4\pi/3$ so  as $X$ is CAT(0) this angle  exactly equals $4\pi/3$. It follows that the disk $\tilde G_j =G_j\bigcup (x,s,t)\cup (t,s,y)$ belongs to $\cF_j^0$, which contradicts the maximality of $G_j$. Hence $\theta_s\geq \pi$.

Let now $s$ be a vertex of $\del D\cap g_1$ which does not belong to $g_0\cap g_1$ and let us show that the internal angle $\theta_s$ of $D$ at $s$ is at least $\pi$ as well. 
As $g_1$ is a simplicial geodesic of $X$ one has $\theta_s\geq 2\pi/3$, so we assume that $\theta_s=2\pi/3$ and argue towards a contradiction. Denote by $(x,s,t)$ and $(t,s,y)$ the corresponding triangles in $D$ so that  $[x,s]$ and $[s,y]$ are two consecutive edges of $g_1$ because $s\notin g_0\cap g_1$. Up to permuting $x$ and $y$ one can write  $g_1=h_0\cup [x,s]\cup [s,y]\cup h_1$ where $h_0$ is a simplicial geodesic in $X$ from $A_0$ to $x$ and $h_1$ is a simplicial geodesic in $X$ from $y$ to $A_n$.
Let $\tilde g_1=h_0\cup [x,t]\cup [t,y]\cup h_1$. As 
\[
\ell(\tilde g_1)\leq \ell(h_0)+\ell(h_1)+2=\ell(g_1)
\]
the  path $\tilde g_1$ is a simplicial geodesic from $A_0$ to $A_n$ in $X$ (where $\ell(g)$ denotes the simplicial length of $g$). Let $\tilde D=D\backslash \{(x,s,t)\cup (t,s,y)\}$ and let $\tilde \cD$ be the union of $\tilde D$ and the disks in $\cD$ which are distinct from $D$. Then $(g_0,\tilde g_1,\tilde \cD)$ is an element of $\cA$ so by minimality of $(g_0,g_1,\cD)$ we get $\tilde g_1\subset C$ and as the interiors of disks in $\cD$ are disjoint from $C$ it follows that $g_0=\tilde g_1$. However this implies that the point $t\in g_0$ has an internal angle in $D$ of $4\pi/3$, which contradicts what was established in the previous paragraph. So $\theta_s\geq \pi$.

It follows that the disk $D$ has internal angles at every point $s\in \del D$  at least $\pi$, except perhaps at the two points $\del D\cap g_0\cap g_1$. But this is a contradiction, no such a disk can exist in a CAT(0) space. Thus $\cH_1$ is empty and it follows that 
\[
\conv(\gamma)\subset C.
\]
Now it easy to show that for any three vertices of a triangle in $C$ there are geodesics of $C$ from $A_0$ to $A_n$ which contains these vertices (it is sufficient to prove the assertion for the flat $G_j$, which is easy). Moreover a geodesic of $X$ from $A_0$ to $A_n$, say of length $\ell_0$, is included in $C$ by the above and so is a geodesic of $C$. Thus all geodesic of $C$ have length $\ell_0$. It follows that $\conv(\gamma)= C$ and the lemma is proved.
\end{proof}

\begin{lemma}\label{emboit}
Let  $\gamma$ be a geodesic segment between two vertices $A$ and $B$ of $X$ and let $I$ be a vertex of $\conv(\gamma)$. Then the simplicial convex closure of the geodesic segment $[A,I]$ is included in $\conv(\gamma)$ and the union of any two simplicial geodesics in $X$ from $A$ to $I$ and  $I$ to $B$ respectively is a simplicial geodesic in $X$ from $A$ to $B$.  
\end{lemma}

\begin{proof}
By definition of the convex closure $I$ belongs to a simplicial geodesic $g$ of $X$ from $A$ to $B$, so $g=g_0\cup g_1$ where $g_0$ and $g_1$ are simplicial geodesic of $X$ from $A$ to $I$ and $I$ to $B$ respectively. Now for any geodesic $g'_0$ and $g_1'$ of $X$ from $A$ to $I$ and $I$ to $B$ respectively we have $\ell(g_0'\cup g_1')\leq \ell(g_0)+\ell(g_1)=\ell(g)$ so  $g'_0\cup g_1'$ is a simplicial geodesic of $X$ from $A$ to $B$ and $g'_0,g_1'\subset \conv(\gamma)$.
\end{proof}

Recall that a subset $D$ of $X$ is said to be convex if for any two points $A,B\in D$ the geodesic segment  $[A,B]$ from $A$ to $B$ is included in $D$ (see \cite{BH}).

\begin{lemma}\label{convex}
Let  $\gamma$ be a geodesic segment between two vertices of $X$. Then $\conv(\gamma)$ is convex.
\end{lemma}

\begin{proof}
Let $A,B$ be two points of $\conv(\gamma)$ and suppose that the geodesic segment $[A,B]$ is not included in $\conv(\gamma)$. Denote by $A'$ 
 the closest point from $A$ in $]A,B]$ which does not belong to the interior of $\conv(\gamma)$, and denote by $B'$ the closest point from $A'$ in $]A',B]$ which belongs to $\conv(\gamma)$. Note that both $A'$ and $B'$ are on the 1-skeleton of $X$  so we can choose a simplicial geodesic path $g$ from $A'$ to $B'$ inside $\conv(\gamma)$.
Let $D$ be the unique disk in $X$ with boundary $g\cup [A',B']$ so that, as in Lemma \ref{conv-triangles},  the internal angle of $D$ at every point of $g$ distinct from $A'$ and $B'$ is at least $\pi$. Since $[A',B']$ is a CAT(0)-geodesic the internal angle of $D$ at every point of $]A',B'[$ is at least $\pi$ as well so we get a contradiction.  Thus $[A,B]\subset \conv(\gamma)$.
\end{proof}

\begin{lemma}\label{equi}
Let  $A,B,C$ be three vertices of $X$ and let $\Delta$ be the geodesic triangle of $X$ with vertices $A$, $B$ and $C$. Assume that the angle of $\Delta$ at  $A,B$ and $C$ are non zero and denote by $D$ the   unique closed topological disc in $X$ whose boundary is $\Delta$.  Let $S_\Delta$ be the union of the simplicial convex closures of the segments $[A,B]$, $[B,C]$ and $[A,C]$. Then   the following assertions holds. 
\begin{enumerate}
\item If $D\backslash S_\Delta$ is empty  there is a point $I$ in $D$ and simplicial segments $\gamma_{AI},\gamma_{BI}$ and $\gamma_{CI}$ from $A$ to $I$, $B$ to $I$ and $C$ to $I$ respectively such that the simplicial paths $\gamma_{AI}\cup \gamma_{BI}$, $\gamma_{BI}\cup \gamma_{CI}$ and $\gamma_{AI}\cup \gamma_{CI}$ are simplicial geodesic segments of $X$ from $A$ to $B$, from $B$ to $C$ and from $A$ to $C$ respectively.  
\item If $D\backslash S_\Delta$ is not empty then  there exist a non empty flat equilateral triangle $T=(A', B', C')$ in $D$ whose edges $\gamma_{A'B'}$, $\gamma_{B'C'}$ and $\gamma_{A'C'}$ are singular simplicial geodesic between three distinct vertices  $A'$, $B'$ and $C'$ of $D$, and three simplicial geodesic segments  $\gamma_{AA'},\gamma_{BB'}$ and $\gamma_{CC'}$  in $X$ such that the simplicial paths $\gamma_{AA'}\cup \gamma_{A'B'}\cup\gamma_{B'B}$, $\gamma_{BB'}\cup \gamma_{B'C'}\cup\gamma_{C'C}$ and $\gamma_{AA'}\cup \gamma_{A'C'}\cup\gamma_{C'C}$ are simplicial geodesic segments in $X$ from $A$ to $B$, from $B$ to $C$ and from $A$ to $C$ respectively.
\end{enumerate}
\end{lemma}

\begin{proof}
Let $S_{AB}$, $S_{BC}$ and $S_{AC}$ be the simplicial convex closures of $[A,B]$, $[B,C]$ and $[A,C]$, so $S_\Delta=S_{AB}\cup S_{BC}\cup S_{AC}$, and let $n_{AB},n_{BC}$ and $n_{AC}$ be their respective nerves. 
 Assume first that $D\backslash S_\Delta$ is empty and let us distinguish the following two cases. 

 Suppose first that  the disk  $D$ is included in the union of the simplicial convex closures of two of the edges of its boundary, say $D\subset S_{AB} \cup S_{AC}$. In particular $[B,C]\subset S_{AB} \cup S_{AC}$. As simplicial convex closures are unions of simplexes of $X$, the nerve $n_{BC}$ of $S_{BC}$ is included in $S_{AB} \cup S_{AC}$ and $n_{BC}\cap  S_{AB} \cap S_{AC}$ is a non empty union of simplexes. Therefore we can find a vertex $I$ in $n_{BC}$ which belongs to  $S_{AB} \cap S_{AC}$. Choose two simplicial geodesics  $\gamma_{BI}$ and $\gamma_{CI}$ in  (the boundary of)  $n_{BC}$ from $B$ to $I$ and $C$ to $I$ respectively. It  follows from Lemma \ref{conv-triangles} that $\gamma_{BI}$ and $\gamma_{CI}$ are simplicial geodesic in $X$  from $B$ to $I$ and $C$ to $I$ respectively. Consider a simplicial geodesic $\gamma_{AI}$ in $X$ from  $A$ to $I$. As  $I$ belongs to both $S_{AB}$ and $S_{AC}$ one has that $\gamma_{AI}\subset S_{AB} \cap S_{AC}$ and that $\gamma_{AI}\cup \gamma_{BI}$ and $\gamma_{BI}\cup \gamma_{CI}$ are simplicial geodesic from $A$ to $B$ and $B$ to $C$  respectively, as  follows from Lemma \ref{emboit}. This shows that (1) holds in that case.

 Suppose now that  $D$ is not included in the union of the simplicial convex closure of any two  edges of its boundary. In particular the closure $D'$ of $D\backslash (S_{AB} \cup S_{AC})$ has non empty interior. We claim that there is a point $I$ in the interior of $D$ which belongs to  $\del D'\cap S_{AB}\cap S_{AC}$. 
Indeed assume there is no such a point. Then every point of $P=\del D'\backslash [B,C]$ is either in $S_{AB}\backslash S_{AC}$ or in  $S_{AC}\backslash S_{AB}$.
This implies that every connected component of $P$ is entirely included in  $S_{AB}\backslash S_{AC}$ or in  $S_{AC}\backslash S_{AB}$ (as $P=(P\cap S_{AB})\amalg (P\cap S_{AC})$ and both sets are closed in $P$). It follows that $P$ is connected,  as any disk whose boundary is included in a given convex closure, say $S_{AB}$, is  actually entirely included  in $S_{AB}$. So suppose for instance that  $P$ is included in $S_{AB}\backslash S_{AC}$. Then the closure of $P$ is included in $S_{AB}$ and so the geodesic $D'\cap [B,C]$ is included in $S_{AB}$ as well by Lemma \ref{convex}. Hence $\del D'\subset S_{AB}$ which implies $D'\subset S_{AB}$ and contradicts the fact that  $D'$ is non empty. This proves the claim. So let $I\in \del D'\cap S_{AB}\cap S_{AC}$ 
(in fact this point is unique). Since $D'$ is non empty and $D\backslash S_\Delta$ is empty, one has  $I\in S_{BC}$. Now choose simplicial geodesic $\gamma_{AI},\gamma_{BI}$ and $\gamma_{CI}$ of $X$, from $A$ to $I$, $B$ to $I$ and $C$ to $I$ respectively.  Lemma \ref{emboit} readily implies that the required conditions in (1) are satisfied.  
This conclude the proof of (1). 

So for the remaining part of the proof, we assume that $D\backslash S_\Delta\neq \emptyset$. 

Denote by $n_\Delta=n_{AB}\cup n_{BC}\cup n_{AC}$ the  nerve of $S_\Delta$.  Note that every triangle $t$ of $n_{AB}$ is divided into two parts by the geodesic $[A,B]$ and that exactly one of this part  has a non trivial intersection $t^+$ with the interior of $D$, and similarly for $n_{BC}$  and $n_{AC}$. We write $n_{\Delta}^0$ for the subset of $n_{\Delta}$ defined by
\[
n_{\Delta}^0=\bigcup_{s\subset n_{\Delta}} s~\cup~ \bigcup_{t\subset n_{\Delta}} \overline{t^+} 
\]
where $s$ runs over the singular segments of $n_{\Delta}$ and $t$ over all its triangles, and $\overline{t^+}$ is the closure of $t^+$.  Then the closure $D_0$  of $D\backslash n_{\Delta}$  is a closed topological disk included in $D$ which is a non empty union of closed triangles and whose boundary $\del D_0$ consists of the points in the boundary of $n_{\Delta}^0$ which are  in the interior of $D$, and the singular segments  of $n_{\Delta}^0$.

Let us construct by induction a  finite decreasing sequence $D_0\supset D_1\supset \ldots\supset D_m$, $m\in \NI$,  of closed  topological disks which are non empty unions of closed triangles such that  the two following conditions, henceforth  referred to as Property $P_k$, are satisfied for every non negative integer $k\leq m$:
\begin{enumerate}
\item[($P_k^0$)] the pairwise intersection of the simplicial segments 
\[
\gamma_{AB}^k=\del D_k\cap S_{AB},~~~~ \gamma_{BC}^k=\del D_k\cap S_{BC}~~~~\mathrm{and}~~~~\gamma_{AC}^k=\del D_k\cap S_{AC}
\]
 is reduced to a point,  say,
 \[ 
 \gamma_{AB}^k\cap \gamma_{AC}^k=A_k,~~ \gamma_{AB}^k\cap \gamma_{BC}^k=B_k~~\mathrm{and}~~ \gamma_{AC}^k\cap \gamma_{BC}^k=C_k,
 \]
\item[($P_k^1$)]  $\gamma_{AB}^k$, $\gamma_{BC}^k$ and $\gamma_{AC}^k$ are included into simplicial geodesic of $X$ from $A$ to $B$, from $B$ to $C$ and from $C$ to $A$ respectively. 
\end{enumerate}
and such that the internal angle of the disk $D_m$ at every point $s\in \del D_m$ distinct from $A_m, B_m$ and $C_m$  is at least $\pi$.

Note that Property $P_0$ is satisfied ($P_0^1$ is a consequence of Lemma \ref{conv-triangles}). Assume the construction has been done up to some non negative integer $k$. If the internal angle of the disk $D_k$ at every point $s\in \del D_k$ distinct from $A_k, B_k$ and $C_k$  is $\geq \pi$ we  set $m=k$ and stop the construction here. Otherwise there is an $s\in \del D_k$ distinct from $A_k, B_k$ and $C_k$ whose internal angle $\theta_s$ in $D_k$ is $\leq 2\pi/3$. Assume for instance that  $s\in \gamma_{AB}^k$  (the case of $\gamma_{BC}^k$ and $\gamma_{AC}^k$ being similar). By Property $P_k^1$  there exists a simplicial geodesic segment $g$ of $X$ from $A$ to $B$ which contains $\gamma_{AB}^k$. In particular $\theta_s\neq \pi/3$, and so $\theta_s=2\pi/3$. Denote by $(x,s,t)$ and $(y,s,t)$ the two corresponding triangles in $D_k$, where  the two edges $[x,s]$ and $[s,y]$ are included in   $\gamma_{AB}^k$.  Up to permuting $x$ and $y$ we can write $g=g_1\cup [x,s]\cup [s,y]\cup g_2$ so that $\ell(g)=\ell(g_1)+\ell(g_2)+2$. Consider the simplicial curve $g'=g_1\cup [x,t]\cup [t,y]\cup g_2$. Obviously $\ell(g')\leq\ell(g)$ so $g'$ is a simplicial geodesic of $X$ from $A$ to $B$. There are two cases.

\begin{itemize}
\item Suppose that $t$ belongs to the interior of $D_k$. Define
\[
D_{k+1} =D_{k} \backslash \{(x,s,t)\cup (y,s,t)\}  
\]
and $\gamma_{AB}^{k+1}$, $\gamma_{BC}^{k+1}$ and $\gamma_{AC}^{k+1}$ to be the intersection  with $\del D_{k+1}$ of $g'$, of  $\gamma_{BC}^k$ and of $\gamma_{AC}^k$
respectively. Then $D_k$ is a topological disk which is a non empty union of triangles (as  $D\backslash S_\Delta\neq \emptyset$) and property $P_{k+1}$ is easily seen to be satisfied and we can iterate the construction.

\item Suppose  that $t$ belongs to the boundary of $D_k$, say $t\in \gamma_{BC}^k$ for example. Denote by $\gamma_{AB}^{k+1}$ the union of the portion of $\gamma_{AB}^{k}$ between $A_k$ and $x$ and of the segment $[x,t]$, and by $\gamma_{BC}^{k+1}$ the portion of $\gamma_{BC}^{k}$ between $t$ and $C_k$. Set $\gamma_{AC}^{k+1}=\gamma_{AC}^k$. Then the three paths $\gamma_{AB}^{k+1}$, $\gamma_{BC}^{k+1}$ and $\gamma_{AC}^{k+1}$ satisfies Property $P^0_{k+1}$ with $A_{k+1}=A_k$, $B_{k+1}=t$ and $C_{k+1}=C_k$, and they bound a disk $D_{k+1}$ included in $D_k$. Since property $P_{k+1}^1$ is satisfied by construction, we only have to show that $D_{k+1}$ is non empty. Assume towards a contradiction that it is  empty. Then $t$ belongs to $\gamma_{AC}^k$ and hence to the intersection $S_{AB}\cap S_{BC}\cap S_{AC}$. Choose simplicial geodesics of $X$ from $t$ to $A$, $B$ and $C$, say to $\gamma_{At},\gamma_{Bt}$ and $\gamma_{Ct}$ respectively. We  have that $\gamma_{At}\cup \gamma_{Bt}\subset S_{AB}$ (Lemma \ref{emboit}) and so the disk $D_{AB}\subset D$ with boundaries $[A,B]$, $\gamma_{At}$ and $\gamma_{Bt}$ is included in $S_{AB}$. Arguing similarly for $[A,C]$ and $[B,C]$, we conclude that $D\backslash S_\Delta$ is empty, contrary to our standing assumption. Finally $D_{k+1}\neq \emptyset$, and this shows that we can iterate the construction in this case too. 
\end{itemize}

Hence there exists a disk $D_m$ satisfying Property $P_m$ and whose  internal angle of at every point $s\in \del D_m$ distinct from $A_m, B_m$ and $C_m$  is at least $\pi$.
Let us now prove that  $D_m$ defined above  satisfies the conditions in Item (2) of the lemma. We need to  show that $D_m$ is a flat equilateral triangle. 

Recall that for any open topological disk  $D$  in $X$ with piecewise linear topological boundary $\Delta$, one calls  \emph{geodesic curvature of $D$} at a point $s$ of $\Delta$  is the number $\kappa_s=\pi-\theta_s$ where $\theta_s$ is the internal angle of $D$ at $s$  (note that $\kappa_s$ is zero for all but a finite number of $s\in \Delta$). For a point  $x$ be a point in $D$, one calls  \emph{curvature  of $D$ at $x$} is defined as  $\delta_x=2\pi-\phi_x$, where $\phi_x$ the sum of the angles at $x$ in $D$, i.e.\ the angular length of the circle in the link $L_x$ of $x$ in $X$ defined by  $D$ (so the disk $D$ is flat at $x$ if and only if $\delta_x=0$).  

Denote by $T=\del D_m$ the triangle with vertex $A_m$, $B_m$ and $C_m$ and simplicial edges $\gamma_{AB}^m$, $\gamma_{BC}^m$ and $\gamma_{AC}^m$. By  construction a point $s\in \del D_m$ whose internal angle in $D_m$ is $<\pi$, if it exists, is necessarily one of the three vertices $A_m$, $B_m$ or $C_m$. This shows that the total geodesic curvature $\int_{\Delta_m} \kappa$ of the boundary $\Delta_m$ of $D_m$, that is the sum over all points $s\in \Delta_m$ distincts from $A_m$, $B_m$ and $C_m$ of $\kappa_s$, satisfies 
\[
\int_{\Delta_m} \kappa\leq 0.
\]
Applying the Gauss-Bonnet Formula (for domains with piecewise linear boundary) to the disk $D_m$ we get  
\[
(\pi-\theta_{A_m}) 
+(\pi-\theta_{B_m})+(\pi-\theta_{C_m})+~~\sum_{x\in D_m}  \delta _x~~=~~2\pi- \int_{\Delta_m} \kappa\geq 2\pi.
\] 
where $\theta_{A_m}$, $\theta_{B_m}$ and $\theta_{C_m}$ are the internal angle of $D_m$ at $A_m$, $B_m$ and $C_m$ respectively, and so
\[
\theta_{A_m} +\theta_{B_m} +\theta_{C_m}-\sum_{x\in D_m}  \delta _x\leq \pi.
\]
Now the fact that $D_m$ is simplicial (with non empty interior) implies that the values $\theta_{A_m}$, $\theta_{B_m}$ and  $\theta_{C_m}$ are at least $\pi/3$.
It follows that $\theta_{A_m}=\theta_{B_m}=\theta_{C_m}=\pi/3$ and that the disk $D_m$ is flat. Thus $T$ is a flat equilateral triangles with simplicial edges. Furthermore Property $P_m^1$  implies that $A_m$ belongs to $S_{AB}\cap S_{AC}$ and that $[A_m,B_m]$ is included in a simplicial geodesic in $X$ from $A$ to $B$. Fixing a simplicial geodesic $\gamma_{AA_m}$ in $X$ from $A$ to $A_m$ (which is included in $S_{AB}\cap S_{AC}$) and arguing similarly for $B_m$ and $C_m$, we conclude that all conditions in (2) hold. This proves the lemma.
\end{proof}

\newcommand{\equ}{\mathfrak{equ}}

We now prove Theorem \ref{th1}.

\begin{proof}[Proof of Theorem \ref{th1}]
Let $\G$ be a group acting properly and isometrically on  $X$ and let $A_0\in X$. Identify $\G$ to the orbit $\G A_0$ and consider the length $\ell$ induced by the 1-skeleton of $X$. We prove that $\G$ has property RD with respect to $\ell$ by using Section \ref{reminder} and the above lemmas.

Assume first the action of $\G$ to be free and transitive. For $z\in \G$ and $r\in \RI_+$ let $U_z^r$ the set triples $(a,u,b)$ with $z=bua$ such that the points $aA_0$, and $uaA_0$ belong to a simplicial geodesic from $A_0$ to $zA_0$ and such that the length of $[A_0,aA_0]$ and $[aA_0,uaA_0]$ is less than $r$. By Lemma \ref{conv-triangles} the family $U=(U_z^r)_{z\in \G, r\in \RI}$ has polynomial growth (in the sense of Definition \ref{HK}). 
Let $\equ$ be the family of equilateral triangles in $\G$ (i.e.\ the family of triangles $(x,y,z)$ in $\G$ such that $(A_0,xA_0,zA_0)$ is equilateral in $X$).  Lemma \ref{equi} shows that $\equ$ is a retract of the set of all triangles  of $\G$  along $U$ (see Definition \ref{rd-retr}, here $p_2(t)=t$). Indeed let $A$, $B$ and $C$ be three points of $X$ and denote by $ABC$ the corresponding geodesic triangle.  Note that $[A,B]\cap [A,C]$ is either reduced to $A$ or a geodesic segment of $X$. In the latter case we let  $A'$ be the extremity of this segment which is distinct from $A$. Otherwise set $A'=A$ and define similarly $B'$ and $C'$. If $A'=B'=C'$ then $ABC$ is a tripode triangle and it reduces to $A'$ along (any choice of) simplicial geodesics from $A'$ to $A$, $B$ and $C$.  If $A'$, $B'$ and $C'$ do not coincide then we get a geodesic triangle $A'B'C'$ whose angles at $A'$, $B'$ and $C'$ are non zero. Observe now that for any two vertices $A', B'$ of $X$ which are ordered $A\leq A'\leq B'\leq B$ on a given geodesic segment $[A,B]$ of $X$, the concatenation of any three simplicial geodesics from $A$ to $A'$,  $A'$ to $B'$, and $B'$ to $B$ respectively is a simplicial geodesic from $A$ to $B$. This together with Lemma \ref{equi} shows that $\equ$ is a retract of the set of all triangles  of $\G$  along $U$ (the metric requirements  
on the retraction being readily satisfied from description of simplicial convex closure in Lemma \ref{conv-triangles}).
Hence by Lemma \ref{rd-retr} property RD for $\G$ is equivalent to property RD for $\equ$.   Lemma \ref{rd-dual} then shows that property RD for $\equ$ is equivalent to property RD for the dual set $\equ^*$ (see  Section \ref{RRS-dual} and note that $\equ$ is obviously balanced). As all triangles in $\equ^*$ are tripods Lemma \ref{rd-retr} applies.   

As noted in \cite{RRS} the same proof works for free isometric action provided we replace everywhere `groups acting transitively' by `transitive groupoids' (see Subsection \ref{groupoids}). Next, as noted in \cite[Section 3.1] {Laf-rd}, the proof   works also for proper isometric actions too provided we replace the space $X$ by the disjoint union of all stabilizers $\G_A$, $A\in X$, of the action of $\G$ on $X$, where each $\G_A$ is endowed with the complete graph structure and there is an edge between $\tilde A\in \G_A$ and $\tilde B\in \G_B$ if and only if there is an edge between $A$ and $B$ in $X$. Then $\G$ acts freely isometrically on $\tilde X$ (which is quasi-isometric to $X$) and the required polynomial growth conditions all are satisfied for $\tilde X$ because  $\sup_{A\in X}\# \G_A<\infty$, as is easily seen.
\end{proof}




Let us conclude this section with the fact that polynomial branching in the sense of Definition \ref{def2} and geodesic polynomial branching in the sense of Definition \ref{polrk}  coincide in case of triangle groups, which  essentially follows from the above lemmas. 

\begin{proposition}\label{polrk-same}
Let $X$ be a triangle polyhedron. Then $X$ has polynomial branching in the sense of Definition \ref{def2} if and only the space $(X^{(0)},\ell)$  has  geodesic polynomial branching in the sense of Definition  \ref{polrk}, where $\ell$ is any length on the vertex set $X^{(0)}$ of $X$ induced by the 1-skeleton of $X$ from some base point $A_0\in X^{(0)}$ (i.e.\ $\ell(A)$ is the simplicial distance between $A_0$ and $A\in X^{(0)}$). 
\end{proposition}

\begin{proof}
Assume that $X$ doesn't have polynomial branching in the sense of definition \ref{def2} and let $(\gamma_{r_n})_{n\in \NI}$ be a sequence of simplicial geodesic segments in $X$ of length $r_n\in \RI_+$ such that  the growth of the number of flat equilateral triangle in $X$ with base $\gamma_{r_n}$ is faster than any given polynomial function of $r_n$. Then   the edges of every equilateral triangle $ABC$ of base $\gamma_{r_n}$ are simplicial geodesic segment, so for any fixed $\delta>0$, the set of $\delta$-geodesics of $ABC$ (for the simplicial structure) between the vertices $A$, $B$ and $C$ lies in the $\delta$-neighbourhood of the boundary of ABC.  Hence letting $T_n$ be the set of equilateral triangle of base $\gamma_{r_n}$, the set  of retractions of  triangles in $T_n$ along any set of $\delta$-geodesic as in Definition \ref{polrk} contains as many triangles  as $T_n$, up to a constant number depending only on $\delta$. This shows that $X$ does not have geodesic polynomial branching.  

Conversely assume that $X$ has polynomial branching in the sense of definition \ref{def2} and consider the retraction along the family $C$ as in the proof of Theorem \ref{th1} (which consists of simplicial geodesic of $X$, i.e.\  $\delta=0$). Then as shown in the lemmas above,  retractions along $C$ consist of flat simplicial equilateral triangles. In particular $X$ has geodesic polynomial branching (with  $\delta=0$) in the sense of Definition \ref{polrk}.    

Note that this proof works as well for geodesic subexponential branching and geodesic exponential branching.
\end{proof}

\section{Rank \sq}\label{class-rg74}

A triangle polyhedron is said to have \emph{order $q\in \NI^*$}  if each of its edges is contained into $q+1$ triangles. 
In this section $q=2$. A graph is \emph{ample} if the length of its smallest cycle is 6.
The following proposition is well-known (it is contained for example in the Foster census).

\begin{proposition}\label{trivalent} The smallest ample trivalent  graph has 14 vertices. Moreover there exists a unique ample symmetric trivalent graph with
\begin{enumerate}
\item  14 vertices: this is the incidence graph $L_2$ of the Fano plane $P^2(\FI_2)$. 
\item  16 vertices: this is $L_{\sqm}$ from the introduction (the Moebius--Kantor graph). 
\item  18 vertices: this is the Pappus graph $L_{3\over 2}$. 
\end{enumerate}
(Recall that the number of vertices of a trivalent graph is ${2\over 3}\times$ the number of its edges and so it is even.) Every ample trivalent graph with at most $16$ vertices is symmetric.
\end{proposition}

The graph $L_2$  is a spherical building. By  \cite{Tits74}  a triangle polyhedron of order 2 is an Euclidean building if and only if its link at each vertex is isomorphic to $L_2$, and there is a correspondence between  local rank 2, i.e.\ links being spherical buildings, and global rank 2, i.e.\ being an affine building.  We refer to \cite{these,rg2} for more details. 

In \cite{cras} the first-named author considered the graph $L_{3\over 2}$ of item (3) and constructed a compact complex $P$ of dimension 2 with one vertex  whose link  is isometric to $L_{3\over 2}$    (see Section 1 in \cite{cras}), and proved the Haagerup property for its fundamental group. The universal cover $\tilde P$ of $P$ has isolated flats and thus polynomial branching.

Proposition \ref{trivalent} shows that $L_{\sqm}$ is in some sense\footnote{In a later paper \cite{bp2010}, we introduced a general notion of local rank for CAT(0) polyhedra of dimension 2. It is based on a numerical invariant $\mathrm{rk}$ which has ample metric graphs as input, and returns the values 2, \sq, $3\over 2$ when applied to the above three graphs  (in the given order). So the numerical values agree. The graph $L_{3\over 2}$ satisfies the more restrictive property of being of rank \td\ (which is a local flat isolation property) as defined in the introduction of the present paper.} the `canonical' graph  which occupies an intermediate position between rank $3\over 2$ and rank 2.

\begin{definition} \label{comp74}
A triangle polyhedron is said to be of rank \sq\ if the link at each of its vertices is $L_{7/4}$. A triangle group is said to be of rank \sq\ if it admits a proper cocompact action on a  polyhedron of rank \sq. We call \emph{complex of rank \sq}\ a compact CW-complex of dimension 2 with equilateral faces and whose universal cover is polyhedron of rank \sq\ (for the usual piecewise linear metric). A complex of rank \sq\ is said to be \emph{orientable} if there is a coherent orientation of its faces.
\end{definition}

\begin{remark} To further interpolate the local rank between $3\over 2$, \sq, and 2, one can (for instance) increase the order $q$ (this will not be studied in the present paper), or mix together various links at vertices whenever possible (see Section \ref{meso}). 
\end{remark}

We now give a proof of the last assertion of Proposition \ref{trivalent}, which is what we need for the present paper and since it is not contained in the Foster census (the Foster census classifies cubic symmetric graphs of low order). Let $L$ be an ample trivalent graph, $n$ be the number of vertices and $m={3\over 2} n$ be the number of edges. The universal cover  $T$ of $L$ is the standard trivalent tree. Fix a vertex $*$ in $T$  and let $B_3$ be the open ball of center $*$ and radius 3 in $T$. The ampleness condition shows that $B_3$ is included in a (connected) fundamental domain $F$ of the action of the fundamental group $\G=\pi_1(L)$ on $T$.  Thus $m\geq 21$ and $n\geq 14$. If $n=14$, then the 12 vertices in the boundary $\del B_3$ of $B_3$ must be identified 3 by 3 without creating cycles of length $<6$. One readily checks that there is indeed a possibility (namely $L_2$) which is unique up to isomorphism. Assume that $n=16$. Then the ball $B_3$ contains all but 3 edges of $F$, so $B_3/\G$  has 15 vertices and $\del B_3/\G$ consists of 5 vertices: 2 of valence 3 (denoted $A$ and $B$) and three of valence 2 (in $B_3/\G$). Write $*'$ for the vertex of $F$ which is not in $B_3$ and call branch a connected component of $B_3\backslash \{*\}$. Then any two vertices $*_A$ and $*_B$ projecting down to $A$ and $B$ respectively cannot be at distance 2 in a same branch (otherwise identifications of the two other extremities of this branch with $*'$ under $\G$ would create a small cycle). Thus each branch must contain (points in the orbits of) $*_A$ and $*_B$ at a distance equaling 4. But then the remaining identifications from a branch to the other are uniquely determined up to isomorphism. The only possibility is the graph $L_\sqm$, which is indeed symmetric.

\subsection{Existence and classification results} 

In this section we  prove Theorem \ref{class} stated below. More precisely, we prove all the assertions of this theorem  except for the fact that there are \emph{at most} 13 complexes of rank \sq\ with one vertex, which is postponed until Subsection \ref{prog} (this is a computer assisted proof). In course of the proof we give  an explicit description of all rank \sq\ complexes appearing in item (1). The verification that these complexes are indeed of rank \sq\ will be straightforward by checking (a), (b) and (c) below. 

\begin{theorem}\label{class}
 Let $V$ be CW-complex of dimension 2 with triangular faces. Assume that 
\begin{enumerate}
\item[(a)] $V$ has 8 faces, 
\item[(b)]  each edge of $V$ is incident to 3 faces,  
\item[(c)] the link at each vertex of $V$ is ample.
\end{enumerate}
Then $V$ is a complex of rank \sq\ with one vertex. In particular $\pi_1(V)$ is a triangle group of rank \sq.
Moreover, 
\begin{enumerate}
\item there are precisely 13 orientable complexes of rank \sq, 
\item the homology group $H_1(V,\ZI)$  of such a $V$  can have rank 0 (i.e.\ can be torsion), rank 1 and rank 2. 
\end{enumerate}
In particular there exist groups of rank \sq\ which do not have Kazhdan's property T.
\end{theorem}

\begin{proof}
Let $V$ be CW-complex of dimension 2 with triangular faces. Let us first show that (a), (b) and (c) implies that $V$ is a complex of rank \sq\ with one vertex. Let $\{s_1,\ldots s_n\}$ be the vertex set of $V$. By (b) and (c) the link $L_{s_i}$ at each vertex $s_i$, $i=1\ldots n$, is an ample trivalent graph, so  Proposition \ref{trivalent}  implies that 
 $|L_{s_i}^{(0)}|\geq 14$  for every $i=1\ldots n$, where $L_{s_i}^{(0)}$ is the vertex set of $L_{s_i}$. By (a) one has 
 $|L_{s_i}^{(1)}|\leq 24$, denoting by $L_{s_i}^{(1)}$  the edge set of $L_{s_i}$. So Proposition \ref{trivalent} again implies that $L_{s_i}=L_2$ or $L_{s_i}=L_{7\over 4}$. This forces $n=1$ which in turn implies $L_{s_1}=L_{7\over 4}$ by (a) again.
 Thus $V$ is a complex of rank \sq\ with one vertex.

To construct such a complex we thus have to understand how  8 triangles can be glued together on a base point $*$ so that the link $L_*$ at $*$ is a trivalent ample graph. This  should be compared to the case of triangle buildings as derived in  \cite{CMSZ} and \cite{these}. 

So consider a bouquet $B_8$ of 8 oriented circles on which these triangles will be glued along their edges. Such a triangle $t$ will be denoted $[x,y,z]$ where $x,y,z$ are numbers in $\{1,\ldots, 8\}\cup \{1^-,\ldots, 8^-\}$ corresponding to the circles in $B_8$ on which its consecutive (for some fixed orientation of $t$) edges are attached. A minus sign occurs when the orientation of the circle in $B_8$ on which the edge of $t$ is glued is opposite to that of $t$. Thus a complex of rank \sq\ is entirely described by a list of 8 triples $[[x_1,y_1,z_1],\ldots [x_8,y_8,z_8]]$ with $x_i,y_i,z_i\in \{1,\ldots, 8\}\cup \{1^-,\ldots, 8^-\}$.   Furthermore,  condition (b) is equivalent to the fact that each number in $\{1,\ldots, 8\}$ appears exactly 3 times. This 8-list is called a presentation of the space $V$. 

Assume that $V$ is orientable. Then it has  a presentation as above with  $x_i,y_i,z_i\in \{1,\ldots, 8\}$. By the first part of Theorem \ref{class} the only condition to  check for $V$ to be a complex of rank \sq\ is that $L_*$ is an ample graph. This can be further simplified by the following lemma. 

\begin{lemma}\label{hgf}
If $L_*$ has no cycle of length 2 and 4, then $L_*$ is ample. 
\end{lemma}
\begin{proof}[Proof of Lemma \ref{hgf}]
Let $[[x_1,y_1,z_1],\ldots [x_8,y_8,z_8]]$ be a presentation of $V$ with $x_i$, $y_i$, $z_i \in \{1,\ldots, 8\}$. Fix some small $\varepsilon_0$ so that the link $L_*$ coincide with the sphere of radius $\varepsilon_0$ and center $*$ in $V$.  For every $i\in \{1,\ldots 8\}$ let $i^\flat$ and 
$i^\sharp$ be the two points (as order by the orientation in $B_8$) of the  $i$-th circle of $B_8$ at distance $\varepsilon_0$ from $*$. The edges in $L_*$ are then  $[x_i^\flat, y_i^\sharp]$ $[y_i^\flat, z_i^\sharp]$  and $[z_i^\flat, x_i^\sharp]$ for $i=1\ldots 8$. In particular there is no cycle of length 3 or 5.
\end{proof}

Thus, the classification of orientable complexes of rank \sq\ can done in two steps: 
\begin{enumerate}
\item[($\alpha$)]  list all admissible presentations, i.e.\ presentations for which $x_i,y_i,z_i\in \{1,\ldots, 8\}$ and each number in $\{1,\ldots, 8\}$ appears exactly 3 times, and 
\item[($\beta$)]  check in each case that the corresponding link has no cycle of length 2 and 4. 
\end{enumerate}
More details  concerning the implementation of this procedure will be given in Section \ref{prog}. 
After extensive computations we get 
 the following list of the 13 complexes of rank \sq\ announced in (1), together with their first homology groups and their fundamental groups.
They are coming in 6 classes, according to the number of \emph{adjacent identifications} in their link (equivalently the number of $[x,x,\cdot]$ in their presentation). 
The most symmetric ones (in particular $V_0$ all of whose identifications occur at distance 3 in $L_*$) can be found by hand from the link.  

\bigskip

\noindent \emph{Type I (no adjacent identification).}   There are four orientable complexes in this class:
\begin{align*}
V_0=[[1,2,6],[2,3,7],[3,4,8],[4,5,1],[5,6,2],[6,7,3],[7,8,4],[8,1,5]]\\
V^1_0=[[1, 2, 3], [1, 4, 5], [1, 6, 4], [2, 6, 8], [2, 8, 5], [3, 6, 7], [3, 7, 5], [4, 8, 7]]\\
V^2_0=[[1, 2, 3], [1, 4, 5], [1, 6, 7], [2, 4, 6], [2, 8, 5], [3, 6, 8], [3, 7, 5], [4, 8, 7]]\\
\check V_0^2=[[1, 2, 3], [1, 4, 5], [1, 6, 7], [2, 6, 4], [2, 8, 5], [3, 6, 8], [3, 7, 5], [4, 8, 7]]
\end{align*}
with respective first homology groups $H_1(~\cdot~ ,\ZI)$:
\[
\ZI/15\ZI,\ \ \ \ \ \ \ \ZI/3\ZI\times \ZI^2, \ \ \ \ \ \ \ (\ZI/3\ZI)^3, \ \ \mathrm{and} \ \ \ \ \ (\ZI/3\ZI)^3
\]
Note that $V_0$ is the only polyhedron satisfying the additional following condition: for any $x=1\ldots 8$, there is $y=1\ldots 8$ such that $[x,y,\cdot]$ and  $[y,x,\cdot]$ are faces. The fundamental group $\pi_1(V^1_0)$ is a group of rank \sq\  with two generators $s,t$ and two relations (where $t^s=s^{-1}ts$):
\[
\pi_1(V_0)=\langle s, t\mid  t  s  t^2  s^{-1}  t   =st^2st^2t^s,\   s  t^{-1}  s^{-1}  t^{-2}  s    = tst^2 t^sst\rangle.
 \]
In the other cases the fundamental group has 3 generators and 3 relations:
\begin{align*}
\pi_1(V^1_0)&=\langle u, v, w\mid u= v  w^{-1}  v^{-1}  u  w,\ 
u  v = wvwuw,\ 
  u  v  u  w  = wuvu\rangle.\\
\pi_1(V^2_0)&=\langle u, v, w\mid u  v = wvwuw,\ 
  w  u  = uvwv^2,\ 
  v   =uvw^{-1}uwu\rangle.\\
\pi_1(\check V^2_0) &=\langle u, v, w\mid  u  v = wvwuw,\ 
w  u  =vuvwv,\ 
    wu^2v=   v  u^{-1}  w\rangle.
\end{align*}

\noindent \emph{Type II (one adjacent identification).} A single orientable polyhedron  in this class:
\[
V_1=[[1, 1, 2], [1, 3, 4], [2, 5, 6], [2, 7, 8], [3, 5, 7], [3, 6, 5], [4, 6, 8], [4, 8, 7]]
\]
with $H_1(V_1,\ZI) =\ZI/3\ZI\times \ZI$, and 
\[
\pi_1(V_1)=\langle s,t\mid s^{4}  t  s^{-3}  t  s  =ts^2t,\  t = s^2  t  s^{-1}  t^{-1}  s^2  t^{-1}  s^{-2}  t^2  s^2\rangle.
\]

\noindent \emph{Type III (2 adjacent identifications).}  There are four orientable complexes in this class:
\begin{align*}
V_2^1=[[1, 1, 3],  [2, 2, 3], [1, 4, 5], [2, 7, 8], [3, 5, 7], [4, 6, 8], [4, 7, 6], [5, 8, 6]]\\
V_2^2=[[1, 1, 3],  [2, 2, 4], [3, 7, 4], [1, 4, 6], [2, 5, 3],  [5, 7, 8], [5, 8, 6], [6, 8, 7]]\\
V_2^3=[[1, 1, 3],  [2, 2, 4], [1, 5, 2], [3, 6, 4] ,  [3, 7, 6], [4, 6, 8], [5, 7, 8], [5, 8, 7]]\\
V_2^4=[[1, 1, 3],  [2, 2, 4], [1, 5, 2], [3, 6, 5], [3, 7, 8], [4, 5, 8], [4, 6, 7], [6, 8, 7]]
\end{align*}
with respective first homology groups:
\[
\ZI/24\ZI,\ \ \ \ \ \ \ \ (\ZI/3\ZI)^2,\ \ \ \ \ \ \ \ZI/3\ZI\times \ZI, \ \ \mathrm{and} \ \ \ \ \ \ZI/24\ZI
\]
and respective fundamental groups
\begin{align*}
\pi_1(V_2^1) =& \langle s,t\mid
s  t^2=    tst^2s^{-1}ts^2 ,\
    t  s^{-1}  t^2  = s^{3}ts^2t^{-2}st^{-1}sts^2  
\rangle  \\
\pi_1(V_2^2) =&\langle s,t\mid        s^{2}=  t^2 s^4  t^{-3}  s  t ,\
   st= t^2  s^{-2}  t^{-4}  s^3  
 \rangle \\
\pi_1(V_2^3) =&\langle s,t\mid           s^2t^3s^3t^2=   t^2  s^2,\
    t^2=s^2t^2s^2t^2s^{-2} ts
 \rangle\\ 
\pi_1(V_2^4)=&\langle s,t\mid               s^2t^2sts^{2}=  t^3  s  t,\
   s= t^3stst^{-1}s^2tst 
 \rangle
\end{align*}

\noindent \emph{Type IV (3 adjacent identifications).} There is a single orientable complex  in this class:
\[
V_3 =[[1, 1, 4], [2, 2, 4], [3, 3, 5], [1, 3, 6],  [2, 5, 7],  [4, 7, 8], [5, 8, 6], [6, 8, 7]]
\]
with $H_1(V_3,\ZI)=\ZI/6\ZI$ and $\pi_1(V_3)=\langle s,t\mid
    st^3st=  t^2  s  t  s^2,\
  s^2=  t^2  s  t  s^{-2}  t  s^{-1}  t^2  s  t^3  
\rangle.$

\noindent \emph{Type V (4 adjacent identifications).} Two orientable complexes in this class:
\begin{align*}
V_4^1=[[1, 1, 5], [2, 2, 5], [3, 3, 6], [4, 4, 6], [1, 3, 8],  [2, 7, 4],  [5, 8, 7], [6, 7, 8]]\\
V_4^2=[[1, 1, 5], [2, 2, 5], [3, 3, 6], [4, 4, 7], [1, 3, 8],  [2, 7, 6],  [4, 8, 6], [5, 8, 7]]
\end{align*}
with $H_1(V_4^1,\ZI)=(\ZI/2\ZI)^2\times  \ZI/12\ZI$ and $H_1(V_4^2,\ZI) =\ZI/66\ZI$, and 
\begin{align*}
\pi_1(V_4^1)=&\langle u,v,w\mid
    v^2 = u^{2},\
    w^{2} uw= u  w  u^2,\
    w^2vuwvu^3wu^{2} = e
\rangle\\
\pi_1(V_4^2)=&\langle s,t\mid
    t^2  s  t^3  s  t  s  t  s^2=e,\
    t^2  s^{-2}  t^{-1}  s^{-1}  t^2  =  s^{3}ts^2
\rangle 
\end{align*}

\noindent \emph{Type VI (5 adjacent identifications).} Only one complex in this class:
\begin{align*}
V_5=[[1,1,2],[3,3,2],[4,4,2],[5,5,6],[7,7,8],[1,8,6],[3,6,7],[5,8,4]]
\end{align*}
with $H_1(V_5,\ZI)=(\ZI/6\ZI)^2$, and 
\begin{align*}
\pi_1(V_5)=&\langle s,t\mid
    s^3t^{-1}s^2t^{-3}st^{-2}=e, s^{-1}t^2s^{-1}t^2s^{-2}t^{-2}s^{-2}t^{-2}=e
\rangle
\end{align*}

This, together with Subsection \ref{prog}, concludes the proof of Theorem \ref{class}. 
\end{proof}

\begin{remark}[compare \cite{FH}]   The first non zero eigenvalue of $L_{7\over 4}$ is 
\[
\lambda_1(L_{7\over 4})=1-\frac 1{\sqrt 3} = 0.42...<\frac 1 2.
\]
Indeed let $\sigma_i$, $i=0,\ldots, 7$ be the $i$-th matrix of cyclic permutation of the set of 8 elements ($\sigma_0=\id$).
The random walk operator $D$  on $L$ in the canonical basis of vertices has the form 
\[
{1\over 3}\left (\begin{matrix}0 & A \\ A^t & 0\end{matrix} \right )
\] 
where $A$ is the $8\times 8$ matrix defined by $A=\id +\sigma_2+\sigma_7$. One has $AA^t=A^tA=2\Id + P +\sigma_4$ where $P^2=8P$. Thus the eigenvalues of $AA^t$ are 1, 3 and 9 and those of $D$ are $\pm{1\over 3}, \pm{1\over \sqrt 3}$ and $\pm1$ (of order 3, 4 and 1 respectively).
\end{remark}

\begin{remark}\label{barv} Here is an example of a complex of rank \sq\ which is not orientable:
\[
\bar V =[[3,1^{-},2], [3,2^{-},4], [2,6,3^{-}],[5,1^{-},6^{-}],[7,4^{-},5],[8,6^{-},7],[5,8,7],[1,4,8^{-}]]
\] 
It has a torsion-free $H_1(\bar V,\ZI)=\ZI$ and 
\[
\pi_1(\bar V) = \langle s,t\mid s^2t=t^2  s^2  t^2  s^{-1}  t^{-1}  s t^{-1} s ,~~
    t=s t s^{-1}  t  s^{-1}  t^2  s^{-1}  t^{-2}  s  t^{-1} 
    s\rangle.
    \] 
\end{remark}


\subsection{Homogeneity and structure of local flats.}\label{localsq} We now study the local rank structure in rank \sq\ polyhedra. 
 
 Endow $L_\sqm$ with the uniform length (edges have length $1$). The following concept  is important to describe the local behavior of flats of rank \sq\ polyhedra. 

\begin{definition}
Let $\alpha,\beta$ be two vertices of $L_\sqm$ at distance $3$ in $L_\sqm$. The  couple $(\alpha,\beta)$ is said to be \emph{of type $3\over 2$} if there are exactly two distinct simplicial paths of length $3$ in $L_\sqm$ with extremities $\alpha$ and $\beta$. If not, then there are exactly three such paths, in which case we call $(\alpha,\beta)$  \emph{of type 2}.
\end{definition}

The bipartite structure of $L_{\sqm}$ gives a partition of its vertex set into two sets of cardinal 8. We call the vertices of the first set (resp. the second set)  of type 0 (resp. type 1). The following is straightforward. 

\begin{proposition}
Let $\alpha$ be a vertex of $L_\sqm$. There are 5 vertices at distance $3$  from $\alpha$ in $L_\sqm$ (if $\alpha$ is of type $i=0,1$, these are the 5 vertices of $L_{\sqm}$ of type $1-i$ which are not adjacent to $\alpha$).  Three of them are vertices $\beta$ such that the couple $(\alpha,\beta)$ is of type $3\over 2$. For the two others, $(\alpha,\beta)$ is of type  2. The diameter of $L_\sqm$ is $4$ and there is a unique vertex of $L_\sqm$ at distance $4$ from $\alpha$.   
\end{proposition}

Let $G=\Aut(L_\sqm)$ be the automorphism group of $L_\sqm$.  

\begin{proposition}\label{autL}
The group $G$ is transitive on the tripods of $L_{\sqm}$. The subgroup of $G$ fixing pointwise  a tripod in $L_{\sqm}$ is trivial.  Its stabilizer is isomorphic to $\ZI/3\ZI\rtimes \ZI/2\ZI$. In particular $|G|=96$.
\end{proposition} 

It is interesting to compare this proposition to the proof of property T for (some) triangle buildings in \cite{Cart}. One can see that it is the lack of transitivity of stabilizers of tripods which explains that their proof doesn't apply to the present situation (what we already know from Theorem \ref{class}). 



\begin{proof}[Proof of Proposition \ref{autL}]  
Let $G_0$ the subgroup of $G$ fixing pointwise the tripod of $L_{\sqm}$ at a vertex $\alpha$. Then the vertex $\beta$ of $L_\sqm$ at distance 4 from $\alpha$ is fixed by $G_0$. Let $S$ be the complement of the (open) tripod of $\alpha$ in $L_\sqm$.  This is a connected graph on which $G_0$ acts with at least 4 fixed points. From this one easily infers that the tripod of $\beta$ is fixed by $G_0$. Thus $G_0$ is trivial. One can then check that the stabilizer of a tripod is isomorphic to $\ZI/3\ZI\rtimes \ZI/2\ZI$. Furthermore it readily seen on Figure \ref{fig1} that $G$ contains the group of dihedral symmetries (of order 16), which respects the type of vertices, as well as reflections which exchange the type of vertices. It follows that $G$ is simply transitive on the tripods (and so  $|G|=6\times 16=2^5\times 3$).
\end{proof}

Let us describe in more details the graph $L_{\sqm}$ from the local rank point of view.
The signification of the following result is that, already at the local level in any polyhedron $X$ of rank \sq, the link $L_\sqm$  \emph{imposes certain directions of (non-)branchings} in $X$. This has to be compared to Section \ref{meso}.

\begin{proposition}
Let $\Pi$ be a 6-cycle in $L_\sqm$. Then one of the 3 couples $(\alpha,\beta)$ of points at distance 3 in $\Pi$ is of type $3\over 2$, and the two others are of type 2. 
The group $G$ is transitive on the 6-cycles of $L_\sqm$ and has two orbits on the flags $f\subset \Pi$ where $f$ is an edge of $L_\sqm$ and $\Pi$ a 6-cycle.
\end{proposition}

\begin{proof} One easily checks on Figure \ref{fig1} that there exists a 6-cycle $\Pi_0$ which satisfies the assertion of the Proposition, i.e.\ one of the 3 couples  of points at distance 3 in $\Pi$, say $(\alpha_0,\beta_0)$ is of type $3\over 2$, and the two others are of type 2. In fact one can further assume   that there is a vertex $\delta_0$ at distance 2 from $\alpha_0$ on $\Pi_0$ such that, for the unique vertex $\beta_0'$ at distance 1 from $\delta_0$ which does not belong to $\Pi_0$, the couple $(\alpha_0,\beta_0')$ is of type 2. 

Proposition \ref{autL} implies that $G$ is transitive on the flags $A\subset \gamma$ where $\gamma$ is a (simplicial) path of length 2 and $A$ is an extremity of $\gamma$.  Let $\Pi$ be a 6-cycle, $(\alpha,\beta)$ be two vertex at distance 3 in $\Pi$, and let $\delta$ be a point at distance 2 from $\alpha$ on $\Pi$. Then there exists $g\in G$ such that $g(\alpha)=\alpha_0$  and $g(\delta)=\delta_0$. Then $g(\Pi)$ is a 6-cycle and we have either $g(\beta)=\beta_0$ or $g(\beta)=\beta_0'$.

 Thus there are 3 possibilities for  $g(\Pi)$ (all containing $\alpha_0$ and $\delta_0$). 
 One readily checks that these three cycles satisfies the first assertion of the proposition. Then, choosing in $g(\Pi)$ the unique couple $(\alpha_1,\beta_1)$ of points at distance 3 of type $3\over 2$, and mapping the couple $(\alpha_1,\delta_1)$, where $\delta_1$ is at distance 2 from $\alpha_1$ in $g(\Pi)$, to the couple $(\alpha_0,\delta_0)$, one sees that there is $h\in G$ such that $hg(\Pi)=\Pi_0$. 
 
Finally the fact that $G$ has two orbits on the flags $f\subset \Pi$ where $f$ is an edge of $L_\sqm$ and $\Pi$ a 6-cycle comes from the two possibilities for position of $f$ relatively to the couple of rank $3\over 2$ in $\Pi$.

This proves the proposition.     
\end{proof}

\subsection{On the asymptotic structure of flats}\label{asymsq}

We now study how the local analysis of the previous subsection `integrates' to polyhedra of rank \sq. 

The proof of the following proposition relies on techniques developed in the proof of theorem \ref{th7} and  is deferred to the end of Section \ref{pfth7}.

\begin{proposition}\label{altern}
Let $V$ be a complex of rank \sq\ and $\G=\pi_1(V)$ be the fundamental group of $V$. Then for any copy of 
the free abelian group $\ZI^2$ in $\G$, there is a $\gamma\in \G$ such that the pairwise intersection of the subgroups  $\gamma^n \ZI^2\gamma^{-n}$, $n\in \ZI$,  is reduced to the identity.
\end{proposition}

This proposition expresses the following dichotomy:  either there is no copy of $\ZI^2$ in $\G$, or there is abundance of $\ZI^2$ all of whose copies sit `weakly malnormaly' in  $\G$, the latter being a strong structural property of $\G$.  
We don't know if there are always copies of $\ZI^2$ in groups of rank \sq. In fact this is precisely what first led us to study polyhedra of rank \sq.  

More precisely, one of our motivations for studying these  polyhedra was the following well-known open problem in geometric group theory (see e.g.\ Question 1.1 in \cite{Best}, that we formulate here in the non-positively curved case and  dimension 2, as in the paragraph following Q 1.1. in \cite{Best}).

\begin{question}\label{Qbest}
Let $\G$ be a countable group admitting a non-positively curved finite $K(\G,1)$ of dimension 2. If the universal cover of $K(\G,1)$ contains a flat, does $\G$ contains $\ZI^2$?   
\end{question}

This question is especially intriguing for  polyhedra of rank \sq\ because of their local structure (as described in Subsection \ref{localsq}).  This  motivated our classification of complex of rank \sq\  in Theorem \ref{class}, where our objective was to study  the \emph{most symmetric} cases, i.e., polyhedra that are transitive on vertices. 
Relying upon this classification we now clarify the issue of Question \ref{Qbest} in this particular case.

\begin{proposition}\label{expsq}
Let $V$ be one of the 13 orientable complexes of rank \sq\ of with one vertex (see Theorem \ref{class}) and $\G=\pi_1(V)$ be the fundamental group of $V$. Then $\G$ contains copies of $\ZI^2$. 
\end{proposition}

\begin{proof}
We shall give full details for one of these complexes  (we chose  $V_0^2$) so as to present the arguments. The other cases can be derived similarly. 

Recall that $V_0^2$ admits the presentation
\[
V^2_0=[[1, 2, 3], [1, 4, 5], [1, 6, 7], [2, 4, 6], [2, 8, 5], [3, 6, 8], [3, 7, 5], [4, 8, 7]].
\]  
Let $a_1,\ldots,a_8$ be the generators $\G_0^2=\pi_1(V_0^2)$ of corresponding to $1,\ldots, 8$. Then one readily checks that the two elements $x=a_1^3$ and $y=a_3a_4$ commute in $\G_0^2$. Hence they generate a subgroup $\Lambda$ isomorphic  to $\ZI^2$.  
\end{proof}

We round up this section by providing some more information on the question of flat branching for groups of rank \sq, by showing that the complex $\tilde V_0^2$ considered in the proof above has exponential flat branching growth. This will be done  by exhibiting free  transverse semigroups. Keeping the notation of the proof, 
 the fundamental domain $P$ of the action of $\Lambda$ on its corresponding flat in $\tilde V$ contains 12 triangles which are respectively (from left to right and bottom to top after a suitable embedding of $P$ in $\RI^2$):
\begin{align*}
[1, 2, 3], [2, 8, 5], [4, 8, 7],  [1, 6, 7]\\
[1, 4, 5], [4, 8, 7],  [3, 6, 8], [1, 2, 3]\\
[1, 6, 7], [3, 6, 8], [2, 8, 5], [1, 4, 5]
\end{align*}
In fact, one can use this flat to show that $\tilde V_0^2$ has exponential branching (as mentioned in the introduction, this holds for several other examples of rank \sq\ polyhedra). For each vertex $A\in \tilde V_0^2$ consider the flat parallelogram $P_A$ in $\tilde V_0^2$ associated to $x$ and $y$ and let $Y=\cup_{A\in  \tilde V_0^2} P_A$. Then $Y$ is a union of flats and from the description of $P$ above, and the definition of exponential branching, it is not hard to check that it is enough to prove that the semi-group of $\G$ generated by the three elements 
\[
u=a_3a_4,~~ v=a_5a_6,~~\mathrm{and} ~ w=a_7a_2,
\] 
has exponential growth. Write $u'=vwu$. One easily check  on the presentation of $V_0^2$ that the for any vertex $A_0\in \tilde V_0^2$ the three points $A_0$, $uA_0$ and $u'uA_0$ (resp. $A_0$, $u'A_0$ and $uu'A_0$) are on geodesic of $\tilde V_0^2$. But this implies that the semi-group generated by $u$ and $u'$ is free in $\G$. Hence the (semi-)group of $\G$ generated by $u$, $v$ and $w$ has exponential growth. Note that this argument relies on the specific structure of $\tilde V_0^2$, and it is unclear what other groups of rank \sq\ (infinitely many?) satisfy the exponential  branching property.

\subsection{End of the proof of Theorem \ref{th4}} \label{prog}

This subsection is devoted to the (computer assisted) proof that there are \emph{at most} 13  orientable compact complexes of rank \sq\ with one vertex, as asserted in Theorem \ref{class}.  For us the program below was primary used to obtain a representative list of examples of polyhedra of rank \sq, beyond the few ones we found by hand.  
    
    We   only explain below what  procedures are  relevant  to understand the source code and to check that its mathematical part was correctly implemented. Other  procedures (e.g.\ displaying, sorting, memory management,...)  are routine and omitted. 
Recall from the proof of Theorem \ref{class} that  polyhedra of rank \sq\ can be represented by a list of numbers. 
At the mathematical level two procedures are important: 
  \begin{enumerate}
  \item[(a)]  the  test that  a list representing a polyhedron has the correct link, i.e., $L_{\sqm}$ (see \verb1test_link1),
  \item[(b)] the iteration process that enumerates all possible list representing orientable polyhedra of rank \sq\ with 1 vertex (see \verb2case_12).
  \end{enumerate}    
Let us briefly explain how they  are inserted in the body of the program (the code reproduced below is written in C++). We first need a class, called \verb1polyhed1, which describes a generic orientable polyhedron of rank \sq\ with 1 vertex:

\smallskip

 class \verb1polyhed1\{
  private:  char index[8][3];

 \hspace{.4cm} public:
     polyhed(); polyhed(char [][3], char);$\sim$polyhed();char length;\

 \hspace{.8cm}     void \verb1add1 (char, char, char); void \verb1display1();
     void \verb1reorder1();

 \hspace{.8cm}       char \verb1compare1(polyhed);        bool \verb1edge1(char, char );
     bool \verb1two_path1(char, char );

 \hspace{.8cm}     bool \verb1test_link1(char, char ,char );
 \};

\smallskip

The function  \verb1add1 add the face (char, char, char) to the given polyhedron (represented by the $8\times 3$ array \verb1index1)  and increments its length \verb1length1 of 1.
The function \verb1reorder1  reorders  \verb1index1 in a  canonical way (i.e.\ removes the ambiguities arising in the process of coding of the actual polyhedra as an array of number).  The function \verb1compare1 is a basic comparison of two polyhedra.

The function \verb1edge1 tests, given two numbers $(k,l)$ between 1 and 8,  whether $(i^\flat,k^\sharp)$ is already an edge of the link of the polyhedron $P$ (see the proof of Theorem \ref{class} for the notations $^\flat\, ^\sharp$ ). 
Similarly \verb1two_path1 tests, given $(k,l)$,  whether there is a length 2 path in the link of $P$ between $k^\flat$ and $l^\sharp$ or $k^\sharp$ and $l^\flat$.
These functions are straightforward to implement. They are used in the last function \verb1test_link1, which  tests whether adding the  face $(i,j,k)$ to a polyhedron $P$ (with \verb1length1$< 8$ faces) creates an admissible portion of the link:

\smallskip

 bool polyhed::\verb1test_link1(char i, char j, char k)\{ 
  if (i==j and j==k) return true;

\hspace{.1cm}  if (\verb1edge1(i, j) or \verb1edge1(j, k) or \verb1edge1(k, i)) return true;  

\hspace{.1cm}  if (i==j)\{if (\verb1edge1(k, k) or \verb1two_path1(i, k)) return true;\};  

\hspace{.1cm}  if (j==k)\{if (\verb1edge1(i, i) or \verb1two_path1(i, j)) return true;\};

\hspace{.1cm}  if (k==i)\{if (\verb1edge1(j, j) or \verb1two_path1(i, j)) return true;\};

\hspace{.1cm}  if (\verb1edge1(i, k) and \verb1edge1(j, j)) return true; 

\hspace{.1cm}  if (\verb1edge1(j, i) and \verb1edge1(k, k)) return true;

\hspace{.1cm}  if (\verb1edge1(k, j) and \verb1edge1(i, i)) return true;

\hspace{.1cm}  for (int a=0; a$<$length; a$++$)\{ 

  \hspace{1cm} if (\verb1edge1(index[a][0], j) and \verb1edge1(i, index[a][1])) return true;

  \hspace{1cm}   if (\verb1edge1(index[a][1], j) and \verb1edge1(i, index[a][2])) return true;

  \hspace{1cm}   if (\verb1edge1(index[a][2], j) and \verb1edge1(i, index[a][0])) return true;

  \hspace{1cm}   if (\verb1edge1(index[a][0], k) and \verb1edge1(j, index[a][1])) return true;

  \hspace{1cm}   if (\verb1edge1(index[a][1], k) and \verb1edge1(j, index[a][2])) return true;

   \hspace{1cm}  if (\verb1edge1(index[a][2], k) and \verb1edge1(j, index[a][0])) return true;

  \hspace{1cm}   if (\verb1edge1(index[a][0], i) and \verb1edge1(k, index[a][1])) return true;

   \hspace{1cm}  if (\verb1edge1(index[a][1], i) and \verb1edge1(k, index[a][2])) return true;

   \hspace{1cm}  if (\verb1edge1(index[a][2], i) and \verb1edge1(k, index[a][0])) return true;\}
   
\hspace{.1cm} return false;\}

\smallskip

Next one needs a class \verb1list1 to record the collected polyhedra. The function \verb1add1 in this class,  given some polyhedron, reorders it in a `canonical' order, compares it to the elements already present in the list, and adds a copy of it to the list if it was found to be new.

\smallskip

 class \verb1list1\{
 private:    long int max; polyhed **K;

  \hspace{.7cm} public:    List();$\sim$List();
     void \verb1display1();
     void \verb1add1 (polyhed *);
     \};

\smallskip

Finally we need to enumerate  all possible 8-tuples presentations of polyhedra. 
In fact doing this enumeration directly would have been too long, so we divided it  into the following  several subcases which  take  \emph{a priori} into account  (some of) the ambiguity in the representation of a polyhedron as a list. 

We first treat the situation where there exists a face, say $(1,2,3)$, which corresponds to pairwise non adjacent edges in the link $L$. Here, we claim that it is sufficient to enumerate all lists of the following form to exhaust all  the existing polyhedra:

Case 1:
$[(1,2,3),(4,4,5), (1,.,.),(2,.,.),(3,.,.),(1,.,.),(2,.,.),(3,.,.)]$

Case 2:
$[(1,2,3),(4,5,6), (1,.,.),(2,.,.),(3,.,.),(1,.,.),(2,.,.),(3,.,.)]$

Case 3:
$[(1,2,3),(1,3,2), (1,4,.),(2,.,.),(3,.,.),(.,.,.),(.,.,.),(.,.,.)]$

Case 4:
$[(1,2,3),(1,3,4), (3,5,.),(1,.,.),(2,.,.),(2,.,.),(.,.,.),(.,.,.)]$

Case 5:
$[(1,2,3),(1,3,4), (3,5,.),(2,1,.),(2,.,.),(.,.,.),(.,.,.),(.,.,.)]$

Case 6:
$[(1,2,3),(1,3,4), (2,1,.),(3,2,.),(.,.,.),(.,.,.),(.,.,.),(.,.,.)]$

The proof is easy: once the face $(1,2,3)$ is chosen, a new face may be disjoint from it (case 1 and 2) or attached to it (case 3$\sim$6). The former involves a new letter 4 which belongs to a face that may (case 1) or may not (case 2) have an adjacent identification. Case 3, which treats the situation where both $(1,2,3)$ and $(1,3,2)$ is a face, could have been divided into further subcases, thus further reducing the number of missing numbers for these cases, but this was not necessary. Finally, Cases $4\sim 6$ correspond to the last situation where $(1,2,3)$ and $(1,3,4)$ are faces.

Here is the (recursive) procedure \verb2case_12 which implements the first of the above cases. 
It takes as argument the list *K of all previously found polyhedra, the currently analyzed polyhedron (or rather its initial segment, of length l),
the remaining possible edges to be added, encoded in Res (r is the length of Res). 
 
\smallskip

 void \verb2case_12(List *K,  polyhed *t, char l, char Res [][2], char r)\{

\hspace{.1cm} char X[r][2];
 char x,u,v;

\hspace{.1cm} for(int a=0; a$<$r; a++)\{
   u=Res[a][0];

\hspace{.5cm}   if (Res[a][1]==1)\{
     for(int i=0; i$<$a; i++)\{X[i][0]=Res[i][0];X[i][1]=Res[i][1];\}
     
  \hspace{.9cm}    for(int i=a+1; i$<$r; i++)\{X[i-1][0]=Res[i][0];X[i-1][1]=Res[i][1];\}
    x=r-1;\}
    
\hspace{.5cm}   else\{
    for(int i=0; i$<$r; i++)\{X[i][0]=Res[i][0];X[i][1]=Res[i][1];\}
     
   \hspace{.9cm}   x=r;X[a][1]--\,--;\}
   
\hspace{.5cm}   for(int b=0; b$<$x; b++)\{
           char Y[x][2];
        char y; 
      v=X[b][0];
      
  \hspace{.9cm}         if (X[b][1]==1)\{
         for(int i=0;i$<$b;i++)\{Y[i][0]=X[i][0];Y[i][1]=X[i][1];\}
         
     \hspace{1.4cm}        for(int i=b+1;i$<$x;i++)\{Y[i-1][0]=X[i][0];Y[i-1][1]=X[i][1];\}
                 y=x-1;\}
                 
  \hspace{.9cm}            else\{
        for(int i=0;i$<$x;i++)\{Y[i][0]=X[i][0];Y[i][1]=X[i][1];\}
           y=x;Y[b][1]--\,--;\}
           
 \hspace{.9cm}         t$\to$length=l;
 
 \hspace{.9cm}         if (t$\to$\verb1test_link1(l\%3+1,u,v)==false)\{
            t$\to$\verb1add1(l\%3+1,u,v);
            
\hspace{1.4cm}             if (l$<$7)\{\verb2case_12(K,t,l+1,Y,y);\} else \{K$\to$\verb1add1(t);\}\}\}\}\}

\smallskip

\noindent The function \verb2case_12 is initialized in the \verb1main1 procedure as:
 
\smallskip

 int \verb1main1 () \{

\hspace{.1cm}        char a [2][3]=\{\{1,2,3\},\{4,4,5\}\};
       polyhed  t(a,2); List L;

        \hspace{.1cm} char S[5][2]=\{\{4,1\},\{5,2\},\{6,3\},\{7,3\},\{8,3\}\};

\hspace{.1cm}        \verb2case_12(\&L,\&t,2,S,5);
  L.display();        return 0;
      \}

\smallskip

\noindent The functions implementing the six remaining cases are  of a similar nature.

At this stage one last case remains, namely the situation where a face $(1,2,3)$ corresponding to pairwise non adjacent edges in the link $L$ doesn't exist. In such a polyhedron, \emph{every} face contains a pair of adjacent edges in $L$. Let show first that such a polyhedron must have 5 adjacent identifications.  If it has 6, then it admits a presentation of the form
$[(1,1,7),(2,2,7), (3,3,7),(4,4,8),(5,5,8),(6,6,8),(a,b,c),(d,e,f)]$ where $\{a,b,c,e,d,f\}=\{1,2,3,4,5,6\}$ (where  both $(7^\flat,7^\sharp)$ and $(8^\flat,8^\sharp)$ have  rank 2 and are at distance 3 in $L$). Furthermore, one of the two missing faces contains two numbers in either $\{1,2,3\}$ (corresponding to 7) or $\{4,5,6\}$ (corresponding to 8). But this creates a cycle of length $<6$ in $L$ and gives a contradiction. On the other hand, if the polyhedron $V$ has at most 4 adjacent identifications, then at least one of its faces correspond to pairwise nonadjacent edges in $L$. Indeed, it is clear that such a complex can't have $<4$ adjacent identifications, and if it has exactly 4, then it would admit a presentation of the form  $[(1,1,\cdot),(2,2,\cdot), (3,3,\cdot),(4,4,\cdot),(1,\cdot,\cdot),(2,\cdot,\cdot),(3,\cdot,\cdot),(4,\cdot,\cdot)]$. But then, again, any completion of the face $(1,\cdot,\cdot)$ (for example) will contradict the girth condition. Therefore, the only configuration that remains to be considered takes the form $[(1,1,6),(2,2,6), (3,3,\cdot),(4,4,\cdot),(5,5,\cdot),(7,\cdot,\cdot),(\cdot,\cdot,\cdot),(\cdot,\cdot,\cdot)]$. A computer function was implemented in this case too, and it turns out that there is only one such polyhedron, $V_5$, with 5 adjacent identifications. 

\newcommand{\rng}{ring}
\newcommand{\rngs}{rings}

\section{Proof of Theorem \ref{th7}} \label{pfth7}

\begin{proof}
Let $V$ be a compact complex of rank \sq\ (Definition \ref{comp74}). Let $\G$ be its fundamental group, $X$ its universal cover and $\pi: X\to V$ the covering map.
Fix a vertex $A_0$ of $X$ and identify $\G$ with its orbit $\G A_0$ in $X$.
By  Dykema and de la Harpe's Theorem 1.4 in \cite{DH}, it is enough to show that $\G$ has the free semigroup property and the $\ell^2$-spectral radius property. The latter follows from property RD (see Proposition \ref{spectr}) which itself follows from Theorem \ref{th1} so we aim  to prove the free semi-group property in the present section. Recall \cite{DH} that a group $\G$ is said to have the \emph{free semi-group property} if for every finite subset $F$ of  $\G$ there is $u\in \G$ such that the set $uF=\{ua, ~a\in F\}$ is semi-free. A finite subset $F$ of $\G$ is said to be \emph{semi-free} if for every $n,m\in \NI$, every $x_1,\ldots x_n,y_1\ldots y_m\in F$, the equality $x_1\ldots x_n=y_1\ldots y_m$ implies that $m=n$ and $x_i=y_i$ for every $i=1\ldots n$.  

In order to prove the free semi-group property for $\G$, which is achieved in Lemma \ref{sf}, we need to introduce one more concept for rank \sq\ polyhedra. Negative curvature is used extensively  below via the Gauss-Bonnet formula.

\begin{definition}\label{def-ana}
We call \emph{analytic geodesic} of $X$  a singular CAT(0) geodesic (i.e.\ included in the 1-skeleton, see Section \ref{rd-triangles}) whose angle at each  vertex equals $4\pi/3$. 
\end{definition}

As $V$ is compact  the family of analytic geodesics of $X$ projects under $\pi$ to a finite set of closed geodesics of $V$.  We call these closed geodesics the \emph{\rngs} of $V$. 
Note that   every analytic geodesic is periodic of period the length of its corresponding \rng\ in $V$ (because $\G$ acts freely on $X$ with compact quotient $V$).

\begin{example} One can show that the complex 
$\bar V$ of Remark \ref{barv} has a single ring $r=841^{-}6537^{-}2$.  Thus all analytic geodesics in the universal cover of $\bar V$ have period 8. Moreover the image of $r$ in $H_1(\bar V,\ZI)=\ZI$ is  equal to 8.
\end{example}

\begin{lemma}[Analyticity]\label{anal}
Let $\gamma_1$ and $\gamma_2$ be two analytic geodesics of $X$. Then exactly one of the following cases occurs:
\begin{enumerate}
\item $\gamma_1$ and $\gamma_2$ are disjoint;
\item the intersection of $\gamma_1$ and $\gamma_2$ is reduced to a point;
\item $\gamma_1=\gamma_2$.
\end{enumerate}
\end{lemma}

\begin{proof}
Assume that $I=\gamma_1\cap\gamma_2$ contains at least two distinct points and let us prove that $\gamma_1=\gamma_2$. One has  $I=[A,B]\cap X$ for some points $A,B\in \bar X$, where $\bar X$ is the disjoint union of $X$ and  its boundary. Assume that $A\in X$. Then $A$ is a vertex of $X$ and as $B\neq A$ by assumption one of the edges of $X$ containing $A$, say $[A,A']$,  is included in $I$. However  the link $L_A$ (being isometric to $L_\sqm$) contains a unique point $A''$ which is at (angular) distance $4\pi/3$ from $A'$. By definition of analyticity $A''$ belongs to both $\gamma_1$ and  $\gamma_2$, which contradicts the definition of $I$. Thus $A\in\bar X\backslash X$ and similarly $B\in\bar X\backslash X$ so $\gamma_1=\gamma_2$.   
\end{proof}

\begin{lemma}\label{s}
Let $\gamma$ be an analytic geodesic of period  $t\in \G$ in $X$.  There exists an $s$ in $\G$ such that for $A\in X$ the geodesic segment $\eta$ from $A$ to $tstA$ contains  $tA$, $stA$, and is not included in an analytic geodesic. 
\end{lemma} 

\begin{proof}
Let $A_0$ be a point of $\gamma$ and $g$ be an analytic geodesic of $X$ such that $\gamma\cap g=\{A_0\}$. Let $u$ be the period of $g$ and write $B_0=u^2A_0$. Then the unique analytic geodesic  $\gamma'$ which contains $[B_0,tB_0]$ does not intersect $\gamma$. Indeed assume it does and denote $C_0=\gamma\cap \gamma'$. Then by analyticity the angles of the geodesic triangle $\Delta =(A_0B_0C_0)$ at $A_0,B_0,C_0$ are at least $\pi/3$. In particular $\Delta$ bounds a topological disk $D$. With the notations of Lemma \ref{equi},  the Gauss-Bonnet formula for $D$ gives
\[
\int_{\Delta} \kappa+ (\pi-\theta_{A_0}) 
+(\pi-\theta_{B_0})+(\pi-\theta_{C_0})+~~\sum_{x\in D}  \delta _x~~=~~2\pi.
\] 
Recall that $\sum_{x\in D}  \delta _x $ is the internal curvature of $D$, so  $\sum_{x\in D}  \delta _x\leq 0$, and   $\int_{\Delta} \kappa$ is total geodesic curvature of $D$ on its boundary, that is $\int_{\Delta} \kappa=\int_{]A_0,B_0[} \kappa+\int_{]B_0,C_0[} \kappa+\int_{]A_0,C_0[} \kappa$ where
\[
\int_{]A_0,B_0[} \kappa= \sum_{s\in]A_0,B_0[} (\pi-\theta_s)
\]
with $\theta_s$ is the internal angle of $D$ at $s\in ]A_0,B_0[$ (and similarly for the two other sides). Thus 
\[
\int_{\Delta} \kappa\leq \int_{]A_0,B_0[} \kappa= - \pi/3
\]
by analyticity of $[A_0,B_0]$. However the inequality $(\pi-\theta_{A_0}) +(\pi-\theta_{B_0})+(\pi-\theta_{C_0})\leq 2\pi$ gives
\[
\int_{\Delta} \kappa\geq 0
\]
which is a contradiction. 

Let $A_1=t^{-3}A_0$ (resp. $B_1= t^3B_0$) and denote $A_1'$ (resp. $B_1'$) the point of $[A_1,A_0]\cap [A_1,B_1]$ (resp. $[B_0,B_1]\cap [A_1,B_1]$) such that $|A_1-A_1'|$ (resp. $|B_1-B_1'|$) is maximal.  We now show that $A_1\neq A_1'$ and $B_1\neq B_1'$.

So assume toward a contradiction that $A_1=A_1'$. Let $h=[A_1,A_0] \cup [A_0,B_0]\cup [B_0,B_1']$. By the above $h$ is a piecewise geodesic path in $X$ from $A_1$ to $B_1'$ without self-intersection. Let $B_1''$ be the closest point from  $A_1$ on $h$ which belongs to $[A_1,B_1']$ and let $h'$ be the part of $h$ going from $A_1$ to $B_1''$.
Consider the topological disk $D'$ whose boundary is the piecewise geodesic closed curve $\Delta'=[A_1,B_1'']\cup h'$.  Note that, the internal angles of $D'$ at $A_1$ being non zero, the point $B_1''$ does not belong to $\gamma$ and so $A_0\in h'$. The point $B_0$ might or might not belong to $h'$; in both cases the following follows from the Gauss-Bonnet formula:
\begin{align*}
\int_{\Delta'} \kappa&\geq 2\pi - (\pi-\theta_{A_1}) 
-(\pi-\theta_{A_0})-(\pi-\theta_{B_0})-(\pi-\theta_{B_1''})-~~\sum_{x\in D'}  \delta _x\\
&>-(\pi-\theta_{A_0})-(\pi-\theta_{B_0})\geq- 2\pi/3.
\end{align*}
Then analyticity of $[A_0,A_1]$ implies 
\[
\int_{\Delta} \kappa\leq \int_{]A_0,A_1[} \kappa=-\pi/3(|A_0-A_1|-1)\leq - 2\pi/3
\]
which is a contradiction.  

Thus $A_1\neq A_1'$ and similarly $B_1\neq B_1'$. 
But this shows that the element $s=t^3u^2t^3$ of $\G$ satisfies that, for any $A\in X$, the geodesic segment $\eta$ from $A$ to $tstA$ contains both $tA$, $stA$. By construction,  $\eta$ is not included in an analytic geodesic. 
\end{proof}

Let $F$ be a finite subset of $\G$ and let $\alpha= \max_{a\in F}|a|$. Fix a analytic geodesic $\gamma$ of period $t$ in $X$ and  let $s$ be as in Lemma \ref{s}. Let $\beta$ be the length of $s$ and denote 
 \[
u=t^{4\alpha+7}st^{9\alpha+\beta+19}.
\]

 \begin{lemma}\label{lem40}
 Let $a\in F$ and $A\in X$. Then $t^{\alpha+3}aA$ is on the geodesic segment  from $A$ to $uaA$ in $X$.
 \end{lemma}

\begin{proof}
Let $B=aA$, $C=uaA$ and  consider the points $A',B',C'$ of $[A,B]\cap [A,C]$, $[B,A]\cap [B,C]$ and $[C,A] \cap [C,B]$ respectively such that $|A-A'|$, $|B-B'|$, and $|C-C'|$ are maximal. By assumption we have $|A-B|\leq \alpha$ so $|B-B'|\leq \alpha$. Let $D$ be the unique disk whose boundary is the geodesic triangle  $\Delta=[A',B']\cup [B',C']\cup [C',A']$. If $\Delta$ is reduced to a point the lemma is clear so we assume this is not the case. If not we apply the Gauss-Bonnet formula to $D$:
\[
\int_{\Delta} \kappa+ (\pi-\theta_{A'}) 
+(\pi-\theta_{B'})+(\pi-\theta_{C'})+~~\sum_{x\in D}  \delta _x~~=~~2\pi.
\] 
which gives 
\[
\int_{\Delta} \kappa> 2\pi -3\pi - \sum_{x\in D}  \delta _x\geq -\pi
\]
as  $(\pi-\theta_{A'}) +(\pi-\theta_{B'})+(\pi-\theta_{C'})<3\pi$. On the other hand
\[
\int_{\Delta} \kappa\leq \int_{]B',C'[} \kappa= \sum_{s\in]B',C'[} (\pi-\theta_s)
\]
As $|B-B'|\leq \alpha$  the point $C''=t^{\alpha+4}B$ belongs to  $]B',C[$ and $[B',C'']$ is an analytic geodesic of $X$.  Write $]B',\tilde C[$ for the initial segment $]B',C'[\cap ]B',C''[$. We have   
\[
\int_{]B',C'[} \kappa\leq  \sum_{s\in]B',\tilde C[} (\pi-\theta_s)=-{\pi\over 3} \min\{|B'-C'|-1,|B'-C''|-1\}
\]
because $]B',C'[$ is geodesic and $]B',\tilde C[$ is analytic. 
Therefore
\[
-{\pi\over 3} \min\{|B'-C'|-1,|B'-C''|-1\}\geq \int_{]B',C'[} \kappa \geq \int_{\Delta} \kappa >-\pi
\]
and so $\min\{|B'-C'|-1,|B'-C''|-1\}< 3$. As $|B'-C''|\geq 4$ it follows that $\tilde C=C'$ and $|B'-C'|\leq 3$ so we have $t^{\alpha+3}B\in [C',C]$. This proves the lemma.
 \end{proof}

\begin{lemma}\label{lem41}
Let $(a_{1},\ldots a_{n})$ be a sequence of elements of $F$ of length $n\in \NI$. For $k=0\ldots n-1$ define recursively   points $x_k$ in $\G$ (viewed as a subset of $X$) by  
\[
x_{k+1}=u a_{{k+1}} x_k,
\] 
where $x_0=e$. Let $\eta_k$ be the geodesic segment from $x_k$  to $x_{k+1}$ in $X$. For any $k=0\ldots n-1$ there exist two points $x_k^+<x_{k+1}^-$ on $\eta_k$ such that $|x_k^+-x_k|\leq (2\alpha+6)|t|$, $|x_{k+1}^--x_{k+1}|\leq (2\alpha+6)|t|$, and such that  $[x_k^+,x_{k+1}^-]$ are pairwise disjoint consecutive geodesic segment on the geodesic segment from $x_0$ to $x_n$.  
\end{lemma}

\begin{proof}
For $i=1\ldots n-1$ let $\tilde x_i$ be the point of $\eta_{i-1}\cap\eta_i$ for which $|x_i-\tilde x_i|$ is maximal. Let us first show that $|x_i-\tilde x_i|\leq (2\alpha+3)|t|$.
By Lemma \ref{lem40}, the point $z=t^{\alpha+3}a_ix_i$ belongs to $\eta_{i}$ and $|x_i-z|\leq (2\alpha+3)|t|$. Assume toward a contradiction that  $z\in ]x_i,\tilde x_i[$. Then the intersection of the segments $[z,x_{i}]\subset \eta_{i-1}$ and $[z,x_{i+1}]\subset \eta_{i}$  contains at least an edge which is readily seen to be in the analytic parts of $\eta_{i-1}$ and $\eta_i$.  From the definition of $u$ and Lemma \ref{anal} we infer that $\eta_{i-1}\cap \eta_i\supset [x_i,z']$
where $z'=t^{-2\alpha-5}x_i$. Thus as $z=t^{\alpha+3}a_ix_i\in \eta_{i-1}$ we have
\[
z\in[t^ez',t^{e+1}z']
\]
for some index $2\leq e\leq 2\alpha+5$. However $z'\in \eta_i$ as well and as  $z=t^{\alpha+3}a_ix_i\in \eta_i$ we also have  
\[
z'\in [t^fz, t^{f+1}z]
\]
for some index $2\leq f\leq 2\alpha+5$. In particular $|tz-tz'|<|z-z'|$ which contradicts the fact that, being its period, $t$ acts isometrically on $\gamma$. 
Thus $\tilde x_i\in [x_i,z]$ and $|x_i-\tilde x_i|\leq (2\alpha+3)|t|$.
Set $\tilde x_0=x_0, \tilde x_n=x_n$ and let $\tilde \eta_i$ be the geodesic from $\tilde x_{i}$ to $\tilde x_{i+1}$, $i=1\ldots n-1$ (so  $\tilde \eta_i\subset \eta_i$).
  
  Let us prove that $\tilde \eta_i$ intersects $\tilde \eta_j$ if and only if they are consecutive (i.e.\ $i=j+1$ or $j=i+1$) for $i,j=1\ldots n-1$.  Suppose on the contrary that there is a $i\in [0,n-2]$ and a $j>i+1$ such that $\tilde \eta_i\cap \tilde \eta_j\neq \emptyset$. We can further assume that $j$ is the smallest index $j>i+1$ satisfying this condition. Let $z\in \tilde \eta_j$ be the closest point from $\tilde x_j$ which belongs to $\tilde \eta_i$. Denote  $\tilde \eta_i'=[z,\tilde x_{i+1}]$ and $\tilde \eta_j'=[\tilde x_j,z]$. Then there is a topological disk $D$ in $X$ whose boundary is the  piecewise geodesic simple curve
 \[
\Delta= \tilde \eta_i'\cup (\bigcup_{i<k<j}\tilde \eta_k)\cup \tilde \eta_j'.
 \]
The Gauss-Bonnet formula for $D$ reads 
\[
\int_{\Delta} \kappa+ (\pi-\theta_{z}) +\sum_{i< k \leq j}(\pi-\theta_{\tilde x_k})
+~~\sum_{x\in D}  \delta _x~~=~~2\pi.
\] 
As $(\pi-\theta_{z}) +\sum_{i< k \leq j}(\pi-\theta_{\tilde x_k})<(j-i+1)\pi$ and $\sum_{x\in D}  \delta _x\leq 0$ we deduce that
\[
\int_{\Delta} \kappa>2\pi -(j-i+1)\pi=(i+1-j)\pi.
\]
On the other hand
\[
\int_{\Delta} \kappa=\sum_{i<k<j} \sum_{s\in \tilde \eta_k} (\pi-\theta_s)+\sum_{s\in \tilde \eta_i'} (\pi-\theta_s) +\sum_{s\in \tilde \eta_j'}+ (\pi-\theta_s) \leq \sum_{i<k<j} \sum_{s\in \tilde \eta_k} (\pi-\theta_s).
\]
So consider the geodesic $\tilde \eta_k\subset \eta_k$ for $k\in ]i,j[$. As $|x_k-\tilde x_k|\leq (2\alpha+3)|t|$ and $|x_{k+1}-\tilde x_{k+1}|\leq (2\alpha+3)|t|$ we can find a segment $S_k$ of length at least $3\alpha$ in $\tilde \eta_j$ which is analytic. In particular
\[
\sum_{s\in \tilde \eta_k} (\pi-\theta_s)\leq \sum_{s\in S_k} (\pi-\theta_s)\leq 3\alpha \cdot (-\pi/3)=-\alpha \pi  
\]
Thus we get
\[
\int_{\Delta} \kappa\leq -\alpha(j-i-1) \pi
\]
and so
\[
(i+1-j)\pi<-\alpha(j-i-1) \pi
\]
It follows that $\alpha<1$, which is a contradiction. This shows  that $\tilde \eta_i$ intersects $\tilde \eta_j$ if and only if $i=j+1$ or $j=i+1$ for $i,j=1\ldots n-1$.

We now prove the lemma. Let $g=\cup_{k=0\ldots n-1} \tilde \eta_k$. By the above $g$ is a piecewise geodesic curve in $X$ from $x_0$ to $x_n$ without self-intersection. 
We proceed by recurrence.  

Set $x_0^+=x_0$ and denote by $x_1^-$ the unique point of $\tilde \eta_0\cap [x_0,x_n]$ such that $]x_1^-,x_n[\cap \tilde \eta_0$ is empty. If $x_1^-=\tilde x_1$ and a neighbourhood  of $\tilde x_1$ in $\tilde \eta_1$ is included in $[x_0,x_n]$ then we let $x_1^+=\tilde x_1$ and the conditions of the lemma are satisfied at $x_1$ (as $|x_1-\tilde x_1|\leq (2\alpha+3)|t|$). Otherwise  let $x_1^+>x_1^-$ be the first point of $g$ distinct from $x_1^-$ which belongs to $[x_0,x_n]$. Note that $x_1^+\not \in \tilde \eta_0$ as this would imply $x_1^-=x_1^+$. As $g$ is a simple curve there is a non empty disk $D$ whose boundary 
\[
\Delta=g_{[x_1^-,x_1^+]}\cup [x_1^-,x_1^+]
\]
consists of the part $g_{[x_1^-,x_1^+]}$ of $g$ which is in between $x_1^-$ and $x_1^+$ and the geodesic segment $[x_1^-,x_1^+]\subset [x_0,x_1]$. Let $j\geq 1$ such that $x_1^+\in \tilde \eta_j$. The Gauss-Bonnet formula for $D$ shows that
\[
\int_{\Delta} \kappa+ (\pi-\theta_{x_1^-}) +(\pi-\theta_{x_1^+})+ \sum_{1\leq  k \leq j}(\pi-\theta_{\tilde x_k})
\geq ~~2\pi.
\] 
(where one could  a priori  have $x_1^-=x_1$ or $x_1^+=\tilde x_i$ for some $i>1$). Reiterating our argument above we get that $j=1$ and thus, $x_1^+\in \tilde \eta_1$. Hence the preceding formula implies
\[
\int_{\Delta} \kappa> 2\pi - 3\pi=-\pi.
\] 
If $x_1^+\in [x_1,t^{\alpha+3}ax_1]$ the condition $|x_1^+-x_1|\leq (2\alpha +6)|t|$ of the lemma is satisfied. It follows that  we  can assume $x_1^+\in [t^{\alpha+3}ax_1,x_2]$. Then
\begin{align*}
\int_{\Delta} \kappa&=\sum_{s\in ]x_1^-,\tilde x_1[}(\pi-\theta_s)+\sum_{s\in ]\tilde x_1,x_1^+[}(\pi-\theta_s)+\sum_{s\in ]x_1^-,x_1^+[}(\pi-\theta_s)\\
&\leq\sum_{s\in ]x_1^-,\tilde x_1[}(\pi-\theta_s)+\sum_{s\in ]\tilde  x_1,x_1^+[}(\pi-\theta_s)\\
&\leq\sum_{s\in ]x_1^-,\tilde x_1[}(\pi-\theta_s)+\sum_{s\in ]t^3a x_1,x_1^+[}(\pi-\theta_s)\\
&\leq \min\{4\alpha|t|,|x_1^--\tilde x_1|\}\cdot (-\pi/3) + \min\{4\alpha|t|,|x_1^+-t^{\alpha+3}a x_1|\}\cdot (-\pi/3)
 \end{align*}
where the last inequality comes from the fact that  the geodesics  $[x_1^-,\tilde x_1]$ and $[t^{\alpha+3}ax_1,x_2]$ are analytic on a segment of length at least $4\alpha|t|$, starting from $\tilde x_1$ and $t^{\alpha+3}ax_1$ respectively. Thus
\[
\min\{4\alpha|t|,|x_1^--\tilde x_1|\} +\min\{4\alpha|t|,|x_1^+-t^{\alpha+3}a x_1|\}\leq 3  
\]
As $4\alpha|t|>3$ this implies $|x_1^--\tilde x_1|\leq 3$ and $|x_1^+-t^{\alpha+3}a x_1|\leq 3$. Thus $|x_1^--x_1|\leq (2\alpha +6)|t|$ and $|x_1^+-x_1|\leq (2\alpha +6)|t|$ as asserted. 

The construction of $x_k^-$ and $x_k^+$ (the latter provided $k<n$, where $x_n^-=x_n$) and the proof that $|x_k^--x_k|\leq (2\alpha +6)|t|$ and $|x_k^+-x_k|\leq (2\alpha +6)|t|$  can then be done exactly as for $x_1^-$ and $x_1^+$ above so  we will omit the details. 

This concludes the proof of Lemma \ref{lem41}.
\end{proof}

Recall that $u=t^{4\alpha+7}st^{9\alpha+\beta+19}$ has been defined in the paragraph preceding Lemma \ref{lem40}.

\begin{lemma}\label{sf}
The finite set   $u F$ is semi-free in $\G$.
\end{lemma}

\begin{proof}
Let $(a_{1},\ldots a_{n})$ and $(b_{1},\ldots b_{m})$ be a sequence of elements of $F$ of length $n,m\in \NI$ respectively. Assume that 
\[
ua_{n}\ldots ua_{1}=ub_{m}\ldots ub_{1}
\]
in $\G$. Let $(x_0,\ldots,x_{n})$ and $(y_0,\ldots,y_{m})$ be the sequences of point of $X$ associated to $(a_{1},\ldots a_{n})$ and $(b_{1},\ldots b_{m})$ respectively as in Lemma \ref{lem41} (so $x_0=y_0=e$ and $x_n=y_m=ua_{n}\ldots ua_{1}$ assuming $\G\subset X$). Associated to $(x_0,\ldots,x_{n})$ (resp. $(y_0,\ldots,y_{m})$) we can find  points  $x_0^+,x_1^{\pm}\ldots, x_{n-1}^\pm,x_n^-$  (resp. $y_0^+,y_1^{\pm}\ldots, y_{n-1}^\pm,y_n^-$) on the geodesic $[x_0,x_n]$  which satisfy the conclusion of Lemma \ref{lem41}. 
Fix $k\in \{0\ldots n-1\}$. Since we have $|x_{k+1}^--x_{k+1}|\leq (2 \alpha+6)|t|$ and $|x_k^+-x_k|\leq (2\alpha+6)|t|$  the point $w_k=t^{3\alpha+\beta+6}a_kx_k$ and the point $st^{6\alpha+13}w_k$ both belong to  $[x_k^+,x_{k+1}^-]$. Note that $w_0<w_1<\ldots <x_{n-1}$ on $[x_0,x_n]$. Similarly, define points $w_k'$ on $[x_0,x_n]$, $k=0\ldots m-1$, relative to the segments $[y_k^+,y_{k+1}^-]$ and satisfying analogous properties.

We claim that if $w$ is a point on $[x_0,x_n]$ such that $st^{6\alpha+13} w\in [x_0,x_{n}]$ then $w=w_k$ for some $k=0\ldots n-1$.

Let us first show that this implies the lemma. Indeed by the claim  there is an increasing injection $j : [0,\ldots ,m-1]\to [0,\ldots ,n-1]$ such that $w'_{k}=w_{j(k)}$. Thus  by symmetry we  obtain that  $n=m$ and $w_k'=w_k$ for $k=0\ldots n-1$. On the other hand we have  
$x_{k+1}=t^{4\alpha+7}st^{6\alpha+13}w_k$ and so 
\[
y_{k+1}=t^{4\alpha+7}st^{6\alpha+13} w_k'=x_{k+1}.
\]
 Thus 
\[
ua_{k+1}=x_{k+1}x_k^{-1}=y_{k+1}y_{k}^{-1}=ub_{k+1}
\]
for $k=0\ldots n-1$. This shows that $uF$ is semi-free.

Let us now prove the claim. Let $w$ is a point on $[x_0,x_n]$ such that $w'=st^{6\alpha+13} w\in [x_0,x_{n}]$. We first show that $w\in [x_k^+,st^{6\alpha+13} w_k]$  for some index $k\in [0,n-1]$. If not then $w\in [st^{6\alpha+13} w_{k-1},x_{k}^+]$ for some $k\in [1,n-1]$. As $|x_{k}-x_{k}^+|\leq  (2\alpha+6)|t|$ and $|st^{6\alpha+13} w_{k-1}-x_{k}|\leq (4\alpha+7)|t|$ we have 
\[
|x_{k}^+-st^{6\alpha+13} w_{k-1}|\leq  6\alpha|t|+13|t|.
\]
Lemma \ref{lem40} and Lemma \ref{lem41} then show that
that the point $t^{6\alpha+13}w$  of $[x_0,x_n]$ is of the form 
$t^\ell a_kx_k$ for some index $\ell>\alpha + 3$, and as $|x_{k}-x_{k}^+| \leq (2\alpha+6)|t|$ we also have 
\[
\ell\leq  (6|t|\alpha +13|t|)+(2\alpha|t| +6|t|)=8\alpha|t|+19|t|.
\]
As  $8\alpha+19<9\alpha+19$ our choice of $u$ (more precisely the exponent $9\alpha+\beta+19$) shows that the points $t^{6\alpha+13}w$ and $st^{6\alpha+13}w$ are extremities of an analytic subsegment  of $[x_0,x_n]$. But this  contradicts  Lemma \ref{s} (i.e.\ the non analyticity of $s$).

Therefore $w\in [x_k^+,st^{6\alpha+13} w_k]$  for some index $k\in [0,n-1]$. Now the segment $[w,t^{6\alpha+13}w]$ is an analytic subsegment of $[x_0,x_n]$ and thus is disjoint from $[t^{6\alpha+13}w_k,st^{6\alpha+13}w_k]$ by construction of $s$. Thus $[w,t^{6\alpha+13}w]\subset [x_k^+,t^{6\alpha+13}w_k]$. However by definition of $w_k$ the segment $[x_k^+,t^{6\alpha+13}w_k]$ is analytic as well, hence disjoint from $[t^{6\alpha+13}w,st^{6\alpha+13}w]$. This implies that $st^{6\alpha+13}w_k=st^{6\alpha+13}w$ and thus $w_k=w$, proving the claim.

This concludes the proof of Lemma \ref{sf}.
\end{proof}

Theorem \ref{th7}  follows from Lemma \ref{sf} and the paragraph preceding Definition \ref{def-ana}.
\end{proof}

As announced in Section \ref{asymsq}, we conclude this section with the proof of Proposition \ref{altern}.

\begin{proof}[Proof of Proposition \ref{altern}]
Let $V$ be a complex of rank \sq,  $X=\tilde V$, and $\G=\pi_1(V)$. Let $\Lambda$ be a subgroup of $\G$ isomorphic to $\ZI^2$. 
Fix  an analytic geodesic $\gamma$ in $X$ of period $t\in \G$ an set 
\[
u=t^{|s|+1}st^{|s|+1}
\]
where $s$ is given by Lemma \ref{s} and $|s|$ is its length.

Let us show that the pairwise intersection of the subgroups  $u^n \Lambda u^{-n}$, $n\in \ZI$,  is reduced to the identity.
Indeed assume that there is $n,m\in \ZI$, and $\lambda_1,\lambda_2\in \Lambda$ such that 
\[
u^n \lambda_1 u^{-n}=u^m \lambda_2 u^{-m}
\]
and let us show that $n=m$ or $\lambda_1=\lambda_2=e$. We have
\[
u^{n-m} \lambda_1 u^{m-n}= \lambda_2 \in \Lambda
\]
so consider the geodesic parallelogram $P$ in $X$ with vertices $A_0$, $A_1=\lambda_1A_0$, $A_2=u^{n-m}\lambda_1A_0$ and $A_3=u^{n-m}A_0$. We assume that $n\neq m$ and prove that $\lambda_1=\lambda_2=e$.

As the segment $[A_0, A_1]$ is flat by assumption, the point $A_0'\in [A_0, A_1]\cap [A_0, A_3]$ such that $|A_0-A_0'|$ is maximal is at distance at most 1 from $A_0$ by definition of $u$. Similarly  define (using flatness of $[A_0, A_1]$ and $[A_2, A_3]$) points  $A_1'$, $A_2'$, $A_3'$ in  $[A_1, A_0]\cap [A_1, A_2]$, $[A_2, A_1]\cap [A_2, A_3]$, and $[A_3, A_2]\cap [A_3, A_0]$ respectively, sharing analogous properties. Using the Gauss-Bonnet formula one can prove that $A_0'=A_1'$ and $A_2'=A_3'$. We shall omit the details as they are similar to that appearing in the middle of the proof of Lemma \ref{lem41}. Then  arguing as in the proof of Lemma \ref{sf} (using the element $s$ in the definition of $u$) we get that $A_0=A_1$ and $A_2=A_3$. (We leave the details to the reader as well.) The action of $\G$ on $X$ being free, we deduce that $\lambda_1=\lambda_2=e$. This proves the proposition. 
\end{proof}

\section{Mesoscopic  rank}\label{meso}

We first prove the statements in the first part of Subsection \ref{meso-intr}.

\begin{proposition} Let $X$ be a locally compact triangle polyhedron.
\begin{enumerate}
\item If $X$ is hyperbolic then the support of $\varphi_A$ is a relatively compact subset of $\RI$ for any vertex $A\in X$.
\item If $X$ has local rank $\leq 3/2$, then  for  every vertex $A\in X$ the support of $\varphi_A$ is included in a compact subset of $\RI$, excepts perhaps for a union of semi-open intervals   $I_1, I_2, I_3,\ldots$ of $\RI$ of the form 
\[
I_p=\left ](2p+1){\sqrt 3\over 2},\sqrt{3p(p+1)+1}\right],~~p\in \NI,
\] 
 where, observe, $|I_p|\to_p 0$. If moreover the order of $X$ (maximal number of triangle adjacent to an edge)  is bounded, then $\varphi_A$ is bounded.
\item If $X$ has local rank 2 (i.e.\ if $X$ is a triangle building) then $\varphi_A$ vanishes identically for any $A\in X$.
\end{enumerate}  
\end{proposition}

\begin{proof}
(1)  Assume that there is some vertex $A\in X$ such the support of $\varphi_A$ is unbounded. Let  $S\subset \RI$ be the support $\varphi_A$ and  $(r_n)_n$  be a sequence of points in $S$ converging to $\infty$. For each $n$ let $D_{r_n}$ be a  flat disk in $X$ with center $A$ and radius $r_n$.
Note that by the local compactness assumption, for every $r\in \RI$ the number of flat disks of radius $r$ in $X$ with center $A$ is finite. In particular there exists an infinite subsequence $S_1$ of $(r_n)_n$ such that every disk $D_r$ for $r\in S_1$ coincide on the ball of radius 1 and center $A$. Iterating, there exists for every $k\in \NI^*$ a infinite subset $S_{k+1}$ of $S_k$ such that every disk $D_r$ for $r\in S_{k+1}$ coincide on the ball of radius $k+1$ with some fixed flat disk $F_k$ of center $A$ and radius $k$. Now the increasing union $F=\cup_k F_k$ is a flat in $X$.  Thus $X$ contains a flat and hence is not hyperbolic by the no flat criterion.

(2) Assume now that the local rank of  $X$ is $\leq \tdm$.  Let $A$ be a vertex of $X$. 
We  claim that there is a finite number of flat hexagon in $X$ (simplicially isometric to hexagons in $\RI^2$ endowed with the tessellation by equilateral triangles) centered at a given vertex $A$, which are not included in a flat of $X$. 

Indeed let us prove  that any two hexagons of $X$ of same simplicial radius which coincide on the ball of simplicial radius $2$ are equal.  We argue by recurrence.  So let $H_1$ and $H_2$ be two flat hexagons of $X$ and let $B_n$ be a simplicial hexagon of radius $n$ in $X$ on which $H_1$ and $H_2$  coincide. We assume that  $n\geq 2$ and prove that $H_1$ and $H_2$ coincide on an hexagon of radius $n+1$ (provided $n+1$ is no greater than the common radius of $H_1$ and $H_2$). 
As $n\geq 2$ there is a vertex $A'$ in the boundary of $B_n$ whose internal angle in $B_n$ is $\pi$. In the link $L_{A'}$ of $X$ at $A'$ the two hexagons $H_1$ and $H_2$ generate two path of length $\pi$ creating with the path corresponding to $B_n$ two cycles of length $2\pi$ in $L_{A'}$.   By our local rank assumption (see the definition in the introduction), this cycles must coincide. But this is easily seen to imply that $H_1$ coincide with $H_2$ up to radius $n+1$. Thus $H_1=H_2$. 

This shows that the number of hexagons of $X$ of radius $n$ is bounded by the number of hexagons of $X$ of radius 2. Consider
 the family $\cF_A$ of maximally flat hexagons (i.e.\ not included in larger flat hexagons) centered at $A$ of radius at most 2. The preceding paragraph shows that the map on $\cF_A$ which associated to a maximally flat hexagon its simplicial sphere of radius 2 is injective.  This proves the claim.
 
Let    $r_0=r_0'+{\sqrt 3\over 2}$, where $r_0'$ is the maximal radius of the elements in $\cF_A$. 
 
Let $D$ be a flat disk of radius  $r>r_0$ and center $A$ in $X$ which is not included in a flat, and let $H'\supset D$ be the union of all triangle of $X$ whose interiors have non empty intersection with   $D$, so $H'$ is flat as $D$ is.  Let $H$ be a maximal hexagon of center $A$ in $H'$. Then $H$ is included in a flat $\Pi$ by definition of $r_0$.    As $D$ is not included in flat there is a triangle $t$ of $H'\backslash H$ which is not included in $\Pi$.  Let $U$ and $V$ be the two vertices of $t$ which belong to $H$. Our assumption that both links $L_U$ and $L_V$ have rank $\leq \tdm$ shows that the disk $D$ must actually be included in $H\cup t$. On the other hand the interior of $t$ has non empty intersection with $D$ by definition. A 2-dimensional computation provides the strong constraint on the radius of $D$ stated in the lemma, i.e.\  $r\in I_p$ for some $p\in \NI$ where
\[
I_p=\left ](2p+1){\sqrt 3\over 2},\sqrt{3p(p+1)+1}\right].
\] 
But this also shows that the subdisk of $D$  of radius $(2p+1){\sqrt 3\over 2}$  is contained in $\Pi$. By \cite{rigidite} there is only a finite number of distinct flat in $X$ containing $A$. As the number of triangle on every edge is uniformly bounded, we conclude that $\varphi_A$ is indeed bounded on $[r_0,\infty[$. More precisely  one has 
\[
\varphi_A \leq Nq^6
\]
where $N$ is  the number of flats and $q+1$ is the maximal valency of edges (i.e.\ $q$ is the order of $X$).  

The above proof can be easily adapted to the case of $\varphi_A$ for any point $A$ of $X$ (i.e.\ not necessarily a vertex). Then the intervals $I_p$ vary accordingly depending   
on the position of $A$ on its face.

(3) Assume finally first that $X$ has local rank 2.  We have to show that any flat disk in $X$ is included in a flat of $X$. Let $H_0\supset D$ be the union of all triangle of $X$ whose interiors have non empty intersection with   $D$. Then $H_0$ is flat as $D$ is, and the internal angle at every point in the boundary of $H_0$ is at most $4\pi/3$. By  applying local rank 2 a finite number of times (at most the number of vertices in the boundary of $H_0$), we deduce that there is a flat simplicial set $H_1$ whose interior contains  $H_0$ and which has interior angle  at most $4\pi/3$ at every point of its boundary. Iterating this we get a sequence of flat simplicial sets $H_0\subset H_1\subset H_2\ldots $ converging to a flat containing $D$.    
\end{proof}

\begin{remarks}
\begin{itemize} \item[(a)] The estimates in the proof of Assertion (2) can be made precised. Explicit computations in the case of the polyhedron of \cite{cras} (which is of rank \td) are summarized on Figure \ref{profile}. 
\item[(b)] The jumps in $\varphi_A$ may be considered as side effects as they are inherent of the fact we chose Euclidean disks in the definition of mesoscopic rank---what we did so as to have a definition applying to any CAT(0) polyhedron of dimension 2 (and higher).  In the triangle case they can be removed be choosing hexagons instead, as we saw along the proof. Observe however the mesoscopic rank behavior of Theorem \ref{th11} are \emph{not} side effects and cannot be removed by choosing hexagons.  
\item[(c)] The above proofs can be generalized to any CAT(0) space of dimension 2 under suitable uniform boundedness geometry assumptions.       
\item[(d)] The converse of Assertion (3) is true. In fact it is enough to assume  that $\varphi_A({1 \over 2})=0$ for any middle point $A$ of any edge of $X$, as is easily seen.  
\end{itemize}
\end{remarks}

In the remaining part of this Section we prove Theorem \ref{th11}. We start with the rank \sq\ case for  we  already are familiar with it from Section \ref{class-rg74}. 

\subsection{Proof of  Theorem \ref{th11}, Item  (b)} Recall that the complex $V_0^1$ admits the following presentation:
\[
V^1_0=[[1, 2, 3], [1, 4, 5], [1, 6, 4], [2, 6, 8], [2, 8, 5], [3, 6, 7], [3, 7, 5], [4, 8, 7]].
\]
We write $a_1,\ldots, a_8$ for the corresponding  generators of the fundamental group $\G$ of $V_0^1$. The 1-skeleton of  the universal cover $X$ of $V_0^1$ coincide with the Cayley graph of $\G$ with respect to  $\{a_1,\ldots a_8\}$ and each (oriented) edge in $X$ is labelled by the index in $\{1,\ldots,8\}$ of its generator.
A singular geodesic of $X$ is said to be of the form $i^\infty$, $i\in \{1,\ldots, 8\}$, if all its edges are labelled by $i$.

Recall (cf.\ \cite[p. 182]{BH}) that a simplicial subset $S$ of $X$ is called a (flat) \emph{strip}  if it is isometric to a product $I\times \RI\subset \RI^2$ where $I$ is a compact interval of $\RI$. The boundary of $S$ is a union of two (parallel) geodesics, say $g$ and $h$, and is denoted $(g,h)$. The \emph{height} of $S$ is the simplicial distance between $g$ and $h$. We say that a strip is \emph{periodic} if there is a non trivial $\gamma\in \G$ such that $\gamma(S)=S$.  Then there is $\gamma\neq e$ of smallest length satisfying this condition, and we call \emph{period} of $S$ the simplicial length of $\gamma$.

\begin{lemma}\label{strips}
There are 3 distinct strips of height 1 and period 1 in $X$ whose boundaries are of the form $(5^\infty,6^\infty)$.  
\end{lemma}

\begin{proof}
A immediate computation shows that the identities 
\[
a_1a_5=a_6a_1,~~ a_2a_5=a_6a_2,~~ a_3a_5=a_6a_3,
\]
and 
\[
a_5a_4=a_4a_6,~~ a_5a_8=a_8a_6, ~~a_5a_7=a_7a_6, 
\]
holds in $\G$. This gives exactly 3 distinct strips of height 1 and period 1 whose boundaries are of the form $(5^\infty,6^\infty)$ in $X$. 
\end{proof}

In particular there are uncountably many flats in $X$.
Extending the notation $i^\infty$ we encode simplicial paths in $X$ (more precisely their classes modulo  $\G$) by their corresponding sequence of labels 
and write $g^\infty$, given such a path $g$, for   the bi-infinite path obtained by juxtaposing of $g$ to itself on both sides infinitely many times.

\begin{lemma}\label{cl58}
Consider a simplicial path of $X$ of the form 
\[
g=271834.
\]
Then $g^\infty$ is geodesic and there are 2 distinct strips  of height 1 and period 6 in $X$ whose boundaries are of the form $(g^\infty, g^\infty)$.
\end{lemma}

\begin{proof}
Let $u=a_1a_8$, $v=a_2a_7$ and $w=a_3a_4$. Then a computation shows that $u$, $v$, and $w$ commute to $a_6$ in $\G$. The lemma readily follows.
\end{proof}

\begin{figure}[htbp]
\centerline{\includegraphics[width=14cm]{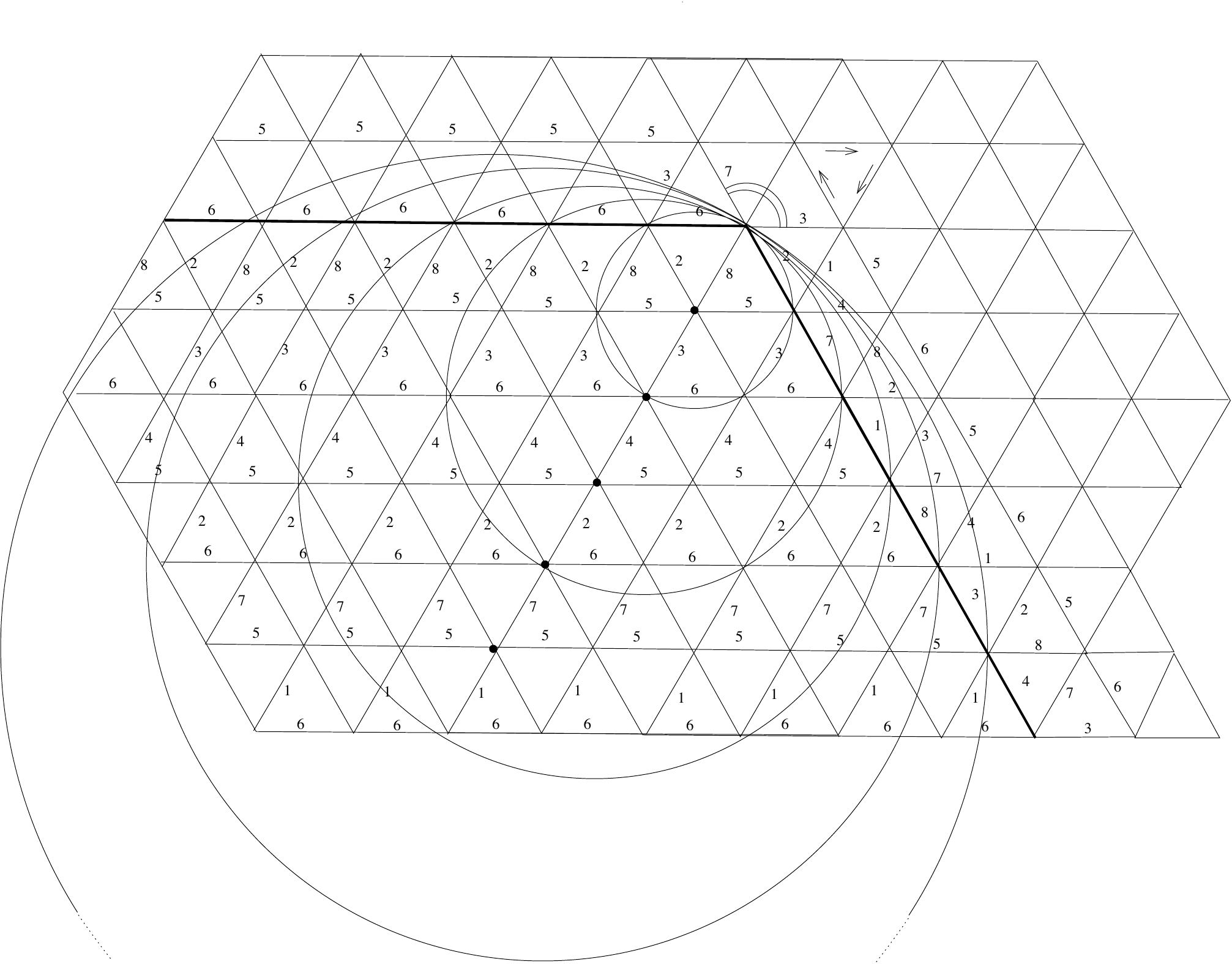}}
\caption{Exponential mesoscopic rank for $V_0^1$}\label{bibifurq}
\end{figure}

In particular the geodesic  $g^\infty$ is included in a flat $\Pi$ of $X$ on which  the group
\[
\Lambda= \langle wvu, a_6\rangle\simeq \ZI^2
\]
acts. The quotient space $\Pi/\Lambda$ has 12 triangles. Part of the  flat $\Pi$ and the geodesic segment $g$ are represented on Figure \ref{bibifurq}.

\begin{lemma}\label{cl59}
Keep the notations of Lemma \ref{cl58} and let  $h$ be a simplicial path of $X$ of the form
\[
h=65.
\]  
Then $h^\infty$ is a CAT(0) geodesics of $X$ and there are 3 distinct strips  of height 1 and period 6 in $X$ whose boundaries have the form  $(g^\infty,h^\infty)$.
\end{lemma}

\begin{proof}
Let $w_g=vuw$ and $w_h=a_5a_8$ be the words corresponding to $g$ and $h$. The Lemma follows from the identity  
\[
w_h^3w_g=w_gw_h^3
\] 
in $\G$, which is easily checked.
\end{proof}
 
In particular there are uncountably many semi-flats on the geodesic $g^\infty$. 
The following Lemma establishes Item (b) of Theorem \ref{th11} (as well as exponential mesoscopic rank).

\begin{lemma}\label{lem-meso74}
Let $A$ be a vertex of $X$ and $k$ be an integer greater than ${{\sqrt3}\over {2-\sqrt{3}}}=7.46...$. Then on the interval $]k-\sqrt 3,k]$  
of $\RI_+$ one has 
\[
\varphi_A\geq 2^{2\mu_k-4}
\]
where 
\[
\mu_k=\left \lceil k({2\over \sqrt3 }-1)\right\rceil.
\]
In particular $\varphi_A$ has continuous support starting from $7$, and exponential growth.
\end{lemma}

\begin{proof}
Let $\Pi$ be the flat of $X$ containing $A$ defined after Lemma \ref{cl58}. The geodesics of the form $6^\infty$ and $g^\infty$ intersect with internal angle $2\pi/3$ in $\Pi$, and there is  a unique vertex $B$ of $\Pi$ whose distance to $A$ is $k$, and such that the line segment  $(AB]$ of $\Pi$ is the bisector of this angle (see Figure \ref{bibifurq}). Denote by $d_1$ and $d_2$ the lines of $\Pi$ issued from $B$ corresponding to $6^\infty$ and $g^\infty$ respectively, and by $\Pi_0$ the sector of angle $2\pi/3$ at $B$ whose boundary is included in  $d_1\cup d_2$.

Let $\mu_k$ be the integer defined in the statement and let $\nu_k=2^{\mu_k-1}$.  By Lemma \ref{strips}  one can find $\nu_k$ distinct strips $\{S_1,\ldots S_{\nu_k}\}$ in $X$ of height $\mu_k$  whose boundaries all contain $d_1$, and which contain  the strip of height 1 on $d_1$ in $\Pi$ which is opposite to $\Pi_0$.   On the other hand by Lemma \ref{cl59} one can find $\nu_k$ distinct strips $\{T_1,\ldots T_{\nu_k}\}$ in $X$ of height $\mu_k$ whose boundary contains $d_2$, which contains $d_1$ as well as the strip of height 1 on $d_1$ defined in Lemma \ref{cl59}.  

Let $i\in \{1,\ldots,\nu_k\}^2$ and consider the subset $\Pi_i$ of $X$ defined by $\Pi_i=\Pi_0\cup S_{i_1}\cup T_{i_2}$. Then  our choice of $\mu_k$ shows, by an elementary exercise in  Euclidean geometry (in $\RI^2$),  that 
\begin{itemize}
\item the set $D_i$ of points of $\Pi_i$ at distance  $\leq k$ from $A$ in $\Pi_i$ is a flat disk in $X$ whose boundary contains $B$,
\item  the disks $D_i$ are pairwise distinct when $i$ varies in  $\{1,\ldots,\nu_k\}^2$. 
\end{itemize}

For $r\in [0,k]$ write $D_i^r$ for the concentric disk of radius $r$ in $D_i$. Assume that $\mu_k\geq  2$. Then it is not hard to show that for any fixed $r\in ]k-{\sqrt 3},k]$ the family of disks $\{D_1^r, \ldots D_{\nu_k}^r\}$ contains at least $2^{2\mu_k-4}$ distinct elements.
Observe that $\mu_k\geq 2$ if $k> {{\sqrt3}\over {2-\sqrt{3}}}$.

Let $j\in \{1,\ldots,\nu_k\}^2$ and $r\in ]k-{\sqrt 3},k]$. Let us now  show that the disk $D_j^r$ is not included in a flat. 

For a disk $D_i$, $i\in \{1,\ldots,\nu_k\}^2$, define $\tilde D_i$ to be the sector of center $A$ and angle $2\pi/3$ inside $D_i$  whose bisector is the  segment $[A,B]$.

\begin{claim}\label{clai60}
Every flat that contains $D_j^r$ must contain one of the sectors  $\tilde D_i$ for some $i\in \{1,\ldots,\nu_k\}^2$.  
\end{claim}

\begin{proof}[Proof of Claim \ref{clai60}]
 Note  that as $r>k-\sqrt 3$ the disk $D_j^r$  intersects some strips $T_{i_1}$ and $S_{i_2}$,  $i\in \{1,\ldots,\nu_k\}^2$, up to height 2 at least. Let $\Pi'$ be a flat that contains 
 $D_j^r$ and $\tilde \Pi'$ be  the sector of center $A$ and angle $2\pi/3$ inside $\Pi'$  whose bisector is the  segment $[A,B]$.  A local argument (along $d_1$ and $d_2$) shows that $\tilde \Pi'$   contains $A$ as well as the triangles in $T_{i_1}$ and $S_{i_2}$ adjacent to $d_1$ and $d_2$. Then as the transverse valency of the two sets $\cup  _{i\in \{1,\ldots,\nu_k\}} T_i$   and $\cup  _{i\in \{1,\ldots,\nu_k\}} S_i$ is maximal (i.e.\  equal to 3), the flatness of $\Pi'$ shows that the intersection of $\tilde \Pi'$ with the disk of center $A$ and radius $k$ in $\Pi'$ is of the form $\tilde D_i$ for some $i\in \{1,\ldots,\nu_k\}^2$.   
\end{proof}

Hence it is enough to prove that  $\tilde D_i$ is not included in a flat. We  prove that $\tilde D_i$ is actually maximally flat in the sense that it is not included in any open flat disk of $X$ centered at $A$. Indeed  assume that there is such a  disk $D_i'$. This gives a path of length $2\pi$ in the link $L_A$ at $A$ in $X$. The construction of $S_i$ and $S_0$ shows that $L_A$  contains the path $3251$ (see Figure \ref{bibifurq}, where we identified  $L_A$ with the simplicial sphere of radius 1 at $A$ so $L_A$ is included in the 1-skeleton of $X$ and inherits of its labelling). It is easily seen that there is no such a cycle in $L_A$.    

This shows that $\varphi_A$ is at least equal to $2^{2\mu_k-4}$ on $[k-\sqrt 3,k]$ and proves the lemma.
\end{proof}

\subsection{Proof of Theorem \ref{th11}, Item (a).} \label{meso-frises}
Let us first give some more details on the complex of friezes $\tilde V_{\bowtie}$.

Let $P$ be the complex constructed in Section 3 of \cite{toulouse}. Recall that $P$ is compact with 2 vertices $A_1$, $A_2$ whose universal cover  is an exotic triangle building $\Delta$ of order 2. By a theorem of Tits (\cite{Tits90}, see also Th\'eor\`eme 1 in \cite{toulouse}) there are exactly two isomorphism classes of spheres of radius 2 of this family of triangle buildings. They  correspond  to the 2-sphere of the building of $\PSL_3(K)$ in the two cases $K=\QI_2$ or $K=\FI_2((t))$. It is shown in \cite[Th\'eor\`eme 6]{toulouse} that the sphere of radius 2 at a vertex $A$ of $\Delta$ corresponds to $\QI^2$ if and only if $A$ projects to (say) $A_1$. 

 Let $S$ be the \emph{median section} which is described on Figure 25 of page 599 in \cite{toulouse}. This is a (metric) graph with 6 vertices and 9 edges, and  
  $P\backslash S$ is a disjoint union of two complex with one vertex and boundary isometric to $S$. Let $P_1$ be the closure of the complex which contains $A_1$ (i.e.\   which associated to $\QI^2$). 
  
  \begin{definition} The complex $V_{\bowtie}$ is defined to be the complex obtained by gluing together two copies of $P_1$  along $S$ via the identity map.
  \end{definition}
   
We denote by $O$ and $O'$ its two vertices corresponding to $A_1$ and number from 1 to 6 the vertices corresponding to $S\subset V_{\bowtie}$.

\begin{lemma}
The universal cover $X=\tilde V_{\bowtie}$ of complex $V_{\bowtie}$ is a CAT(0) space.
\end{lemma}

\begin{proof}
It is easily checked that the links at the vertices of $S$ all are isometric to a trivalent graph with two vertices, and 3 edges between these vertices of respective length $2\pi/3$, $4\pi/3$, and $4\pi/3$. Thus all cycles have length $\geq 2\pi$.
\end{proof}

\begin{figure}[htbp]
\centerline{\includegraphics[width=8cm]{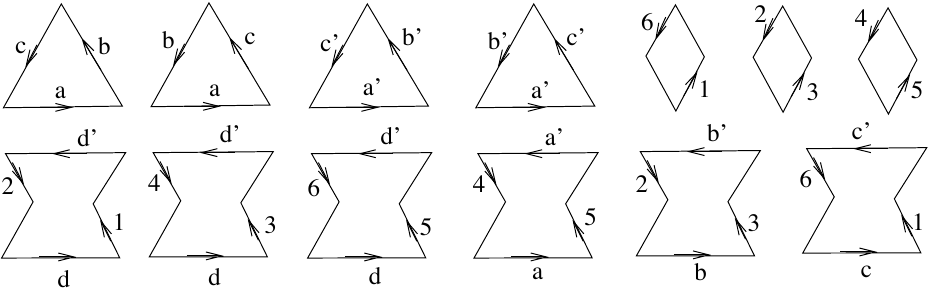}}
\caption{The group of friezes $\G_{\bowtie}$}\label{frise}
\end{figure}

The above surgery creates 3 shapes in $V_{\bowtie}$ (represented on Figure \ref{frise}): a triangle, a lozenge, and a bow tie (the latter we translate a `queue d'aronde' in french). They are unions of 4, 2, and 6  equilateral triangles respectively.   

The complex $V_{\bowtie}$ contains 4 triangles, 3 lozenges and 6 bow ties.  One checks that two of the triangles $t_1=abc$ and $t_2=acb$ are glued on $O$ while the two others $t_1'=a'b'c'$ and $t_2'=a'c'b'$ are glued on $O'$. Three bow ties are adjacent to a loop $d$ on $O$ and a loop $d'$ on $O'$. The three others are glued alongside  to  triangles as indicated on Figure \ref{frise}.

The following shows exponential mesoscopic rank of the group $\G_{\bowtie}$.

\begin{lemma}
Let $A$ be a vertex of $X$ which projects down to $O$ in $V_{\bowtie}$ and  $k\geq 3$ be an integer. Then on the interval $](k-1)\sqrt 3,k\sqrt 3]$  
of $\RI_+$ one has 
\[
\varphi_A\geq 2^{2k-4}.
\]
The function $\varphi_A$ has continuous support starting from ${\sqrt 3\over 2}=0.86...$, and exponential growth.
\end{lemma}

\newcommand{\St}{\mathrm{St}}

\begin{proof}
Let $T$ be a flat sector of angle $\pi/3$ in $X$ built out of $t_1$ and $t_2$, whose boundary is a union of two half-lines $d_1$ and $d_2$ intersecting at a point $B$. We can assume that $d_1$ and $d_2$ are labelled by $a$ and $b$ respectively. Let $[B,\infty[$ be the bisector of $T$. 

Consider the following three strips of height 1 in $X$: 
\begin{itemize}
\item  $\St_1$, corresponding to the fourth bow tie, with boundaries of the form $(a^\infty,a'^\infty)$, 
\item $\St_2$,  corresponding to the triangles $t_1,t_2$,  with boundaries of the form $(a^\infty, a^\infty)$, 
\item   $\St_3$, corresponding to the triangles   $t_1',t_2'$, with boundaries of the form $(a'^\infty,a'^\infty)$,
\end{itemize}
with the labelling of Figure \ref{frise} and the notations of Lemma \ref{lem-meso74}. Then, given a integer $k\geq 1$, a computation shows that there are precisely $\nu_k=2^{k-1}$ flat strips $S_1, \ldots, S_{\nu_k}$ of simplicial height $k$ on the line $d_1$ which start with $\St_1$.

The configuration is symmetric relatively $[B,\infty[$,  i.e.\ denoting  $\St_1^*, \St_2^*, \St_3^*$ the three strips defined analogously to $\St_1^*, \St_2^*, \St_3^*$, there 
 precisely $\nu_k$ flat strips $S_1^*, \ldots, S_{\nu_k}^*$ of simplicial height $k$ on the line $d_2$ which start with $\St_1^*$.

Let $A$ be a vertex of $]B,\infty[$, so that $|A-B|=k\sqrt 3$ for some integer $k\geq 1$ (note that by construction all vertices of $[B,\infty[$ project down to $O$ in $V_{\bowtie}$). We assume that $k>2$, so in particular 
\[
(k-1)\sqrt 3> k{\sqrt 3\over 2}.
\] 
As in the proof of Lemma \ref{lem-meso74}   for $i\in \{1,\ldots, \nu_k\}^2$  we let $\Pi_i$ be the subset of $X$ defined by $\Pi_i=T\cup S_{i_1}\cup S_{i_2}^*$. Then 
\begin{itemize}
\item the set $D_i$ of points of $\Pi_i$ at distance  $\leq k\sqrt 3$ from $A$ in $\Pi_i$ is a flat disk in $X$ whose boundary contains $B$,
\item  the disks $D_i$ are pairwise distinct when $i$ varies in  $\{1,\ldots,\nu_k\}^2$. 
\end{itemize}
Then for $r\in [0,k]$ the family of concentric disks $D_i^r$ of radius $r$ in $D_i$, for any fixed $r\in ](k-1){\sqrt 3},k{\sqrt 3}]$ contains at least $\nu_{k-2}^2$ distinct elements.

But for $j\in \{1,\ldots,\nu_k\}^2$ and $r\in ]k-{\sqrt 3},k]$ the disk $D_j^r$ is not included in a flat, because the disks $D_i$, $i\in \{1,\ldots, \nu_k\}^2$ are not (for a  similar, albeit more geometrical, reason to that of Lemma \ref{lem-meso74}: a singularity arises at the apex $B$ of $T$).

This shows that $\varphi_A$ is at least equal to $\nu_{k-2}^2=2^{2k-4}$ on $[(k-1)\sqrt 3,k\sqrt 3]$. The fact that $\varphi_A$ has continuous support as soon as $r>{\sqrt 3/2}$ is not hard to show. This proves the lemma.
\end{proof}

\begin{remark}
If in the above proof one had enumerated the strips $S_i, S_i^*$ up to simplicial isomorphic (i.e.\ up to `type of frieze') rather than up to equality in $X$, then  we would have found precisely $F_k$ such strips of simplicial height $k$,  where $(F_k)_{k\geq 0}$ is the Fibonacci sequence
\[
0,1,1,2,3,5,8,13,\ldots
\]
\end{remark}








\end{document}